\documentclass[11pt]{article}

\usepackage{fullpage}
\usepackage{amsmath, amsthm, amssymb}
\usepackage{bbm}
\usepackage{authblk}
\usepackage{mathtools}
\usepackage{graphicx}
\usepackage{xcolor}
\usepackage{enumerate}
\usepackage{booktabs}
\usepackage{tabularx}
\usepackage{natbib}
\usepackage{subfig}
\usepackage{multirow}
\usepackage{float}
\usepackage{tikz}
\usepackage{longtable}
\usepackage{indentfirst}
\usepackage{setspace}
\usepackage{xurl}
\usepackage[inline]{enumitem}
\usepackage{arydshln}
\usepackage{rotating}
\RequirePackage[colorlinks = true, linkcolor = blue, citecolor=blue,urlcolor=blue]{hyperref}

\newtheorem{lemma}{Lemma}
\newtheorem{cor}{Corollary}
\newtheorem{theorem}{Theorem}

\newtheorem{prop}{Proposition}

\newcommand{\dd}{\,\mathrm{d}}
\DeclareMathOperator*{\argmin}{argmin}

\DeclareMathOperator{\RCV}{RCV}

\DeclareMathOperator{\sign}{sign}
\DeclareMathOperator{\BL}{BL}
\DeclareMathOperator{\Span}{Span}
\DeclareMathOperator{\SEE}{SEE}

\definecolor{c1}{rgb}{0,  0, 1}
\definecolor{c2}{rgb}{0,  0, 0}
\definecolor{c3}{rgb}{128,  0, 128}

\setlist[itemize]{noitemsep}
\setlist[enumerate]{noitemsep}

\begin{document}
%%%%%%%%%%%%%%%%%%%%%%%%%%%%%%%%%%%%%%%%%%%%%%%%%%%%%%%%%%%%%%%%%%%%%%%%%%%%%%

\title{Robust multidimensional location estimation from discretely sampled functional data}
%\author{Ioannis Kalogridis}
%\affil{School of Mathematics and Statistics, University of Glasgow}
\author[1]{Ioannis Kalogridis}
\author[2]{Stanislav Nagy}
\author[2]{Michele Cavazzutti}
\affil[1]{School of Mathematics and Statistics, University of Glasgow}
\affil[2]{Faculty of Mathematics and Physics, Charles University}
\date{\today}
\maketitle

\begin{abstract}

We introduce and study the first class of robust location estimators for discretely sampled functional data on $\mathbbm{R}^d$. The proposed method is based on M-estimation with repeated measurements and a Sobolev-type penalty, and is suitable for both commonly and independently observed trajectories, even in the presence of heavy-tailed random fields or measurement error. Under mild assumptions, we show that the estimators are well-defined elements of an appropriate infinite-dimensional function space, they are measurable, and admit unique finite-dimensional representations that make them computationally feasible in practice. Moreover, the estimators are equivariant under rigid motions of the domain and are domain-agnostic, in the sense that they do not depend on the domain's shape or boundary. We establish sharp convergence rates in both fixed and random design settings and, in particular, show that in both cases the optimal parametric rate $n^{-1}$ is attainable even when the data are incomplete or irregularly observed. An extensive simulation study and a real-data example demonstrate the practical usefulness of the proposed methods.

\end{abstract}



\section{Introduction}
\label{sec:1}
Functional data analysis (FDA) is the modern branch of Statistics providing tools to summarize, model, and infer properties of random processes $X: (\Omega, \mathcal{A}, \mathbb{P}) \to L_2(\mathcal{I})$ from a random sample $X_1, \ldots, X_n$ at our disposal. Here,  $L_2(\mathcal{I})$ denotes the space of square integrable functions on $\mathcal{I} \subset \mathbbm{R}^d$ with its usual Borel $\sigma$-algebra. Although by now FDA has become a broad field of research with many specialized sub-paths, the estimation of centrality remains a fundamental problem. Of particular interest is the mean function of $X$ defined as the Bochner integral $\mu_0 = \mathbb{E}[X]$. The integral exists whenever $\mathbb{E}[\|X\|_{L_2(\mathcal{I})}] < \infty$ in which case $\mu_0$ provides a concise summary of the overall trend of $X$ over $\mathcal{I}$. But beyond its intrinsic interpretative value, the mean function also plays a central role in key downstream FDA procedures such as covariance estimation, functional principal component analysis, and functional regression. The interested reader is referred to the classical monographs of \citet{Ramsay_Silverman2005}, \citet{Ferraty_Vieu2006}, and \citet{Horvath_Kokoszka2012} for detailed discussions and illustrative examples. 

It is well known that for jointly measurable processes $X$, the mean function $\mu_0$ can  be defined pointwise as $\mu_0(\mathbf{t}) = \mathbb{E}[X(\mathbf{t})], \mathbf{t} \in \mathcal{I}$ \citep[see, e.g.][Chapter 7]{Hsing:2015}. This implies that $\mu_0$ can be estimated in a pointwise manner; in particular, if the random functions $X_i$ were fully observed on $\mathcal{I}$, then $\mu_0$ could be readily estimated with the sample mean $\bar{X}(\mathbf{t})$. In practice, however, complete observations are almost never feasible and instead only discretely sampled and noisy observations can be recorded:
\begin{align}
\label{eq:dfda}
Y_{ij} = X_i(\mathbf{T}_{ij}) + \zeta_i(\mathbf{T}_{ij}), \quad (j=1,\ldots,m_i;\ i=1,\ldots,n),
\end{align}
where ${\mathbf{T}_{ij}} \in\mathcal{I}$ are sampling vectors and $\{\zeta_i\}_{i=1}^n \in L_2(\mathcal{I})$ are noise processes on the same probability space $(\Omega,\mathcal{A},\mathbb{P})$. The noise processes $\zeta_i$ account for uncontrolled factors such as measurement error or environmental variability. In the presence of discretely sampled functional data, classical estimators, such as the sample mean and its M-type generalizations \citep{Sinova:2018}, may perform poorly for two important reasons: (i) there may exist subregions of $\mathcal{I}$ with few or even no design points making these estimators highly variable, and (ii) these estimators incorporate noise from the $\zeta_i(\mathbf{T}_{ij})$; with high levels of noise these estimators yet again become highly variable with extreme local fluctuations.

Location estimation from discretely sampled functional data has been well-studied for data observed on a bounded one-dimensional domain $\mathcal{I} \subset \mathbbm{R}$. We mention, in particular, the landmark contributions of \citet{Cai:2011}, who derived the optimal rates of convergence and showed that these are attained by smoothing spline type estimators; \citet{Zhang:2016}, who investigated general weighting schemes and showed that the optimal rates are also attained by kernel estimators; and \citet{Xiao:2020}, who showed that the optimal rates may also be attained by lower rank penalized spline estimators. More recently, the results of \citet{Cai:2011} were generalized to a wide variety of location estimators, not necessarily for the mean function, by \citet{K:2023}, who showed that generalized (M-type) smoothing spline estimators for data observed on a bounded one-dimensional $\mathcal{I}$ likewise attain the optimal rates of convergence under standard regularity conditions.

The main drawback of existing estimators, spline- and kernel based alike, is their reliance on one-dimensional domains. This is restrictive, as in many contemporary applications, such as oceanography and meteorology, the data are discretely observed on multidimensional domains $\mathcal{I} \subset \mathbbm{R}^d$ so that existing one-dimensional location estimation methods are no longer applicable. Location estimation methods for multidimensional discretely observed data remain scarce, and to the best of our knowledge,  the only relevant contribution to-date is due to \citet{Berger_etal2025}, who proposed local linear kernel estimators with asymptotic guaranties for functional data observed on a lattice.  Nevertheless, their approach has four notable limitations. First, it requires common sampling locations $\mathbf{T}_{ij}$ across trajectories, which excludes longitudinal data. Second, its asymptotic theory assumes a product domain $\mathcal{I} = \prod_{k=1}^d \mathcal{I}_k$ with $\mathcal{I}_k \subset \mathbbm{R}$, thereby ruling out irregular observation domains that frequently arise in environmental applications. Third, it is computationally demanding, as the kernel estimator of those authors must be refitted at every point $\mathbf{t}\in\mathcal{I}$ of interest, making dense coverage of the domain infeasible. Finally, it is based on least squares and therefore lacks robustness: functional data are frequently contaminated by outlying trajectories or heavy-tailed measurement errors, and the sample mean, as well as least-squares-based smoothers, can be arbitrarily distorted by even a few such observations.

To address these shortcomings, we take the spline perspective in this paper and propose a novel, flexible and computationally feasible family of location estimators accommodating (i) bounded multidimensional domains $\mathcal{I}\subset \mathbbm{R}^d$ of arbitrary shape, (ii) fixed or random sampling points $\mathbf{T}_{ij}$, (iii) subject-specific sample sizes $m_i$, and (iv) robustness to contamination, obtained through a general convex loss $\rho$, for example the absolute or Huber loss in place of the quadratic one. Our contributions are as follows:

\begin{enumerate}
    \item We show that the proposed estimators constitute natural generalizations of the one-dimensional spline estimators studied in the functional data literature; establish their formal existence under weak conditions and characterize them as finite-dimensional thin-plate splines \citep[see, e.g.][]{Green:1994, Wend:2005}. We further formalise their equivariance under rigid motions (rotations, reflections and translations) of the domain, which makes these estimators particularly well suited for spatial smoothing.
    \item We study the theoretical properties of the proposed estimators and derive detailed rates of convergence with a useful decomposition of the $L_2(\mathcal{I})$ error under both fixed (common) and random (subject-specific) design settings. Our decomposition points exactly to the sources of error that these estimators incur and to the phase transition phenomenon first observed by \citet{Cai:2011} for one-dimensional domains. Our results hold not only for the popular quadratic loss but very broadly for a wide variety of convex loss functions satisfying our weak regularity conditions.
    \item Using the established finite-dimensional representation, we develop a computationally efficient algorithm implemented in the dedicated \texttt{R} package \texttt{RobustSpline}, featuring a \texttt{C++} back-end built upon the modern OSQP algorithm introduced by \citet{stellato2020osqp}, which, to the best of our knowledge, has not been previously used in this context.
\end{enumerate}
In particular, our theoretical results reveal that if $m$, the number of discretization points in the fixed design or the harmonic mean of the $m_i$ in the random design, exceeds $n^{d/(2r)}$, where $r$ is the smoothness index of $\mu_0$, then the estimators attain the parametric rate $n^{-1}$. Otherwise, the convergence rate is governed by the discretization density or jointly by $(m,n)$, see Section~\ref{sec:3}.  Our theoretical framework, combining techniques from functional analysis, reproducing kernel Hilbert spaces of conditionally positive definite kernels, and empirical process theory, provides new insights into the interplay between discretization, sample size, dimension and smoothness and may be of independent interest.

\paragraph{Notation.} Throughout, $d \in \mathbbm{N}$ and $D \subset \mathbbm{R}^d$ is a bounded, open domain with Lipschitz boundary. For a multi-index $\alpha \in \mathbbm{N}_0^d$ we write $|\alpha| = \sum_{j=1}^d \alpha_j$, $\alpha! = \prod_{j=1}^d \alpha_j!$ and $\partial^{\alpha} = \partial_{x_1}^{\alpha_1}\cdots\partial_{x_d}^{\alpha_d}$, with derivatives understood in the weak sense on $D$ or $\mathbbm{R}^d$ (for continuously differentiable functions they coincide with classical derivatives). For $r \in \mathbbm{N}_0$, $\mathcal{C}^{r}(D)$ denotes the space of $r$ times continuously differentiable real-valued functions on $D$ and $\mathcal{P}_r(D)$ denotes the space of $r$th order polynomials on $D$. When $D = \mathbbm{R}^d$, we abbreviate $\mathcal{P}_r(\mathbbm{R}^d)$ with $\mathcal{P}_r$. The sup-norm  of a bounded function $f$ on $D$ is  $\|f\|_{\infty,D} = \sup_{x \in D}|f(x)|$. For $A \in \{D, \mathbbm{R}^d\}$, $L_2(A)$ is the space of square-integrable real-valued functions on $A$ with inner product $\langle f,g\rangle_{L_2(A)} = \int_A f(x)g(x)\,\mathrm{d}x$. The $r$th order Hilbert--Sobolev space on $D$ is
\begin{align*}
  \mathcal{H}^{r}(D) = \bigl\{ f \in L_2(D) : \partial^{\alpha} f \in L_2(D) \text{ for all } |\alpha| \leq r \bigr\},
\end{align*}
a Hilbert space with inner product $\langle f,g \rangle_{\mathcal{H}^{r}(D)} = \langle f,g\rangle_{L_2(D)} + \sum_{|\alpha| = r} \frac{r!}{\alpha!}\,\langle \partial^{\alpha} f, \partial^{\alpha} g\rangle_{L_2(D)}$. For $\ell \in \mathbbm{N}$ with $\ell > d/2$, the Beppo Levi space of order $\ell$ \citep[Section 10.5]{Wend:2005} is
\begin{align*}
  \BL_{\ell}(\mathbbm{R}^d) = \bigl\{ f \in \mathcal{C}(\mathbbm{R}^d) : \partial^{\alpha} f \in L_2(\mathbbm{R}^d) \text{ for all } |\alpha| = \ell \bigr\}
\end{align*}
with the semi-inner product $( f,g)_{\BL_{\ell}(\mathbbm{R}^d)} = \sum_{|\alpha| = \ell} \frac{\ell!}{\alpha!}\,\langle \partial^{\alpha} f, \partial^{\alpha} g\rangle_{L_2(\mathbbm{R}^d)}$ and seminorm $|f|_{\BL_{\ell}(\mathbbm{R}^d)} = ( f,f)_{\BL_{\ell}(\mathbbm{R}^d)}^{1/2}$, whose null space is the space of polynomials of total degree at most $\ell-1$. Two bilinear forms on $\BL_\ell(\mathbbm{R}^d)$ that we will be using are
\begin{align}
     \langle f, g \rangle_n & := \sum_{i=1}^n \sum_{j=1}^{m_i} f(\mathbf{T}_{ij}) g(\mathbf{T}_{ij}) + ( f, g )_{\BL_\ell(\mathbbm{R}^d)}. \\ 
    \langle f, g \rangle_{\BL_r(\mathbbm{R}^d), \lambda} &: = \langle f, g \rangle_{L_2(\mathcal{I})} + \lambda ( f, g )_{\BL_\ell(\mathbbm{R}^d)}, \label{eq:equivnorm}
\end{align}
where $\lambda>0$ is the penalty parameter. The induced norms are denoted with $\|\cdot\|_n$ and $\|\cdot\|_{\BL_\ell(\mathbbm{R}^d), \lambda}$, respectively. As shown in Section~\ref{sec:2}, both $\langle \cdot, \cdot \rangle_n$ and $\langle \cdot, \cdot \rangle_{\BL_\ell(\mathbbm{R}^d), \lambda}$ are inner products on $\BL_{\ell}(\mathbbm{R}^d)$ and equipped with either $\BL_{\ell}(\mathbbm{R}^d)$ is a reproducing kernel Hilbert space (RKHS). For a function $f$ defined on $A$ and $B \subseteq A$, $f|_{B}$ denotes the restriction of $f$ to $B$; in particular, if $f \in \BL_{\ell}(\mathbbm{R}^d)$, then $f|_{D} \in \mathcal{H}^{\ell}(D)$.

\section{Thin-Plate Splines for Functional Location Estimation}
\label{sec:2}
\subsection{The Proposed Family of Estimators}

To motivate the form of our estimator, it is instructive to revisit the one-dimensional case. Here we assume that $\mathcal{I} \subset \mathbbm{R}$ is open, bounded and connected, say $\mathcal{I} = (0,1)$ without loss of generality. Effective estimators for $\mu_0(t) := \mathbb{E}[X(t)]$ from smoothing spline principles have been proposed by \citet{Cai:2011} and \citet{K:2023}. These are obtained by solving
\begin{align}
\label{eq:1dest}
\widehat{\mu}_n = \argmin_{f \in \mathcal{H}^r((0,1))} \left[ \frac{1}{n} \sum_{i=1}^n \frac{1}{m_i} \sum_{j=1}^{m_i}   \rho \left(Y_{ij} - f(T_{ij}) \right) + \lambda \int_{0}^1 |f^{(r)}(t)|^2 \dd t \right],
\end{align}
for a convex loss function $\rho:\mathbbm{R} \to \mathbbm{R}_{+}$ and a regularization parameter $\lambda>0$. The minimiser need not be unique unless $\rho$ is strictly convex. However, it can be shown \citep[see, e.g.][]{K:2023} that a minimiser $\widehat{\mu}_n \in \mathcal{H}^r((0,1))$ exists and is a piecewise polynomial with the following properties: (i) $\widehat{\mu}_n$ is of order $2r$ (degree $2r-1$) between consecutive $\{T_{ij}\}_{i=1, j = 1}^{n, m_i}$, (ii) $\widehat{\mu}_n$ is of order $r$ on $(-\infty, \min_{i,j} T_{ij}] \cup [\max_{i,j} T_{ij}, \infty)$ and (iii) $\widehat{\mu}_n \in \mathcal{C}^{2r-2}(\mathbbm{R})$. Consider now extending the integration in the penalty from $(0,1)$ to the whole of $\mathbbm{R}$ and define the new objective function $\widetilde{L}_n(f)$
\begin{align*}
\widetilde{L}_n(f) = \frac{1}{n} \sum_{i=1}^n \frac{1}{m_i} \sum_{j=1}^{m_i}   \rho \left(Y_{ij} - f(T_{ij}) \right) + \lambda \int_{\mathbbm{R}} |f^{(r)}(t)|^2 \dd t.
\end{align*}
Since $2r>1$ for every $r \in \mathbbm{N}$, our definitions imply that $\mathcal{H}^r(\mathbbm{R}) \subset \BL_r(\mathbbm{R})$. 
Although the objective functions are different, perhaps surprisingly, $\widehat{\mu}_n$ is also a minimiser of $\widetilde{L}_n$ over $\BL_r(\mathbb{R})$, as we now show. This observation is the main motivation for the estimators herein.
\begin{lemma}
\label{lem:equiv}
Suppose that $\rho: \mathbbm{R} \to \mathbbm{R}_{+}$ is convex and coercive, i.e., satisfies $\lim_{|x|\to \infty} \rho(x) = \infty$. Moreover, suppose that $\# \{ \text{distinct} \ T_{ij}\} \geq r$. Then $\widehat{\mu}_n$ minimises $\widetilde{L}_n$ over $\BL_{r}(\mathbbm{R})$.
\end{lemma}
\begin{proof}
     Denote the objective function leading to $\widehat{\mu}_n$ by $L_n$, viz,
\begin{align*}
    L_n(f) = \frac{1}{n} \sum_{i=1}^n \frac{1}{m_i} \sum_{j=1}^{m_i}   \rho \left(Y_{ij} - f(T_{ij}) \right) + \lambda \int_{0}^{1} |f^{(r)}(t)|^2 \dd t.
\end{align*}
Since $\# \{ \text{distinct} \ T_{ij}\} \geq r$ and $\rho$ is convex and coercive, its minimiser $\widehat{\mu}_n$ is well-defined and fulfills properties (i)--(iii) \citep[see, e.g. Theorem~1]{K:2022a}. But, by property (ii), the $r$th derivative of $\widehat{\mu}_n$ vanishes outside $(0,1)$ and so $\|\widehat{\mu}_n^{(r)}\|_{L_2(\mathbbm{R})} = \|\widehat{\mu}_n^{(r)}\|_{L_2((0,1))}$. In particular, $\widehat{\mu}_n \in \BL_r(\mathbbm{R})$, which implies that for every $f \in \BL_r(\mathbbm{R})$, $f|_{(0,1)} \in \mathcal{H}^r((0,1))$. Therefore, by the definition of $\widehat{\mu}_n$,
\begin{align*}
\widetilde{L}_n(\widehat{\mu}_n) = L_n(\widehat{\mu}_n) \leq L_n(f) \leq \widetilde{L}_n(f),
\end{align*}
for every $f \in \BL_r(\mathbbm{R})$, where the last inequality follows by extending the domain of integration to $\mathbbm{R}$. The result of the lemma follows.
\end{proof}

The condition $\# \{ \text{distinct} \ T_{ij}\} \geq r$ is needed to ensure polynomial unisolvency in $\mathbbm{R}$ and is standard in smoothing spline estimation. Lemma~\ref{lem:equiv} shows that $\widehat{\mu}_n$ solves a more general variational problem than initially assumed and unlocks the generalization to higher dimensions. We shall assume henceforth that $\mathcal{I} \subset \mathbbm{R}^d$ is open bounded and connected, $2r>d$ and $X: (\Omega, \mathcal{A}, \mathbb{P}) \times \mathcal{I} \to \mathbbm{R}$ is a random field. We seek to estimate $\mu_0(\mathbf{t}) := \mathbb{E}[X(\mathbf{t})]$ from the discretely sampled functional data given in \eqref{eq:dfda}. We shall also assume that $\mu_0 \in \BL_r(\mathbbm{R}^d)$, the multidimensional Beppo Levi space. We propose to estimate $\mu_0$ with 
\begin{align}
\label{eq:est}
\widehat{\mu}_n = \argmin_{f \in \BL_r(\mathbbm{R}^d)} \left[ \frac{1}{n} \sum_{i=1}^n \frac{1}{m_i} \sum_{j=1}^{m_i} \rho \left(Y_{ij} - f(\mathbf{T}_{ij}) \right) + \lambda |f|_{\BL_r(\mathbbm{R}^d)}^2 \right],
\end{align}
where $\rho$ is a convex loss function and $\lambda>0$ is a regularization parameter controlling the trade-off between fidelity to the observed values $Y_{ij}$ and smoothness; small values of $\lambda$ lead to estimates close to the $Y_{ij}$ while large values $\lambda$ push the estimates to the null space of $|\cdot|_{\BL_r(\mathbbm{R}^d)}$, which consists of all polynomials of total order $r$ on $\mathbbm{R}^d$, see \citet[Lemma 10.38]{Wend:2005} for a formal proof.  By Lemma~\ref{lem:equiv}, $\widehat{\mu}_n$ is a direct generalization of one-dimensional spline estimators. The intriguing penalty functional $|\cdot|_{\BL_r(\mathbbm{R}^d)}^2$ was introduced by \citet{Duchon:1977} for interpolation in higher dimensions; see \citet[Chapter~2.4]{Wahba:1990} or \citet[Chapter~7]{Green:1994} for discussions from a statistical point of view.

Before we turn to the theoretical properties of $\widehat{\mu}_n$, we make two observations. First, it may seem awkward that $X$ is defined on a bounded set $\mathcal{I} \subset \mathbbm{R}^d$, yet we require that $\mu_0$ is the restriction of an element of $\BL_{r}(\mathbbm{R}^d)$. That no contradiction arises follows from the Sobolev extension theorem \citep[see, e.g.,][Theorem~5.24]{Adams:2003}: since $\mathcal{I}$ has a Lipschitz boundary and $2r>d$ there exists a bounded extension operator $E: \mathcal{H}^r(\mathcal{I}) \to \mathcal{H}^r(\mathbbm{R}^d)$, and $\mathcal{H}^r(\mathbbm{R}^d) \subset \BL_r(\mathbbm{R}^d)$. Thus, starting from any $\mu_0 \in \mathcal{H}^r(\mathcal{I})$, the extension $E(\mu_0)$ belongs to $\BL_r(\mathbbm{R}^d)$ and satisfies $E(\mu_0)|_{\mathcal{I}} = \mu_0$. Secondly, although we focus in this work on the estimation of the mean function $\mu_0$, our framework is much more general and allows for the estimation of other population quantities, including the pointwise median, quantiles and expectiles. These are obtained by choosing appropriate loss functions $\rho$, for example, $\rho(x) = |x|$ for the median, $\rho_{\tau}(x) = (\tau - \mathbbm{1}\{x<0\})x$ for the $\tau$-quantile, and the asymmetric expectile loss $\rho_{\alpha}(x) = |\alpha-\mathbbm{1}\{x<0\}|x^2$ for the $\alpha$-expectile, with $\alpha \in (0,1)$. These estimators are also covered by our existence, measurability and asymptotic results in the following sections.

\subsection{Existence, Characterization and Spatial Equivariance}

We next establish the formal existence of our thin-plate spline estimators. The main technical hurdle is that $|\cdot|_{\BL_r(\mathbbm{R}^d)}$ is only a semi-norm on $\BL_r(\mathbbm{R}^d)$ and so, equipped with it, $\BL_r(\mathbbm{R}^d)$ is not a Hilbert space. But following the development in \citet[Theorem~10.20]{Wend:2005}, we can equip $\BL_r(\mathbbm{R}^d)$ with the bilinear form $\langle \cdot, \cdot \rangle_n$ defined at the end of Section~\ref{sec:1}. Denote $\{\mathbf{T}_{ij}\}_{i = 1, j=1}^{n, m_i}$ with $\mathcal{T}$. To ensure that $\BL_r(\mathbbm{R}^d)$ equipped with $\langle \cdot, \cdot \rangle_n$ is a Hilbert space, $\mathcal{T}$ needs to be unisolvent for the space of polynomials of order $r$ on $\mathbbm{R}^d$, $\mathcal{P}_r$. Specifically, if for any $p \in \mathcal{P}_r$ it needs to hold that $p(\mathbf{T}_{ij}) = 0,  j=1, \ldots, m_i; \ i=1, \ldots, n$ implies that $p = 0$, i.e., $p$ is the zero polynomial. Under polynomial unisolvency, Proposition~\ref{prop:1} shows that $\BL_r(\mathbbm{R}^d)$ equipped with $\langle \cdot, \cdot \rangle_n$ is not only a Hilbert space, but also a reproducing kernel Hilbert space (RKHS) wherein evaluation functionals $\delta_{\mathbf{x}}(f) = f(\mathbf{x})$ are continuous. 
\begin{prop}
\label{prop:1}
Suppose that $2r>d$ and that $\mathcal{T}$ is $\mathcal{P}_r$-unisolvent. Then,  $\BL_r(\mathbbm{R}^d)$ equipped with $\langle \cdot,\cdot \rangle_{n}$ is an RKHS.
\end{prop}
The proof is given in the Appendix. Polynomial unisolvency is a widely satisfied condition.  For $d=1$, it is satisfied if $\mathcal{T}$ contains at least $r$ distinct points, which is a classical condition for the existence of smoothing spline estimators \citep[see, e.g.][Theorem~2.4]{Green:1994}. For $d>1$, the condition is fulfilled if $\mathcal{T}$ contains at least $M := \dim \mathcal{P}_r$ elements satisfying no nontrivial polynomial relationship. This will be case with probability one, for example, if $\mathcal{T}$ contains $M$ independent and identically distributed (i.i.d.) sampling vectors from an absolutely continuous distribution on $\mathcal{I}$. 

Prior to the statement of our existence result in  Theorem~\ref{thm:1} below, we recall that a thin-plate spline of order $r$ with \textit{knots} at $\mathcal{T} := \{\mathbf{T}_{ij}\}_{i = 1, j = 1}^{n, m_i}$ is any function $g: \mathbbm{R}^d \to \mathbbm{R}$ of the form
\begin{align}
\label{eq:TP}
g\left(\mathbf{t}\right) = \sum_{i=1}^n\sum_{j=1}^{m_i} \gamma_{ij} \eta_{dr} \left(\left\|\mathbf{t}-\mathbf{T}_{ij} \right\|_{\mathbbm{R}^d} \right) + \sum_{k=1}^M \delta_k \phi_k\left(\mathbf{t}\right),
\end{align} 
with $\eta_{dr}: \mathbbm{R}_{+} \to \mathbbm{R}$ given by
\begin{align*}
\eta_{dr}(x) = \begin{cases} \frac{(-1)^{r+1+d/2}}{2^{2r-1} \pi^{d/2} (r-1)!(r-d/2)!} x^{2r-d} \log(x) & d \ \text{even}, \\ 
\frac{\Gamma(d/2-r)}{2^{2r} \pi^{d/2} (r-1)!}x^{2r-d} & d \ \text{odd}, \end{cases}
\end{align*}
where $\Gamma(\cdot)$ denotes the Gamma function. The functions $\phi_{1}, \ldots, \phi_M$ are any basis for the space $\mathcal{P}_r$, which has dimension $M = \binom{r+d-1}{d}$. A thin-plate spline $g$ is called a natural thin-plate spline, if, in addition to \eqref{eq:TP}, $|g|_{\BL_r(\mathbbm{R}^d)}<\infty$. In other words, natural thin-plate splines are precisely the thin-plate splines that are elements of $\BL_r(\mathbbm{R}^d)$. It can be shown that the kernel $\Phi_{dr} : = \eta_{dr}(\|\cdot\|_{\mathbbm{R}^d})$ is conditionally positive definite of order $r$ and its induced native space coincides with $\BL_r(\mathbbm{R}^d)$; the interested reader is referred to \citet{Wend:2005} for a thorough development.

For a rigorous theoretical development it is important to know whether, and under which conditions, \eqref{eq:est} is well-defined as an element of $\BL_r(\mathbbm{R}^d)$ and measurable. We take a more general perspective than necessary for this work in Theorem~\ref{thm:1} below and establish the existence of minimisers in $\BL_r(\mathbbm{R}^d)$ for objective functions of the form
\begin{align}
\label{eq:objf}
L_n(f) = \frac{1}{n} \sum_{i=1}^n \frac{1}{m_i} \sum_{j=1}^{m_i} G(f(\mathbf{T}_{ij}),Y_{ij}) + \lambda |f|_{\BL_r(\mathbbm{R}^d)}^2.
\end{align}
Here, $G:\mathbbm{R} \times \mathbbm{R} \to \mathbbm{R}$ is \textit{any} loss function. In the functional location setting, $G(x,y) = \rho(y-x)$ but this formulation allows for more general estimators. For example, we can have $G(x,y) = b(x)-xy$, leading to nonparametric generalized linear models with repeated measurements. 
\begin{theorem}
\label{thm:1}
Suppose that for every $y \in \mathbbm{R}$ the function $x \mapsto G(x,y)$ is continuous, bounded from below, and satisfies $\lim_{|x| \to \infty} G(x,y) = \infty$, and that for every $x \in \mathbbm{R}$ the function $y \mapsto G(x,y)$ is measurable. Suppose further that $2r>d$, $\lambda>0$ and that $\mathcal{T}$ is $\mathcal{P}_r$-unisolvent. Then, the following hold.
\begin{enumerate}
\item[(a)] There exists a minimiser of $L_n$ in $\BL_{r}(\mathbbm{R}^d)$, and every minimiser is a natural thin-plate spline of order $r$ with knots at the distinct points of $\mathcal{T}$.
\item[(b)] If $x \mapsto G(x,y)$ is strictly convex for every $y \in \mathbbm{R}$, the minimiser is unique.
\item[(c)] If the $\mathbf{T}_{ij}$ and $Y_{ij}$ are random variables on $(\Omega,\mathcal{A},\mathbb{P})$, so that the above assumptions hold for every $\omega \in \Omega$, a minimiser $\widehat{\mu}_n$ can be chosen such that $(\omega,\mathbf{t}) \mapsto \widehat{\mu}_n(\omega,\mathbf{t})$ is jointly measurable.
\end{enumerate}
\end{theorem}

The characterization of the estimators in (a) forms the backbone of the computational algorithm in Section~\ref{sec:4}. The conditions of Theorem~\ref{thm:1} are notably weaker than those imposed in the literature. For example, \citet[Theorem~2.9]{Gu:2013} and \citet[Theorem~3.2]{Cox:1985} require the convexity of $x \mapsto G(x,y)$ as well as the uniqueness of the minimisers in the null space of $|\cdot|_{\BL_r(\mathbbm{R}^d)}$, $\mathcal{P}_r$. By contrast, Theorem~\ref{thm:1} allows for a broader family of loss functions provided that they are continuous and coercive. An example of a non-convex loss fulfilling these properties is the negative Cauchy log-likelihood $G(x,y) = \log(1+(y-x)^2)$. But even for convex $x \mapsto G(x,y)$, the uniqueness of the minimiser in the null space of $|\cdot|_{\BL_r(\mathbbm{R}^d)}$ is not guaranteed:  the finite-dimensional estimator with the absolute loss $G(x,y) = |y - x|$ is well-known to exist but need not be unique \citep[see, e.g.][p.~31]{Koenker:2005}. Thus, Theorem~\ref{thm:1} allows us to consider a whole new family of thin-plate spline estimators that, to the best of our knowledge, had not been covered by statistical theory. Theorem~\ref{thm:1} also shows that these estimators can be chosen measurable, which is by itself a novel result.

The existence and measurability of the minimiser and the properties of $\BL_r(\mathbbm{R}^d)$ now jointly imply a special form of equivariance. Let $\mathbf{P} \in \mathbbm{R}^{d \times d}$ denote an orthogonal matrix, i.e., $\mathbf{P} \mathbf{P}^{\top} = \mathbf{P}^{\top} \mathbf{P} = \mathbf{I}_d$ and $\mathbf{b} \in \mathbbm{R}^d$. Define a rigid motion $\phi: \mathbbm{R}^d \to \mathbbm{R}^d: \mathbf{x} \mapsto \mathbf{P}\mathbf{x}+\mathbf{b}$ and consider minimising the objective function
\begin{align*}
L_n^{\phi}(f) =  \frac{1}{n} \sum_{i=1}^n \frac{1}{m_i} \sum_{j=1}^{m_i} G(f \circ \phi (\mathbf{T}_{ij}),Y_{ij}) + \lambda |f|_{\BL_r(\mathbbm{R}^d)}^2. 
\end{align*} 
That is, $L_n^{\phi}(f)$ is the objective function after rotating or reflecting and translating the original sampling vectors $\mathbf{T}_{ij}$. Proposition~\ref{prop:inv} establishes that $\widehat{\mu}_n$ minimises $L_n$ in \eqref{eq:objf} if, and only if, $\widetilde{\mu}_n : = \widehat{\mu}_n \circ \phi^{-1}$ minimises $L_n^{\phi}$.

\begin{prop}
\label{prop:inv}
Suppose $2r>d$ and $\mathcal{T}$ is $\mathcal{P}_r$-unisolvent. Let $\phi(\mathbf{x}) = \mathbf{P}\mathbf{x} + \mathbf{b}$ where $\mathbf{P} \in \mathbbm{R}^{d \times d}$ satisfies $\mathbf{P} \mathbf{P}^{\top} = \mathbf{P}^{\top} \mathbf{P}  = \mathbf{I}_d$ and $\mathbf{b} \in \mathbbm{R}^d$. Then, $\widetilde{\mu}_n = \widehat{\mu}_n \circ \phi^{-1}$ minimises $L_n^{\phi}$ if, and only if, $\widehat{\mu}_n$ minimises $L_n$.
\end{prop}

\begin{proof}
It is a standard exercise \citep[see, e.g., Exercise~4.10 in][]{Gu:2013} to show that if $f \in \BL_r(\mathbbm{R}^d)$, then $f \circ \phi \in \BL_r(\mathbbm{R}^d)$ and $f \circ \phi^{-1} \in \BL_r(\mathbbm{R}^d)$, as well as
\begin{align}
\label{eq:TPinv}
|f\circ \phi^{-1}|_{\BL_r(\mathbbm{R}^d)}^2 = |f|^2_{\BL_r(\mathbbm{R}^d)}= |f \circ \phi|_{\BL_r(\mathbbm{R}^d)}^2.
\end{align}
In particular, $f \mapsto f\circ\phi$ is a bijection of $\BL_r(\mathbbm{R}^d)$ onto itself. Moreover, by \eqref{eq:TPinv}, for every $f \in \BL_r(\mathbbm{R}^d)$,
\begin{align*}
L_n^{\phi}(f\circ\phi^{-1}) = \frac{1}{n} \sum_{i=1}^n \frac{1}{m_i} \sum_{j=1}^{m_i} G\bigl(f(\mathbf{T}_{ij}), Y_{ij}\bigr) + \lambda |f\circ\phi^{-1}|_{\BL_r(\mathbbm{R}^d)}^2 = L_n(f).
\end{align*}
Now suppose $\widehat{\mu}_n$ minimises $L_n$ and let $g \in \BL_r(\mathbbm{R}^d)$ be arbitrary. Since $g\circ\phi \in \BL_r(\mathbbm{R}^d)$,
\begin{align*}
L_n^{\phi}(\widehat{\mu}_n\circ\phi^{-1}) = L_n(\widehat{\mu}_n) \leq L_n(g\circ\phi) = L_n^{\phi}(g).
\end{align*}
Hence $\widehat{\mu}_n\circ\phi^{-1}$ minimises $L_n^{\phi}$. Conversely, if $\widehat{\mu}_n\circ\phi^{-1}$ minimises $L_n^{\phi}$ and $f \in \BL_r(\mathbbm{R}^d)$, then $L_n(\widehat{\mu}_n) = L_n^{\phi}(\widehat{\mu}_n\circ\phi^{-1}) \leq L_n^{\phi}(f\circ\phi^{-1}) = L_n(f)$, so $\widehat{\mu}_n$ minimises $L_n$.
\end{proof}

An immediate consequence of Proposition~\ref{prop:inv} is that the fitted values are unchanged by rigid motions of the coordinates: for corresponding minimisers $\widehat{\mu}_n$ of $L_n$ and $\widetilde{\mu}_n = \widehat{\mu}_n \circ \phi^{-1}$ of $L_n^{\phi}$,
\begin{align*}
\widetilde{\mu}_n\bigl(\phi(\mathbf{T}_{ij})\bigr) = \widehat{\mu}_n \bigl(\phi^{-1}(\phi(\mathbf{T}_{ij}))\bigr) = \widehat{\mu}_n(\mathbf{T}_{ij}) \quad (j=1, \ldots, m_i; \ i=1, \ldots, n).
\end{align*}
The invariance of the fitted values under rotations and reflections has been noted informally in the literature \citep[see, e.g.,][p.~214]{Wood:2017}, but, to the best of our knowledge, Proposition~\ref{prop:inv} is the first formal statement of the full equivariance relation between $\widetilde{\mu}_n$ and $\widehat{\mu}_n$ under rigid motions. It highlights an appealing property of thin-plate estimators: they depend only on the spatial configuration of the data, and not on the choice of origin or orientation of the coordinate system. Other estimators, such as additive ones, are not equivariant in the sense of Proposition~\ref{prop:inv}, because they depend on the coordinate axes.

\iffalse
and additionally require $d$ separate smoothing parameters, as opposed to the single parameter $\lambda$ needed herein. We refer to \citet[Sections 5.6--5.8]{Wood:2017} for an insightful discussion of the advantages and disadvantages of each approach.
\fi
\section{Asymptotic Properties}
\label{sec:3}
\subsection{Fixed Design}

We now derive the rates of convergence of the estimators. As is conventional in the literature \citep[see, e.g.][]{Cai:2011, Xiao:2020}, we differentiate between a \textit{fixed} and a \textit{random} design. While we will not be mentioning the unisolvency condition on $\mathcal{T}$ explicitly in this section, the reader should be reassured that it is satisfied in our set-up for all large $n$. In particular, for fixed designs we will be assuming that the $\{\mathbf{T}_{ij}\}_{i = 1, j = 1}^{n, m_i}$ cover $\mathcal{I}$ increasingly well as $n \to \infty$, from which, via \citet[Theorem~3.14]{Wend:2005}, unisolvency follows for all large $n$. For random designs, we will be assuming that the $\mathbf{T}_{ij}$ are random vectors with a Lebesgue density on $\mathcal{I}$ so that unisolvency will be satisfied with probability one, see the discussion after Proposition~\ref{prop:1}.

We assume throughout this section that $G(x,y) = \rho(y-x)$ so that we are estimating functional location. In fixed designs, the functional data are observed on the same multidimensional grid, viz,
\begin{align*}
\mathbf{T}_{1j} = \mathbf{T}_{2j} = \ldots = \mathbf{T}_{nj}, \quad (j=1, \ldots, m),
\end{align*}
for some $m \in \mathbbm{N}$. To lighten the notation we drop the second subscript and merely write $\mathbf{T}_1, \ldots, \mathbf{T}_m$. We also adopt the shorthand notation
\begin{align*}
\epsilon_{i}(\mathbf{t}) = \left(X_i(\mathbf{t}) - \mu_0(\mathbf{t})\right)+\zeta_i(\mathbf{t}), \quad  (i=1, \ldots, n),
\end{align*}
and write $\epsilon_{ij} = \epsilon_i(\mathbf{T}_{ij})$. Note that, for $i \neq i^{\prime}$, $\epsilon_{ij}$ and $\epsilon_{i^{\prime} j}$ are i.i.d. random variables.  We call $\epsilon_i(\omega, \cdot) = (X_i(\omega,\cdot) - \mu_0(\cdot))+\zeta_i(\omega,\cdot)$ the \textit{error process}. The assumptions that we require are as follows.
\begin{enumerate}[label=(A\arabic*), ref=(A\arabic*)]
\item \label{A1} The discretization points $\mathbf{T}_1, \ldots, \mathbf{T}_m$  are contained in an open, bounded and connected $\mathcal{I} \subset \mathbbm{R}^d$ with Lipschitz boundary.
\item \label{A2} Define $h_{\max} = \sup_{\mathbf{t} \in \mathcal{I}} \inf_{1 \leq j \leq m} \|\mathbf{t}-\mathbf{T}_j\|_{\mathbbm{R}^d}$ and $h_{\min} = \min_{i \neq j} \|\mathbf{T}_i - \mathbf{T}_j \|_{\mathbbm{R}^d}$. There exists a constant $B>0$ such that $h_{\max}/h_{\min} \leq B$.
\item \label{A3} $\rho$ is an absolutely continuous convex and unbounded function on $\mathbbm{R}$ with derivative $\psi$ existing almost everywhere. % \SN{$\rho$ was not used before} \footnote{See opening remark in 2nd paragraph of this section!}
\item \label{A4} There exist finite constants $\kappa$ and $M_1$ such that for all $x \in \mathbbm{R}$ and $|y| \leq \kappa$,
\begin{align*}
|\psi(x+y)-\psi(x)| \leq M_1.
\end{align*}
\item \label{A5}  The error processes $\epsilon_i(\omega,\mathbf{t}) : (\Omega, \mathcal{A}, \mathbb{P}) \times \mathcal{I} \to \mathbbm{R}$ are i.i.d.
\item \label{A6} For $u \to 0$, it holds that
\begin{align*}
\sup_{\mathbf{t} \in \mathcal{I}}\mathbb{E} \left\{\left|\psi(\epsilon_1(\mathbf{t})+u)-\psi(\epsilon_1(\mathbf{t}))\right|^2 \right\} = O(|u|).
\end{align*}
\item \label{A7} $\mathbb{E}\{\psi(\epsilon_{1}(\mathbf{t}))\} =0$ for all $\mathbf{t} \in \mathcal{I}$, \ $\sup_{\mathbf{t} \in \mathcal{I}} \mathbb{E}\{|\psi(\epsilon_{1}(\mathbf{t}))|^2\} < \infty$ and there exists a function $\delta: \mathcal{I} \to \mathbbm{R}_{+}$ satisfying  $0<\inf_{\mathbf{t} \in \mathcal{I}} \delta(\mathbf{t}) \leq \sup_{\mathbf{t} \in \mathcal{I}}  \delta(\mathbf{t}) <\infty$ and
\begin{align*}
\sup_{\mathbf{t} \in \mathcal{I}}\left| \mathbb{E}\{\psi(\epsilon_{1}(\mathbf{t})+u)\}-\delta(\mathbf{t}) u \right| = o(u), \quad \mbox{as $u \to 0$.}
\end{align*}
\end{enumerate}

Assumption \ref{A1} requires that the boundary of the domain $\mathcal{I}$ is locally a graph of a Lipschitz function---this is a standard assumption in the theory of Sobolev spaces. It is satisfied, for example, if $\mathcal{I}$ is an open ball or a rectangle. The quantities $h_{\max}$ and $h_{\min}$ in \ref{A2} are referred to as ``fill" and ``separation" distances, respectively; they are critical elements of approximation theory \citep[see, e.g.][]{Wend:2005}. The requirement that the ratio $h_{\max}/h_{\min}$ is bounded is referred to as ``quasi-uniformity" and ensures that $\mathcal{I}$ is well-covered by the $\{\mathbf{T}_j\}_{j=1}^m$.

Assumption \ref{A3} specifies the type of a loss function that may be used; many interesting loss functions, such as the square, Huber, absolute and logistic losses, satisfy this assumption. Moreover, as \ref{A3} does not require the symmetry of $\rho$, asymmetric loss functions such as the quantile and expectile losses, are also permitted within our framework. The derivative of convex functions exists almost everywhere and is a monotone function. Assumption \ref{A4} requires that the sizes of its jumps (the only type of discontinuity that monotone functions can have) are uniformly bounded. This assumption is trivially satisfied if $\psi$ is bounded or uniformly continuous, e.g., if $\psi(x) = \sign(x)$ or $\psi(x) = x$. Assumption \ref{A5} is very general and is fulfilled, for example, if the $X_i$ and $\zeta_i$ are i.i.d. and independent of each other, as was assumed by \citet{Cai:2011} and \citet{Xiao:2020}.
 
Assumption \ref{A6} requires that $\psi(\epsilon_{1}(\mathbf{t})+u)$ is uniformly mean-square Lipschitz continuous over $\mathbf{t} \in \mathcal{I}$ at $u= 0$. It is a minimal condition that is satisfied (with the stronger $O(|u|^2)$ on the right-hand side) if $\psi$ is a Lipschitz function, e.g., $\psi(x) = x$ or $\psi_k(x) = \max(-k,\min(x,k))$ for some $k>0$. However, assumption \ref{A6} is general enough to cover discontinuous $\psi$-functions by trading regularity between $\psi$ and the one-dimensional distributions of $\epsilon_1(\mathbf{t})$ at each $\mathbf{t} \in \mathcal{I}$. For example, for $\psi(x) = \sign(x)$, if $\epsilon_1(\mathbf{t})$ at each $\mathbf{t} \in \mathcal{I}$ has CDF $F_{\mathbf{t}}$ with positive and bounded Lebesgue density $f_{\mathbf{t}}$ in a neighborhood about zero, for $u \to 0$ we find
\begin{align*}
\sup_{\mathbf{t} \in \mathcal{I}}\mathbb{E} \left\{\left|\psi(\epsilon_1(\mathbf{t})+u)-\psi(\epsilon_1(\mathbf{t}))\right|^2 \right\} \leq 4 |u| (\sup_{\mathbf{t} \in \mathcal{I}} \sup_{|u|\leq \delta} f_{\mathbf{t}}(u)),
\end{align*}
for some small $\delta>0$. If $\sup_{\mathbf{t} \in \mathcal{I}} \sup_{|u|\leq \delta} f_{\mathbf{t}}(u) < \infty$ for some small $\delta>0$, this is clearly $O(|u|)$.

Assumption \ref{A7} is similar in spirit to \ref{A6}. It is required that the expectation of the score function $\psi(\epsilon_1(\mathbf{t}))$ is zero at all $\mathbf{t} \in \mathcal{I}$ and the squared expectation is finite, both of which are well-known conditions in the asymptotics of maximum likelihood and M-estimators. For bounded $\psi$, the squared expectation will always be finite irrespective of the distribution of $\epsilon_1(\mathbf{t})$. The second part of \ref{A7} is a uniform linearization condition of $\mathbb{E}\{\psi(\epsilon_{1}(\mathbf{t})+u)\}$ for small $u>0$. It can be verified, for example, if $\psi \in \mathcal{C}^2(\mathbbm{R})$ with $\psi^{\prime \prime}$ bounded and $0<\mathbb{E} \{\psi^{\prime}(\epsilon_1(\mathbf{t}))\} <\infty$ for all $\mathbf{t} \in \mathcal{I}$, in which case we can take $\delta(\mathbf{t}) = \mathbb{E} \{\psi^{\prime}(\epsilon_1(\mathbf{t}))\}$. For $\psi(x) = x$, we can simply take $\delta(\mathbf{t}) := 1$. However, like \ref{A6}, condition \ref{A7} has been designed with non-smooth $\psi$-functions in mind. For example, for $\psi(x) = \sign(x)$ if $F_\mathbf{t}(0) = 1/2$, that is, $\epsilon_1(\mathbf{t})$ has median at zero but is not necessarily symmetric, then choosing $\delta(\mathbf{t}) = f_{\mathbf{t}}(0)$, the condition is fulfilled provided that
\begin{align*}
\lim_{u \to 0} \sup_{\mathbf{t} \in \mathcal{I}} \left| \frac{F_{\mathbf{t}}(u)-F_{\mathbf{t}}(0)}{u} - f_{\mathbf{t}}(0)  \right| = 0.
\end{align*}
In other words, the condition is fulfilled if the family of CDFs $\{F_{\mathbf{t}}, \mathbf{t} \in \mathcal{I}\}$ is uniformly first order differentiable at $u=0$ and $\{f_{\mathbf{t}}(0), \mathbf{t} \in \mathcal{I}\}$ (the derivatives/densities of $F_{\mathbf{t}}$ at $0$), are bounded away from zero and infinity. These conditions generalize classical conditions for the consistency of quantile regression estimators, compare with \citet[Chapter 4]{Koenker:2005}.  It is worth highlighting that these conditions only need to hold for small $|u|$, not globally. 

As noted previously, the restriction of $f \in \BL_r(\mathbbm{R}^d)$ to a bounded $\mathcal{I} \subset \mathbbm{R}^d$ is an element of the Sobolev space of functions on $\mathcal{I}$, $\mathcal{H}^r(\mathcal{I})$. Thus, to derive our theoretical results we also need to work with the space $\mathcal{H}^{r}(\mathcal{I})$. It is convenient to equip $\mathcal{H}^r(\mathcal{I})$ with the penalty-weighted inner product
\begin{align}
\label{eq:sobip}
    \langle f, g \rangle_{\mathcal{H}^r(\mathcal{I}), \lambda} := \langle f, g \rangle_{L_2(\mathcal{I})} + \lambda \sum_{|\alpha| = r} \frac{r!}{\alpha!} \langle \partial^{\alpha} f, \partial^{\alpha} g \rangle_{L_2(\mathcal{I})}.
\end{align}
The induced norm is denoted with $\|\cdot \|_{\mathcal{H}^r(\mathcal{I}), \lambda}$. For $\lambda>0$ it is easy to see that $\|\cdot\|_{\mathcal{H}^r(\mathcal{I})}$ and $\|\cdot\|_{\mathcal{H}^r(\mathcal{I}), \lambda}$ are equivalent norms and so $\mathcal{H}^r(\mathcal{I})$ equipped with $\langle \cdot, \cdot \rangle_{\mathcal{H}^r(\mathcal{I}), \lambda}$ is also an RKHS. The smallest embedding constant of $\mathcal{H}^r(\mathcal{I}) \hookrightarrow \mathcal{C}(\mathcal{I})$ with the former space equipped with $\langle \cdot, \cdot \rangle_{\mathcal{H}^r(\mathcal{I}), \lambda}$  is critical in our proofs. Proposition~\ref{prop:3} shows that this constant (constant in terms of $f \in \mathcal{H}^r(\mathcal{I})$)  is proportional to $\lambda^{-d/(4r)}$ for all $\lambda \in (0,1]$. 
\begin{prop}
\label{prop:3}
Suppose that \ref{A1} is satisfied. Then, for $2r>d$ there exists a constant $c_0>0$ such that, for all $\lambda \in (0,1]$ and $f \in \mathcal{H}^r(\mathcal{I})$,
\begin{align}
\label{eq:emb}
\sup_{ \mathbf{t} \in \mathcal{I}} \left|f(\mathbf{t}) \right| \leq c_0 \lambda^{-d/(4r)} \|f\|_{\mathcal{H}^r(\mathcal{I}), \lambda}.
\end{align}
\end{prop}
\noindent
The condition that $\lambda \in (0,1]$ is not restrictive, as in our theoretical results below we assume that $\lambda \to 0$, as $n \to \infty$, and our results are asymptotic in nature. With this preparation in place, we now state the main result of this section. 

\begin{theorem}
\label{thm:2}
Suppose that $2r >d$ and that assumptions \ref{A1}--\ref{A6} are satisfied. If $h_{\max} \to 0$, $\lambda \sim C h_{\max}^{2r}$ for some large $C>0$ and  $n h_{\max}^{4d} \to \infty$ as $n \to \infty$, then
\begin{align*}
\int_{\mathcal{I}} |\widehat{\mu}_n(\mathbf{t})-\mu_0(\mathbf{t})|^2 \dd \mathbf{t} = O_{\mathbb{P}}\left(\frac{1}{n}  + h_{\max}^{2r} \right), \quad \mbox{as $n \to \infty$.}
\end{align*}
\end{theorem}

The proof of Theorem~\ref{thm:2} relies on an asymptotic linearisation of the objective function. We show that a local minimiser exists in the interior of the ball centered at $\mu_0$ with radius proportional to $n^{-1}+h_{\max}^{2r}$, which then, by convexity, is a global minimiser $\widehat{\mu}_n$. Our proof makes use of the equivalent norm in $\|\cdot\|_{\BL_r(\mathbbm{R}^d),\lambda}$ defined in \eqref{eq:equivnorm}, with the equivalence of $\|\cdot\|_{\BL_r(\mathbbm{R}^d),\lambda}$ and $\|\cdot\|_n$ following from unisolvency of the $\mathbf{T}_{ij}$. The condition $n h_{\max}^{4d} \to \infty$ ensures that the remainder of this expansion is asymptotically negligible without imposing any distributional assumptions within the blocks $\{\epsilon_{ij}\}_{j=1}^{m_i}$; it is the price to pay for the generality of our treatment. It can be relaxed if $\psi$ is more regular; for example, if $\psi$ is Lipschitz then $n h_{\max}^{3d} \to \infty$ suffices and if $\psi(x) = x$, the non-linear remainder of the linearisation vanishes and the condition may be omitted entirely.

The convergence rate in Theorem~\ref{thm:2} reveals a sharp phase transition between the sample size $n$ and the discretization of the domain $\mathcal{I}$, as captured by the fill distance $h_{\max}$. Specifically, the $L_2(\mathcal{I})$ error decomposes as the sum of $n^{-1}$ and $h_{\max}^{2r}$ implying that the dominant term depends on the relative magnitude of these quantities. When the grid is relatively sparse, i.e., $h_{\max} \geq n^{-1/(2r)}$ it is the discretization error $h_{\max}^{2r}$ that dominates the asymptotics and $\|\widehat{\mu}_n-\mu_0\|_{L_2(\mathcal{I})}^2 = O_{\mathbb{P}}(h_{\max}^{2r})$. However, if the discretization is sufficiently fine, i.e., $h_{\max} < n^{-1/(2r)}$, Theorem~\ref{thm:2} leads to
\begin{align*}
\int_{\mathcal{I}} |\widehat{\mu}_n(\mathbf{t})-\mu_0(\mathbf{t})|^2 \dd \mathbf{t} = O_{\mathbb{P}}\left(\frac{1}{n} \right),
\end{align*}
that is, $\widehat{\mu}_n$ converges to $\mu_0$ at the parametric rate. 

While parametric rates of convergence for functional location estimators are not new in the literature \citep[see, e.g.][]{Horvath_Kokoszka2012, Ch:2014}, most estimators require completely observed noiseless functional data. By contrast, the present estimator allows for both discretely sampled and noisy data and still attains the parametric convergence rate. In this light, Theorem~\ref{thm:2} is a positive result showing that the fastest possible rate of convergence may even be obtained in the real world where the data is almost never observed in its entirety or without error. The theory even points to the threshold that needs to be maintained by $h_{\max}$ for the parametric rate to occur, i.e.,  $n^{-1/(2r)}$. A useful heuristic, which is exact if $\mathcal{I} = (0,1)^d$, is $h_{\max} \approx m^{-1/d}$ so that to attain the parametric rate it is sufficient that $m>n^{d/(2r)}$. This threshold is larger for high-dimensional data (larger $d$), but smaller for smoother mean functions (larger $r$).

The discretization error $h_{\max}^{2r}$ is the order of the interpolation error for thin-plate splines. Specifically, if $\mu_0 \in \BL_r(\mathbbm{R}^d)$ and $\mu_0(\mathbf{T}_1), \ldots, \mu_0(\mathbf{T}_m)$ were directly observed, an interpolating thin-plate spline $\mu_{\text{interp}}$ would satisfy
\begin{align*}
\int_{\mathcal{I}} |\mu_{\text{interp}}(\mathbf{t}) - \mu_0(\mathbf{t})|^2 \dd \mathbf{t} \leq C h_{\max}^{2r} |\mu_0|^2_{\BL_r(\mathbbm{R}^d)},
\end{align*}
as $h_{\max} \to 0$ with a global constant $C>0$ depending only on $\mathcal{I}$ \citep[Proposition 3]{Duchon:1978}. Of course, in practice $\mu_0$ is unobservable, so in addition to the interpolation error $h_{\max}^{2r}$ we have to contend with the statistical error, which per Theorem~\ref{thm:2} is of order $n^{-1}$. Thus, the rate $n^{-1}+h_{\max}^{2r}$ can be interpreted as a decomposition into variability and interpolation error.

\subsection{Random Design}

We now turn to the random or longitudinal-type of design wherein the discretization vectors $\mathbf{T}_{ij} \in \mathcal{I}$ are random. This captures situations, such as patient follow-up, in which the discretization vectors for each random field may be vastly different. We shall assume henceforth that the $\mathbf{T}_{ij}$ are i.i.d. draws from an absolutely continuous distribution on $\mathcal{I} \subset \mathbbm{R}^d$. To establish sharp rates of convergence, we first note an important consequence of Proposition~\ref{prop:3}: by \eqref{eq:emb}, the evaluation functionals $\delta_\mathbf{t}(f) = f(\mathbf{t})$ are continuous when $\mathcal{H}^r(\mathcal{I})$ is equipped with $\langle \cdot, \cdot \rangle_{\mathcal{H}^r(\mathcal{I}),\lambda}$.  Therefore, by the fundamental existence theorem for reproducing kernel Hilbert spaces \citep[Section 2]{Ar:1950}, there exists a symmetric function  $\mathcal{R}_{\mathcal{H}^r(\mathcal{I}), \lambda}: \mathcal{I} \times \mathcal{I} \to \mathbbm{R}$,  the reproducing kernel, such that for every $\mathbf{y} \in \mathcal{I}$, $\mathbf{x} \mapsto \mathcal{R}_{\mathcal{H}^r(\mathcal{I}), \lambda}(\mathbf{x}, \mathbf{y}) \in \mathcal{H}^r(\mathcal{I})$ and
\begin{align}
\label{eq:rk}
f(\mathbf{t}) = \langle \mathcal{R}_{\mathcal{H}^r(\mathcal{I}),\lambda}(\mathbf{t}, \cdot), f \rangle_{\mathcal{H}^r(\mathcal{I}),\lambda},
\end{align}
for all $\mathbf{t} \in \mathcal{I}$ and $f \in \mathcal{H}^r(\mathcal{I})$. For a complete understanding of the asymptotic behaviour of our estimators, we need to study the reproducing kernel $\mathcal{R}_{\mathcal{H}^r(\mathcal{I}),\lambda}$. Proposition~\ref{prop:4} derives a Mercer-type expansion significantly expanding existing results in the literature,  cf. \citet[Theorem~2.8.4]{Hsing:2015}.

\begin{prop}
\label{prop:4}
Suppose that $\mathcal{I} \subset \mathbbm{R}^d$ is open, bounded and connected, $2r>d$ and $\mathcal{H}^r(\mathcal{I})$ is equipped with $\langle \cdot, \cdot \rangle_{\mathcal{H}^r(\mathcal{I}),\lambda}$. Then, the following holds true.
\begin{enumerate}[label=$(\roman*)$, ref=$(\roman*)$]
\item \label{i} There exists a complete orthonormal system of $L_2(\mathcal{I})$, $\{\phi_j\}_{j=1}^{\infty}$, not depending on $\lambda$, such that 
\begin{align*}
\langle \phi_i, \phi_j \rangle_{L_2(\mathcal{I})} = \delta_{ij}, \quad \text{and} \quad \langle \phi_i, \phi_j \rangle_{\mathcal{H}^r(\mathcal{I}), \lambda} = (1+\lambda \gamma_i) \delta_{ij},
\end{align*}
where $\delta_{ij}$ is the Kronecker delta, the sequence $\{\gamma_j\}_{j=1}^{\infty}$ satisfies $\gamma_1 = \ldots = \gamma_M = 0$ for $M = \binom{r+d-1}{d}$, and
\begin{align*}
C_1 j^{2r/d} \leq \gamma_j \leq C_2 j^{2r/d}, \quad (j = M+1, \ldots),
\end{align*}
for some $0<C_1 \leq C_2<\infty$.
\item \label{ii} The functions $\{\phi_j/(1+\lambda \gamma_j)^{1/2}\}_{j=1}^{\infty}$ form a complete orthonormal system for $\mathcal{H}^r(\mathcal{I})$ and $\{\phi_j\}_{j=1}^M$ form an orthonormal basis for the $M$-dimensional space of polynomials of total order $r$ on $\mathcal{I}$, $\mathcal{P}_r(\mathcal{I})$.
\item \label{iii} The reproducing kernel, $\mathcal{R}_{\mathcal{H}^r(\mathcal{I}), \lambda}(\mathbf{x}, \mathbf{y})$\, may be written as
\begin{align*}
\mathcal{R}_{\mathcal{H}^r(\mathcal{I}), \lambda}(\mathbf{x}, \mathbf{y}) = \sum_{j=1}^{\infty} \frac{\phi_j(\mathbf{x}) \phi_j(\mathbf{y})}{1+\lambda \gamma_j},
\end{align*}
with the series converging (at least) in $L_2(\mathcal{I} \times \mathcal{I})$.
\end{enumerate}
\end{prop}
\noindent
We establish these results through the well-developed theory of interpolation; the interested reader is referred to \citet{triebel1995} and our proofs for the complete arguments. Conditions \ref{A1} and \ref{A2} for the fixed design need to be replaced with \ref{B1}, given below.

\begin{enumerate}[label=(B\arabic*), ref=(B\arabic*)]
\item \label{B1} The discretization vectors $\{\mathbf{T}_{ij}\}_{i = 1, j = 1}^{n, m_i}$ are i.i.d. with Lebesgue density $f$ supported on $\mathcal{I}$, satisfying  $c \leq f(\mathbf{t}) \leq C$ for some $0<c\leq C<\infty$ for all $\mathbf{t} \in \mathcal{I}$. Moreover, the $\{\mathbf{T}_{ij}\}_{i = 1, j = 1}^{n, m_i}$ are independent of the error processes $\{\epsilon_{i}\}_{i=1}^n$.
\end{enumerate}

Condition \ref{B1} is a multidimensional generalization of the respective conditions found in \citet{Cai:2011, Zhang:2016, Xiao:2020} and is natural in our setting. In the statement of Theorem~\ref{thm:3} below, $m$ denotes the harmonic average of the $\{m_i\}_{i=1}^n$, namely, $m= n/(\sum_{i=1}^n m_i^{-1})$. One may think of $m$ as an overall measure of density of the discretization vectors: larger values imply that the random fields $X_1, \ldots, X_n$ are more completely observed on average (or on harmonic average).

\iffalse
Before stating the main result, it is worth highlighting a technical point involving conditional expectations that frequently appears in the proof of Theorem~\ref{thm:3}. Specifically, in the proof of Theorem~\ref{thm:3} we need to deal with expressions such as
\begin{align*}
u \mapsto \mathbb{E}\{\psi(\epsilon_{i}(\mathbf{T}_{ij})+u)|\mathcal{T}\},
\end{align*}
where $\mathcal{T} = \{\mathbf{T}_{ij}\}_{i = 1, j = 1}^{n, m_i}$ is the sequence of knots of our splines (I think we should just stick to calling them one thing - How about \textit{sampling vectors}? IK.). By independence of $\{\epsilon_i\}_{i=1}^n$ from the sampling vectors $\{\mathbf{T}_{ij}\}_{i = 1, j = 1}^{n, m_i}$,the identical distributions of $\epsilon_i$ and the transformation thoerem, these expressions simplify to 
\begin{align*}
u \mapsto \int_{\mathbbm{R}} \psi(y+u) \dd F_{\mathbf{T}_{ij}}(y),
\end{align*}
where $F_{\mathbf{T}_{ij}}$ denotes the distribution of the (real) random variable $\epsilon_1(\mathbf{T}_{ij})$. Assumptions \ref{A6} and \ref{A7} ensure that these maps are uniformly regular by imposing smoothness and linearization simultaneously over all coordinates $\mathbf{t} \in \mathcal{I}$. This uniform regularity is crucial for controlling stochastic fluctuations and establishing rates of convergence for the resulting M-estimators and cannot be relaxed.
\fi

\begin{theorem}
\label{thm:3}
Suppose that $\mathcal{I} \subset \mathbbm{R}^d$ is open, bounded with Lipschitz boundary and connected, $2r>d$ and assumptions \ref{A3}--\ref{A6} and \ref{B1} are satisfied. If $\lambda \to 0$ in such a way that $n \lambda^{2d/r} \to \infty$ then
\begin{align*}
\int_{\mathcal{I}} | \widehat{\mu}_n(\mathbf{t}) - \mu_0(\mathbf{t})|^2  \dd \mathbf{t}= O_{\mathbb{P}}\left(\frac{1}{n} + \frac{1}{nm \lambda^{d/(2r)}} + \lambda \right), \quad \mbox{as $n \to \infty$.}
\end{align*}
\end{theorem}

As in the situation of the fixed design, Theorem~\ref{thm:3} reveals that the parametric rate of convergence is attainable. In particular, taking $\lambda \sim (nm)^{-2r/(2r+d)}$, Theorem~\ref{thm:3} implies
\begin{align*}
\int_{\mathcal{I}} | \widehat{\mu}_n(\mathbf{t}) - \mu_0(\mathbf{t})|^2  \dd \mathbf{t}= O_{\mathbb{P}}\left(\frac{1}{n} + \frac{1}{(nm)^{2r/(2r+d)}} \right).
\end{align*}
Thus, if $m \geq n^{d/2r}$ the $n^{-1}$ component of the error dominates, leading to $\|\widehat{\mu}_n-\mu_0\|_{L_2(\mathcal{I})}^2 = O_{\mathbb{P}}(n^{-1})$. This is an even more remarkable result than before considering how disparately observed the random fields $X_1, \ldots, X_n$ may be. By contrast, if $X_1, \ldots, X_n$ are relatively sparsely sampled, that is, if $m< n^{d/(2r)}$, we find  $\|\widehat{\mu}_n-\mu_0\|_{L_2(\mathcal{I})}^2 = O_{\mathbb{P}}(( nm)^{-2r/(2r+d)})$, which is the well-known optimal nonparametric rate of convergence for the estimation of $\mathcal{H}^r(\mathcal{I})$-functions from $nm$ independent observations. This rather interesting finding shows that, unlike the corresponding situation in the fixed design, the effects of the sample size and the discretization are intertwined.

We note that, although Theorems \ref{thm:2} and \ref{thm:3} were stated with convex loss functions in mind, they can be generalized to non-convex loss functions that are continuous and coercive. For such loss functions, our proofs still go through and establish local optimality. In particular, in the case of a fixed design, the proof of Theorem~\ref{thm:2} yields
\begin{align*}
    \lim_{D \to \infty} \liminf_{n \to \infty} 
    \mathbb{P}\Big(&\text{a local minimiser}\ \widehat{\mu}_n \ \text{exists and satisfies} \  \|\widehat{\mu}_n - \mu_0\|_{L_2(\mathcal{I})}^2 
    \leq D(n^{-1} + h_{\max}^{2r}) \Big) = 1.
\end{align*}
Thus, there exists a local minimiser enjoying the same rates of convergence as in the convex case. An analogous result holds for the random design (Theorem \ref{thm:3}).

\section{Computation and Smoothing Parameter Selection}
\label{sec:4}

We describe our computational algorithm for distinct sampling vectors across the functional data $X_1, \ldots, X_n$ and $m_1, \ldots, m_n$; the case of functional data observed on the same grid is a special case thereof. Denote by $m^{\star} = \sum_{i=1}^n m_i$ the total number of discretized observations, and let $\left\{ \mathbf{T}_j \right\}_{j=1}^{m_0} = \bigcup_{i = 1}^n \bigcup_{j=1}^{m_i} \left\{ \mathbf{T}_{ij} \right\}$ be an enumeration of the $m_0 \leq m^{\star}$ unique observation locations in $\mathcal{I}$ for all functional data. A general thin-plate spline function $f \colon \mathbbm{R}^d \to \mathbbm{R}$ from~\eqref{eq:TP} evaluated at $\mathbf{T}_j$ can be written as
    \[
    f(\mathbf{T}_j) = \sum_{\ell=1}^{m_0} \gamma_{\ell} \eta_{dr} \left(\left\|\mathbf{T}_j-\mathbf{T}_{\ell} \right\|_{\mathbbm{R}^d} \right) + \sum_{k=1}^M \delta_k \phi_k\left(\mathbf{T}_j\right), \quad (j = 1, \dots, m_0),
    \]
which can be expressed in vector notation as
    \begin{equation} \label{eq: f unequal design}
    \left(f\left(\mathbf{T}_1\right), \dots, f\left(\mathbf{T}_{m_0} \right) \right)^\top = \mathbf{K} \boldsymbol{\gamma} + \mathbf{\Phi} \boldsymbol{\delta}.
    \end{equation}
Here, $\mathbf{K} \in \mathbbm{R}^{m_0 \times m_0}$ is the symmetric matrix with entries $\mathbf{K}_{ij} = \eta_{dr}\left( \left\Vert \mathbf{T}_i - \mathbf{T}_j \right\Vert_{\mathbbm{R}^d} \right)$, $i, j = 1, \dots, m_0$, the matrix $\boldsymbol{\Phi} \in \mathbbm{R}^{m_0 \times M}$ is given by $\boldsymbol{\Phi}_{ij} = \phi_j(\mathbf{T}_i)$, $i = 1, \dots, m_0$, $j = 1, \dots, M$, and the vectors of coefficients to be estimated are denoted by $\boldsymbol{\gamma} = \left( \gamma_1, \dots, \gamma_{m_0} \right)^\top \in \mathbbm{R}^{m_0}$ and $\boldsymbol{\delta} = \left( \delta_1, \dots, \delta_M \right)^\top \in \mathbbm{R}^M$.

The requirement for $f$ to satisfy $|f|_{\BL_r(\mathbbm{R}^d)}< \infty$ is equivalent to $\mathbf{\Phi}^{\top} \boldsymbol{\gamma} = \mathbf{0}$ \citep[see, e.g.][Chapter 7]{Green:1994}. This constraint can be incorporated into the problem by reparametrizing $\boldsymbol{\gamma} = \boldsymbol{Q} \boldsymbol{\xi}$, where $\boldsymbol{Q} \in \mathbbm{R}^{m_0 \times (m_0 - M)}$ is a matrix whose columns span the null space of $\boldsymbol{\Phi}^\top$, and $\boldsymbol{\xi} \in \mathbbm{R}^{m_0 - M}$. Under the constraint $|f|_{\BL_r(\mathbbm{R}^d)}< \infty$, \citet{Green:1994} show that $|f|_{\BL_r(\mathbbm{R}^d)}^2 = \boldsymbol{\xi}^\top \boldsymbol{Q}^\top \mathbf{K} \boldsymbol{Q} \boldsymbol{\xi}$. \iffalse To substitute the expression~\eqref{eq: f unequal design} into the objective function~\eqref{eq:est}, we need additional notation.\fi For each $i = 1, \dots, n$, denote by $\pi_i \colon \left\{ 1, \dots, m_i \right\} \to \left\{ 1, \dots, m_0 \right\}$ the map that takes the observation point $\mathbf{T}_{ij}$ corresponding to the response value $Y_{ij}$ to the associated index of the sorted vector of all observation times $\left( \mathbf{T}_1, \dots, \mathbf{T}_{m_0} \right)^\top$. Then $Y_{ij} - f(\mathbf{T}_{ij}) = Y_{ij} - f(\mathbf{T}_{\pi_i(j)})$, and the loss function~\eqref{eq:est} can be expressed as
    \begin{equation} \label{eq:Mest2}
    \begin{aligned}
    L_n(\boldsymbol{\theta}) & = \frac{1}{n} \sum_{i=1}^n \frac{1}{m_i} \sum_{j=1}^{m_i} \rho \left(Y_{ij} - f(\mathbf{T}_{\pi_i(j)}) \right) + \lambda \, \boldsymbol{\xi}^\top \boldsymbol{Q}^\top \mathbf{K} \boldsymbol{Q} \boldsymbol{\xi} \\
    & = 
    \sum_{i=1}^n \sum_{j=1}^{m_i} \frac{1}{n\,m_i} \rho \left(Y_{ij} - \mathbf{Z}_{\pi_i(j)}^{\top} \boldsymbol{\theta} \right) + \lambda \, \boldsymbol{\theta}^{\top} \mathbf{H} \boldsymbol{\theta}.
    \end{aligned}
    \end{equation}
In this expression, $\boldsymbol{\theta} = (\boldsymbol{\xi}^{\top}, \boldsymbol{\delta}^{\top})^{\top} \in \mathbbm{R}^{m_0}$ is the complete vector of parameters to estimate, the penalty matrix $\mathbf{H} \in \mathbbm{R}^{m_0 \times m_0}$ is given by
    \begin{equation*} % \label{eq: H matrix}
    \mathbf{H} = \begin{bmatrix} 
    \mathbf{Q}^{\top} \mathbf{K} \mathbf{Q} & \mathbf{0} \\
    \mathbf{0} & \mathbf{0} 
    \end{bmatrix},
    \end{equation*}
and $\mathbf{Z}_j$ is the $j$th row of the compressed design matrix $\mathbf{Z}_0 = \left[ \mathbf{K} \mathbf{Q}, \mathbf{\Phi} \right] \in \mathbbm{R}^{m_0 \times m_0}$. The complete design matrix $\mathbf{Z} \in \mathbbm{R}^{m^{\star} \times m_0}$ is composed of the rows $\mathbf{Z}_{\pi_i(j)}$ for all $i = 1, \dots, n$ and $j = 1, \dots, m_i$, in the same order as they appear in the complete vector of function values $Y_{ij}$, $i = 1, \dots, n$, $j = 1, \dots, m_i$.

%\subsection{Finite-Dimensional Problem Solution and Smoothing Parameter Selection}
%\label{subsec:IRLS_and_QP}
Having expressed the objective function as in \eqref{eq:Mest2}, we have reduced the variational problem \eqref{eq:est} to a much simpler finite-dimensional weighted and penalized M-estimation problem. The practical minimisation of $L_n(\boldsymbol{\theta})$ depends on the specific choice of the loss function $\rho(x)$. For the square loss $\rho(x) = x^2$, the objective function \eqref{eq:Mest2} admits the familiar closed-form minimiser $\widehat{\boldsymbol{\theta}}_n = \left(\mathbf{Z}^{\top}\mathbf{Z} + \lambda \mathbf{H}\right)^{-1} \mathbf{Z}^{\top} \mathbf{Y}$, where $\mathbf{Y} \in \mathbb{R}^{m^{\star}}$ is the vector of all observed functional values $Y_{ij}$. For other loss functions $\rho$, the objective function $L_n(\boldsymbol{\theta})$ can be minimised iteratively through a weighted version of the iteratively reweighted least squares (IRLS) algorithm \citep[Chapter~9]{Maronna:2019}, adapted analogously from \citet{K:2025}. This approach is classical in the robust regression literature, owing to its ease of implementation \citep{Green:1984};
%, and has a computational cost of $O\big(It_{\mathrm{IRLS}} (m^{\star} m_{0}^{2} + m_{0}^{3})\big)$, where $It_{\mathrm{IRLS}}$ denotes the number of iterations to convergence\footnote{Perhaps this sentence can be omitted - IK}.
however, the IRLS approach requires the loss $\rho$ to be differentiable in order to guarantee convergence, a requirement that fails for the quantile (pinball) loss. When differentiability is lacking, only an approximate (smoothed) version of the problem can typically be solved with IRLS-type algorithms. 

As an alternative to IRLS, we show that the M-estimation problem \eqref{eq:Mest2} admits an exact reformulation as a quadratic program (QP) when $\rho$ is the quantile loss. Let 
\begin{equation}
    \label{eq::quantile_loss}
    \rho_{\tau}(x) = (\tau - \mathbbm{1}\{x<0\})x, \quad x \in \mathbb{R}, \quad \tau \in (0,1), 
\end{equation}
\begin{equation}
\begin{aligned}
    \min_{\boldsymbol{u},\boldsymbol{v},\boldsymbol{\theta}} \quad
    & \sum_{i=1}^{n} \sum_{j=1}^{m_i}
    \frac{1}{n m_i}
    \left\{ \tau u_{ij} + (1-\tau)v_{ij} \right\}
    + \lambda \boldsymbol{\theta}^{\top}\mathbf{H}\boldsymbol{\theta} \\
    \text{subject to} \quad
    & u_{ij} - v_{ij}
    = Y_{ij} - \mathbf{Z}_{\pi_i(j)}^{\top}\boldsymbol{\theta},
    \qquad i=1,\ldots,n,\quad j=1,\ldots,m_i,\\
    & u_{ij} \geq 0,\qquad
    v_{ij} \geq 0,
    \qquad i=1,\ldots,n,\quad j=1,\ldots,m_i.
\end{aligned}
\label{eq:QP_quantile}
\end{equation}
This is a convex QP in $(\boldsymbol{\theta},\boldsymbol{u},\boldsymbol{v})$. The objective is linear in $\boldsymbol{u}$ and $\boldsymbol{v}$ and quadratic only in $\boldsymbol{\theta}$ through the penalty, so its Hessian is positive semidefinite. Such problems are traditionally solved with interior-point and active-sets methods \citep{Koenker:1997}; however, these methods can be numerically fragile when the Hessian is singular and the constraints are degenerate \citep[Section~16]{Nocedal:2006}. To avoid these issues, we instead solve \eqref{eq:QP_quantile} with OSQP \citep{stellato2020osqp}, a modern QP solver based on the alternating direction method of multipliers (ADMM) rather than on interior-point or active-set techniques. OSQP is perfectly suited to our setting because it accommodates convex quadratic programs without requiring a positive definite Hessian, while exploiting the sparsity of the problem for computational efficiency. Moreover, its factorisation caching and warm-starting facilitate repeated solutions over a grid of smoothing parameters $\lambda$ (which we require for cross-validation). These features make OSQP an attractive choice for our robust (functional) estimation problem.

Unlike IRLS, the QP formulation treats the non-differentiable quantile loss exactly, without smoothing, up to the termination tolerance of the solver. The price is memory and cost: the formulation introduces $2m^{\star}$ auxiliary variables $u_{ij},v_{ij}$ in addition to $\boldsymbol{\theta}$, and, as the computation times in Section~\ref{sec: appendix B} of the Supplementary Material show, for smooth losses, the OSQP algorithm is more expensive than its IRLS counterpart. The choice thus involves a trade-off between exactness (and stability of the solver) and computational cost. The same QP machinery may also be employed for the solution of specific differentiable losses $\rho$ like the Huber. In this case, the computational cost is similar to that of the classical IRLS algorithm.

The smoothing parameter $\lambda$ is selected by minimizing the robust generalized cross-validation criterion $\RCV(\lambda)=\tau(\mathbf{r}_{-}(\lambda))^2$, where $\mathbf{r}_{-}(\lambda)$ denotes the vector of exact or approximate leave-one-out residuals and $\tau$ is the robust and efficient $\tau$-scale of \citet{yohai1988high}. We follow the procedure described in \citet{K:2025}, which approximates the leave-one-out residuals using the diagonal of the weighted hat matrix
$\mathbf{H}_{W}(\lambda) := \mathbf{Z} \left(\mathbf{Z}^{\top}\mathbf{W}\mathbf{Z}+\lambda\mathbf{H}\right)^{-1}\mathbf{Z}^{\top}\mathbf{W}$
associated with the fitting procedure. As in \citet{K:2025}, when IRLS is employed, $\mathbf{H}_{W}(\lambda)$ is obtained from the final iteration of the algorithm. For estimators obtained from the quadratic formulation, the same approximation of the leave-one-out residuals is used, with $\mathbf{W}$ given by the weights entering the weighted loss $L_n(\boldsymbol{\theta})$. Since only the diagonal is required, the diagonal of $\mathbf{H}_{W}(\lambda)$ can be estimated stochastically using Hutchinson's method \citep{Bekas:2007}, without explicitly forming the full matrix thereby reducing the computational overhead for large datasets. % \footnote{At this point point the reader to the R package and emphasise the C++ back-end – IK.}
Both these estimation methods (the classical IRLS and the exact QP formulation) are implemented in an efficient and fully documented \textsf{R} package \textsf{RobustSpline}, with a back-end \textsf{C++}.\footnote{Available at \url{https://github.com/NagyStanislav/RobustSpline}.} % in agreement with considerations in Section \ref{sec:4}.

\section{Numerical Results}
\label{sec:5}

We now present numerical results supporting the proposed robust location estimators. We consider four thin-plate spline estimators: % \eqref{eq:1dest}:
\begin{itemize}
        \item[\textsf{S}:] the least-squares estimator with $\rho(x)=x^2$,
        \item[\textsf{A}:] the $L_1$ estimator with $\rho(x) = |x|$,
        \item[\textsf{H}:] the Huber estimator % with $\rho(x)$ given by 
        $\rho(x) = \min( |x|, k )^2 + 2 k \max( 0, |x| - k )$
            % \begin{align}
            % \label{eq:huber}
            %     \rho_{k}(x) = \begin{cases}
            %         x^2 & |x| < k, \\
            %         2k|x|-k^2 & |x| \geq k,
            %     \end{cases}
            % \end{align}
        for $k = 1.345$, mixing the square and $L_1$ losses (this choice of $k$ ensures $95\%$ efficiency in the Gaussian location model, see, e.g., \citealp[][Chapter 2]{Maronna:2019}), and
        \item[\textsf{L}:] the logistic estimator with $\rho(x) = 2x + 4\log(1+e^{-x})$.
    \end{itemize}  
    % For the Huber estimator in~\eqref{eq:huber} we use $k= 1.345$, which ensures $95\%$ efficiency in the Gaussian location model \citep[see, e.g.,][Chapter 2]{Maronna:2019}. 
    These methods are compared with two well-established competitors:
    \begin{itemize}
        \item[\textsf{P}:] the \textsf{pace} estimator based on kernel smoothing of sparse functional data \citep{Yao_etal2005}, as implemented in the \textsf{R} package \textsf{fdapace} \citep[][available only for $d=1$]{R:fdapace}, and
        \item[\textsf{NP}:] the local linear estimator from the \textsf{R} package \textsf{np} \citep{R:np}.
        The bandwidth tuning parameter is selected by cross-validation as proposed in \citet{Li_Racine2004}. 
    \end{itemize} 
    In addition, we consider two naive estimators of $\mu_0$ based on the piecewise linear interpolation of the observed functional data:
    \begin{itemize}
        \item[\textsf{N$_1$}:] the sample mean of the piecewise linear interpolants, and
        \item[\textsf{N$_2$}:] the coordinate-wise median function of the piecewise linear interpolating functions.
    \end{itemize}
    For $d=2$, the piecewise linear interpolation on $\mathcal{I} = (0,1)^2$ in \textsf{N$_1$} and \textsf{N$_2$} is obtained by function \textsf{interp} from the \textsf{R} package of the same name \citep{R:interp}.\footnote{Outside the convex hull of the observed time points, the piecewise linear interpolation is obtained by including additional artificial function values at the four corners $(0,0)^\top$, $(0,1)^\top$, $(1,0)^\top$, and $(1,1)^\top$, with values equal to the closest observed functional values.}

\begin{figure}[htpb]
    \centering
    \includegraphics[width = 0.975\textwidth]{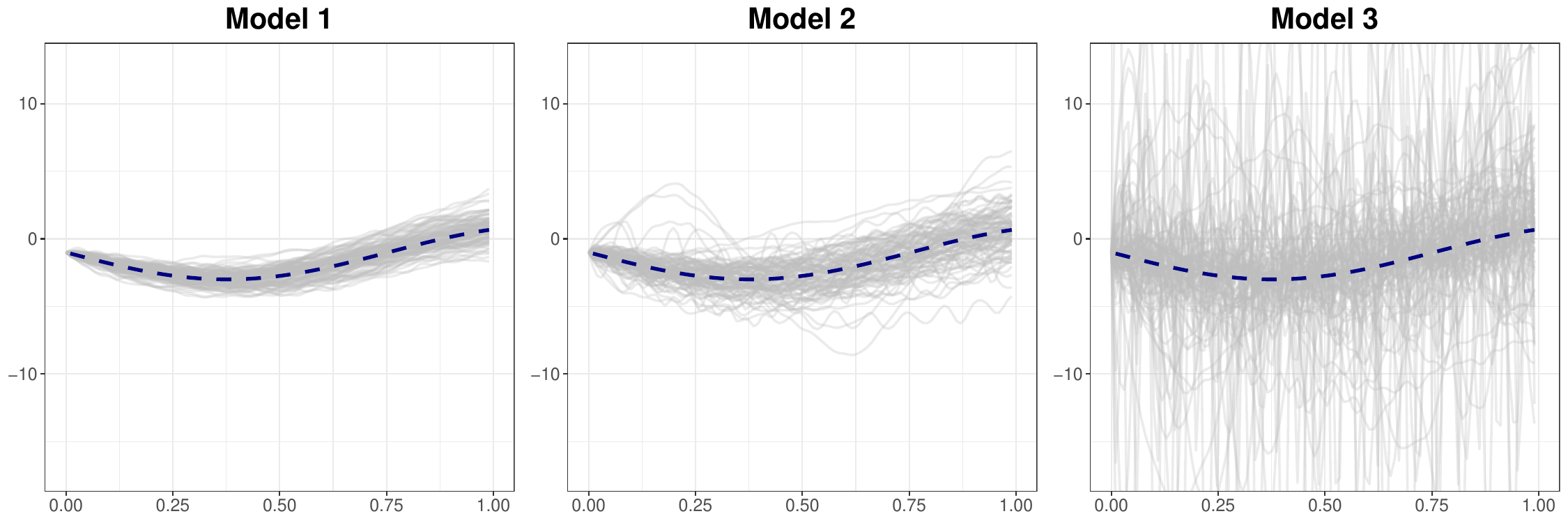}
    \caption{$100$ curves generated with standard Gaussian, $t_2$ and $t_1$ distributed $\{W_{ij}\}_{i = 1,j=1,1}^{100,50}$ on the left, central and right panels, respectively.}
    \label{fig:model_curves}
\end{figure}

We explore the sensitivity of all these functional location estimators to the sample size $n$, the grid size $p$, and the pattern of irregularity in the observation points in $\mathcal{I}$ with varying distributions of $X$. We generate a functional sample $X_1, \dots, X_n$ of the form $X_i(\mathbf{t}) = \mu_0(\mathbf{t}) + \epsilon_i(\mathbf{t})$ for $\mathbf{t} \in \mathcal{I} = (0,1)^d$ for $d = 1, 2$. Here, $\mu_{0}$ is the target function to be estimated, and $\epsilon_1, \dots, \epsilon_n$ are independent, zero-centered noise functions given by their Karhunen-Lo\`eve expansion. In the one-dimensional setting ($d=1$), we consider 
    \[
    \begin{aligned}
    % \label{eq:mu_epsilon_1D}
    \mu_0(t) = - 2 \sin(5\,t/1.2) - 1 % \quad \mbox{for }t \in \mathcal{I} = (0,1); \\
    %\quad \quad \quad
    \quad \mbox{and} \quad \epsilon_i(t) = \sqrt{2} \sum_{j=1}^{50} W_{ij} \frac{\sin \left( \left(j-1/2 \right) \pi t \right)}{\left(j-1/2\right)\pi}, \quad \mbox{for $t \in (0,1)$,}
    \end{aligned}
    \]
while for $d=2$ we use the model
    \[
    \begin{aligned}
    \mu_0(\mathbf{t}) &= -\sin(5\,t_{1}/1.2)/0.5 + 2\,t_{2} -1 \quad \mbox{for } \mathbf{t}=(t_{1}, t_{2})^\top \in \mathcal{I} = (0,1)^2;
    \\
    \epsilon_i(\mathbf{t}) &= 2 \sum_{j=1}^{20} \sum_{k=1}^{20} W_{ijk} \frac{\sin \left( \left(j-1/2 \right) \pi t_{1} \right)}{\left(j-1/2\right)\pi} \frac{\sin \left( \left(k-1/2 \right) \pi t_{2} \right)}{\left(k-1/2\right)\pi}, \quad \mbox{for } i = 1, \ldots, n.
    \end{aligned}
    \]
In both setups, the coefficients $W_{ij}$ (resp. $W_{ijk}$) are independent, distributed as: % zero-centered random variables. We test the estimators in the following three scenarios for the noise processes:
\begin{itemize}[leftmargin=2cm, labelwidth=2cm, labelsep=0.5cm]
    \item[\textbf{Model$^d$~1}] $W_{ij}$ follows the standard Gaussian distribution $\mathsf{N}(0,1)$;
    \item[\textbf{Model$^d$~2}] $W_{ij}$ follows the t-distribution with two degrees of freedom $t_2$;
    \item[\textbf{Model$^d$~3}] $W_{ij}$ has the t-distribution with one degree of freedom $t_1$ (the Cauchy distribution).
\end{itemize}
These scenarios reflect the settings of clean, well-behaved data (\textbf{Model$^{d}$~1}), data contaminated with heavy-tailed, yet integrable noise (\textbf{Model$^{d}$~2}), and data with non-integrable errors (\textbf{Model$^{d}$~3}). Figure~\ref{fig:model_curves} shows $n = 100$ sample curves generated by the three models for $d=1$; for two-dimensional functions ($d=2$), Figure~\ref{fig:d2} in the Supplementary Material shows a similar display. % The heavier the tails of the noise process, the more complicated. We stress that switching from the Gaussian to a t-distribution leads to several outlying measurements among the functional data, which in turn complicate the estimation of $\mu_0$. That effect becomes even more severe when the errors fail to have a finite expectation (\textbf{Model~3}).

The complete functional data $\{ X_{i} \}_{i=1}^{n}$ are represented in a regular grid: for $d=1$ at equidistant points $t_1 < t_2 < \dots < t_{100}$ in the interval $(0,1)$, and for $d=2$ in a regular lattice grid of $30^2 = 900$ points $\left\{ \mathbf{t}_j \right\}_{j=1}^{900} = \left\{ (t_h, t_\ell)^\top \right\}_{h = 1, \ell = 1}^{30, 30} \subset (0, 1)^2$. Two incomplete data settings are considered:
    \begin{itemize}
        \item[] \textbf{Balanced design} (for $d=1$), where all $X_1, \dots, X_n$ are observed on a common grid of randomly chosen $p$ points $0<t_{k_1}<\ldots<t_{k_p} <1$ from $t_1, \ldots, t_{100}$.
        %compute the estimators of $\mu_0$ from the incomplete data $(X_i(t_{k_1}), \ldots, X_i(t_{k_p}))$, $i =1, \ldots, n$.  
        \item[] \textbf{Unbalanced design} (for $d=1, 2$), where for each function $X_i$ we randomly choose $p$ points from the complete lattice, with the measurement points varying among statistical units. %and compute the estimators of $\mu_0$ from the incomplete data $(X_i(t_{i k_1}), \ldots, X_i(t_{i k_p}))$, $i =1, \ldots, n$. In this setup, each function is observed on a different grid of points.
    \end{itemize}  
For $d=1$ we consider sample sizes $n \in \{25, 50, 100\}$ and $p \in \{ 10, 25, 50 \}$; for $d=2$ we have $n \in \{ 10, 25, 50 \}$ and $p \in \{10, 50, 100, 200\}$. The resolution $p = 50$ for $d=1$ (resp. $p = 200$ for $ d = 2$) corresponds to moderately dense functional data. The performance of the estimators is assessed by the squared estimation error %SEE of $\widehat{\mu}_n$
    %\begin{equation} \label{eq: SEE}
    $\SEE =  B^{-1} \sum_{j=1}^{B} \left|\mu_0(\mathbf{t}_{j}) - \widehat\mu_n(\mathbf{t}_{j}) \right|^2$, for $B = 100$ with $d=1$ and $B = 900$ for $d=2$. 
    %\end{equation}
SEE approximates the $L_2$-distance between $\mu_0$ and its estimator $\widehat{\mu}_n$.

\begin{table}[htpb]
\centering
\caption{One-dimensional case, means and standard errors (in brackets) of the squared estimation errors ($\times 100$) for \textbf{Models$^1$ 1--3} in the balanced and unbalanced design; 1000 runs of the simulation study. Values larger than 1000 are truncated. The smallest mean in each setup is in bold.}
\label{Tab:results_model_d1}
\resizebox{\textwidth}{!}{\begin{tabular}{ccc|ccc|ccc|ccc}
\hline
 & & & \multicolumn{3}{c}{$n=25$} & \multicolumn{3}{c}{$n=50$} & \multicolumn{3}{c}{$n=100$} \\
 & & & $p=10$  & $p=25$  & $p=50$  & $p=10$  & $p=25$  & $p=50$  & $p=10$  & $p=25$  & $p=50$ \\
\hline
\parbox[t]{2mm}{\multirow{24}{*}{\rotatebox[origin=c]{90}{\textbf{Balanced}}}} & \parbox[t]{2mm}{\multirow{8}{*}{\rotatebox[origin=c]{90}{\textbf{Model$^{1}$~1}}}} & \textsf{S} & \textbf{6.986 }  \scriptsize{(0.674)} & \textbf{2.121}  \scriptsize{(0.074)} & \textbf{2.104}  \scriptsize{(0.080)} & \textbf{3.639}  \scriptsize{(0.410)} & \textbf{1.187}  \scriptsize{(0.044)} & \textbf{1.036}  \scriptsize{(0.037)} & \textbf{3.349}  \scriptsize{(0.364)} & \textbf{0.571}  \scriptsize{(0.018)} & \textbf{0.479 } \scriptsize{(0.017)} \\
 &  & \textsf{A} & 11.557  \scriptsize{(0.859)} & 3.678  \scriptsize{(0.113)} & 3.200  \scriptsize{(0.113)} & 5.842  \scriptsize{(0.460)} & 2.195  \scriptsize{(0.067)} & 1.736  \scriptsize{(0.053)} & 4.417  \scriptsize{(0.412)} & 1.244  \scriptsize{(0.033)} & 0.970  \scriptsize{(0.026)} \\
 &  & \textsf{H} & 7.063  \scriptsize{(0.689)} & 2.191  \scriptsize{(0.075)} & 2.172  \scriptsize{(0.084)} & 3.674  \scriptsize{(0.410)} & 1.213  \scriptsize{(0.045)} & 1.071  \scriptsize{(0.038)} & 3.518  \scriptsize{(0.375)} & 0.602  \scriptsize{(0.019)} & 0.499  \scriptsize{(0.018)} \\
 &  & \textsf{L} & 7.036  \scriptsize{(0.680)} & 2.143  \scriptsize{(0.074)} & 2.129  \scriptsize{(0.082)} & 3.654  \scriptsize{(0.410)} & 1.191  \scriptsize{(0.044)} & 1.050  \scriptsize{(0.038)} & 3.383  \scriptsize{(0.357)} & 0.587  \scriptsize{(0.018)} & 0.491  \scriptsize{(0.017)} \\
 &  & \textsf{P} & 11.607  \scriptsize{(0.553)} & 2.661  \scriptsize{(0.100)} & 2.135  \scriptsize{(0.079)} & 9.555  \scriptsize{(0.536)} & 1.721  \scriptsize{(0.071)} & 1.071  \scriptsize{(0.038)} & 10.038  \scriptsize{(0.563)} & 1.077  \scriptsize{(0.043)} & 0.501  \scriptsize{(0.017)} \\
 &  & \textsf{NP} & 6.991  \scriptsize{(0.554)} & 2.220  \scriptsize{(0.075)} & 2.175  \scriptsize{(0.081)} & 4.323  \scriptsize{(0.373)} & 1.278  \scriptsize{(0.045)} & 1.084  \scriptsize{(0.038)} & 3.604  \scriptsize{(0.272)} & 0.619  \scriptsize{(0.018)} & 0.506  \scriptsize{(0.018)} \\
 &  & \textsf{N$_1$} & 12.518  \scriptsize{(0.544)} & 2.666  \scriptsize{(0.101)} & 2.133  \scriptsize{(0.079)} & 10.713  \scriptsize{(0.528)} & 1.753  \scriptsize{(0.071)} & 1.070  \scriptsize{(0.038)} & 10.984  \scriptsize{(0.556)} & 1.106  \scriptsize{(0.043)} & 0.501  \scriptsize{(0.017)} \\
 &  & \textsf{N$_2$} & 13.469  \scriptsize{(0.548)} & 3.724  \scriptsize{(0.118)} & 3.280  \scriptsize{(0.103)} & 11.189  \scriptsize{(0.535)} & 2.269  \scriptsize{(0.076)} & 1.605  \scriptsize{(0.048)} & 11.249  \scriptsize{(0.567)} & 1.381  \scriptsize{(0.047)} & 0.772  \scriptsize{(0.022)} \\ \cdashline{2-12}
 & \parbox[t]{2mm}{\multirow{8}{*}{\rotatebox[origin=c]{90}{\textbf{Model$^{1}$ 2}}}} & \textsf{S} & 35.030  \scriptsize{(6.809)} & 18.184  \scriptsize{(1.808)} & 40.549  \scriptsize{(19.728)} & 13.812  \scriptsize{(0.910)} & 21.397  \scriptsize{(6.863)} & 8.213  \scriptsize{(0.568)} & 11.625  \scriptsize{(1.465)} & 5.720  \scriptsize{(0.726)} & 5.308  \scriptsize{(0.526)} \\
 &  & \textsf{A} & 17.984  \scriptsize{(1.044)} & 7.001  \scriptsize{(0.226)} & 6.431  \scriptsize{(0.197)} & 10.365  \scriptsize{(0.709)} & 4.022  \scriptsize{(0.119)} & 3.589  \scriptsize{(0.104)} & 7.715  \scriptsize{(0.843)} & 2.251  \scriptsize{(0.062)} & 1.959  \scriptsize{(0.053)} \\
 &  & \textsf{H} & \textbf{13.130} \scriptsize{(0.880)} & \textbf{5.634} \scriptsize{(0.194)} & \textbf{5.156} \scriptsize{(0.173)} & 8.406  \scriptsize{(0.780)} & \textbf{2.881}  \scriptsize{(0.086)} & \textbf{2.513}  \scriptsize{(0.081)} & \textbf{6.122}  \scriptsize{(0.780)} & \textbf{1.532}  \scriptsize{(0.048)} & \textbf{1.335}  \scriptsize{(0.041)} \\
 &  & \textsf{L} & 14.185  \scriptsize{(0.989)} & 6.240  \scriptsize{(0.224)} & 5.736  \scriptsize{(0.196)} & \textbf{8.072}  \scriptsize{(0.658)} & 3.181  \scriptsize{(0.098)} & 2.833  \scriptsize{(0.092)} & 6.315  \scriptsize{(0.786)} & 1.644  \scriptsize{(0.050)} & 1.500  \scriptsize{(0.048)} \\
 &  & \textsf{P} & 37.181  \scriptsize{(6.643)} & 19.936  \scriptsize{(1.858)} & 41.986  \scriptsize{(19.758)} & 20.441  \scriptsize{(1.298)} & 22.304  \scriptsize{(6.775)} & 8.786  \scriptsize{(0.581)} & 14.968  \scriptsize{(1.151)} & 6.378  \scriptsize{(0.711)} & 5.602  \scriptsize{(0.523)} \\
 &  & \textsf{NP} & 38.586  \scriptsize{(6.946)} & 18.511  \scriptsize{(1.700)} & 40.444  \scriptsize{(19.179)} & 14.094  \scriptsize{(0.929)} & 22.316  \scriptsize{(7.287)} & 8.637  \scriptsize{(0.620)} & 12.695  \scriptsize{(1.521)} & 6.263  \scriptsize{(0.777)} & 5.589  \scriptsize{(0.593)} \\
 &  & \textsf{N$_1$} & 35.750  \scriptsize{(6.622)} & 19.444  \scriptsize{(1.852)} & 41.914  \scriptsize{(19.753)} & 18.852  \scriptsize{(0.845)} & 22.087  \scriptsize{(6.769)} & 8.736  \scriptsize{(0.580)} & 15.440  \scriptsize{(1.142)} & 6.319  \scriptsize{(0.709)} & 5.580  \scriptsize{(0.523)} \\
 &  & \textsf{N$_2$} & 16.140  \scriptsize{(0.631)} & 6.930  \scriptsize{(0.174)} & 6.481  \scriptsize{(0.167)} & 13.657  \scriptsize{(0.656)} & 3.663  \scriptsize{(0.091)} & 3.104  \scriptsize{(0.076)} & 10.696  \scriptsize{(0.528)} & 2.128  \scriptsize{(0.058)} & 1.617  \scriptsize{(0.038)} \\ \cdashline{2-12}
 & \parbox[t]{2mm}{\multirow{8}{*}{\rotatebox[origin=c]{90}{\textbf{Model$^{1}$~3}}}} & \textsf{S} & $>1000$ & $>1000$ & $>1000$ & $>1000$ & $>1000$ & $>1000$ & $>1000$ & $>1000$ & $>1000$ \\
 &  & \textsf{A} & 64.567  \scriptsize{(3.741)} & 32.840  \scriptsize{(4.837)} & 22.255  \scriptsize{(3.260)} & 42.932  \scriptsize{(2.536)} & 13.666  \scriptsize{(0.419)} & 11.593  \scriptsize{(0.929)} & 24.969  \scriptsize{(1.358)} & 8.207  \scriptsize{(0.245)} & 5.836  \scriptsize{(0.158)} \\
 &  & \textsf{H} & 39.427 \scriptsize{(1.572)} & \textbf{20.873} \scriptsize{(0.662)} & \textbf{17.194}  \scriptsize{(0.522)} & 25.201  \scriptsize{(1.276)} & \textbf{10.089} \scriptsize{(0.288)} & \textbf{8.759}  \scriptsize{(0.258)} & 17.326  \scriptsize{(1.183)} & \textbf{5.305}\scriptsize{(0.155)} & \textbf{4.436}  \scriptsize{(0.126)} \\
 &  & \textsf{L} & 43.330  \scriptsize{(1.648)} & 24.537  \scriptsize{(0.769)} & 20.431  \scriptsize{(0.636)} & 26.978  \scriptsize{(1.344)} & 11.706  \scriptsize{(0.329)} & 10.267  \scriptsize{(0.298)} & 18.102  \scriptsize{(1.230)} & 6.011  \scriptsize{(0.166)} & 5.302  \scriptsize{(0.151)} \\
 &  & \textsf{P} & $>1000$ & $>1000$ & $>1000$ & $>1000$ & $>1000$ & $>1000$ & $>1000$ & $>1000$ & $>1000$ \\
 &  & \textsf{NP} & $>1000$ & $>1000$ & $>1000$ & $>1000$ & $>1000$ & $>1000$ & $>1000$ & $>1000$ & $>1000$ \\
 &  & \textsf{N$_1$} & $>1000$ & $>1000$ & $>1000$ & $>1000$ & $>1000$ & $>1000$ & $>1000$ & $>1000$ & $>1000$ \\
 &  & \textsf{N$_2$} & \textbf{37.519}  \scriptsize{(1.050)} & 28.627  \scriptsize{(0.650)} & 28.873  \scriptsize{(0.553)} & \textbf{21.536} \scriptsize{(0.686)} & 13.216  \scriptsize{(0.246)} & 13.826  \scriptsize{(0.245)} & \textbf{15.566} \scriptsize{(0.618)} & 6.765  \scriptsize{(0.125)} & 6.779  \scriptsize{(0.118)} \\ \hline
\parbox[t]{2mm}{\multirow{24}{*}{\rotatebox[origin=c]{90}{\textbf{Unbalanced}}}} & \parbox[t]{2mm}{\multirow{8}{*}{\rotatebox[origin=c]{90}{\textbf{Model$^{1}$~1}}}} & \textsf{S} & 2.612  \scriptsize{(0.075)} & 2.224  \scriptsize{(0.079)} & 2.169  \scriptsize{(0.078)} & 1.447  \scriptsize{(0.042)} & 1.156  \scriptsize{(0.038)} & 1.074  \scriptsize{(0.038)} & \textbf{0.717}  \scriptsize{(0.021)} & \textbf{0.575}  \scriptsize{(0.018)} & 0.549  \scriptsize{(0.020)} \\
 &  & \textsf{A} & 3.421  \scriptsize{(0.100)} & 3.183  \scriptsize{(0.111)} & 3.032  \scriptsize{(0.105)} & 2.025  \scriptsize{(0.058)} & 1.665  \scriptsize{(0.049)} & 1.586  \scriptsize{(0.049)} & 1.177  \scriptsize{(0.032)} & 0.953  \scriptsize{(0.027)} & 0.899  \scriptsize{(0.027)} \\
 &  & \textsf{H} & 2.619  \scriptsize{(0.075)} & 2.289  \scriptsize{(0.083)} & 2.210  \scriptsize{(0.080)} & 1.456  \scriptsize{(0.043)} & 1.174  \scriptsize{(0.038)} & 1.103  \scriptsize{(0.039)} & 0.737  \scriptsize{(0.022)} & 0.598  \scriptsize{(0.020)} & 0.572  \scriptsize{(0.020)} \\
 &  & \textsf{L} & \textbf{2.578} \scriptsize{(0.075)} & 2.243  \scriptsize{(0.081)} & 2.176  \scriptsize{(0.079)} & \textbf{1.433}  \scriptsize{(0.043)} & \textbf{1.155}  \scriptsize{(0.038)} & 1.084  \scriptsize{(0.039)} & 0.718  \scriptsize{(0.022)} & 0.583  \scriptsize{(0.019)} & 0.561  \scriptsize{(0.020)} \\
 &  & \textsf{P} & 2.988  \scriptsize{(0.079)} & 2.340  \scriptsize{(0.080)} & 2.214  \scriptsize{(0.079)} & 1.628  \scriptsize{(0.043)} & 1.220  \scriptsize{(0.039)} & 1.112  \scriptsize{(0.039)} & 0.808  \scriptsize{(0.022)} & 0.638  \scriptsize{(0.020)} & 0.598  \scriptsize{(0.020)} \\
 &  & \textsf{NP} & 3.035  \scriptsize{(0.081)} & 2.372  \scriptsize{(0.080)} & 2.232  \scriptsize{(0.079)} & 1.701  \scriptsize{(0.045)} & 1.243  \scriptsize{(0.039)} & 1.109  \scriptsize{(0.039)} & 0.842  \scriptsize{(0.022)} & 0.631  \scriptsize{(0.019)} & 0.572  \scriptsize{(0.020)} \\
 &  & \textsf{N$_1$} & 5.386  \scriptsize{(0.087)} & \textbf{2.161}  \scriptsize{(0.078)} & \textbf{2.106}  \scriptsize{(0.077)} & 4.460  \scriptsize{(0.056)} & 1.210  \scriptsize{(0.038)} & \textbf{1.048}  \scriptsize{(0.038)} & 3.935  \scriptsize{(0.037)} & 0.677  \scriptsize{(0.018)} & \textbf{0.534}  \scriptsize{(0.020)} \\
 &  & \textsf{N$_2$} & 5.431  \scriptsize{(0.101)} & 3.228  \scriptsize{(0.100)} & 3.217  \scriptsize{(0.103)} & 4.034  \scriptsize{(0.065)} & 1.621  \scriptsize{(0.047)} & 1.575  \scriptsize{(0.047)} & 3.217  \scriptsize{(0.041)} & 0.866  \scriptsize{(0.024)} & 0.808  \scriptsize{(0.026)} \\ \cdashline{2-12}
 & \parbox[t]{2mm}{\multirow{8}{*}{\rotatebox[origin=c]{90}{\textbf{Model$^{1}$~2}}}} & \textsf{S} & 22.849  \scriptsize{(2.737)} & 23.532  \scriptsize{(4.197)} & 15.356  \scriptsize{(1.279)} & 22.124  \scriptsize{(8.567)} & 13.537  \scriptsize{(3.623)} & 9.194  \scriptsize{(0.802)} & 6.517  \scriptsize{(0.573)} & 6.503  \scriptsize{(0.653)} & 4.995  \scriptsize{(0.344)} \\
 &  & \textsf{A} & 7.846  \scriptsize{(0.251)} & 6.686  \scriptsize{(0.212)} & 6.075  \scriptsize{(0.203)} & 4.499  \scriptsize{(0.124)} & 3.470  \scriptsize{(0.101)} & 3.533  \scriptsize{(0.102)} & 2.656  \scriptsize{(0.071)} & 2.216  \scriptsize{(0.059)} & 1.908  \scriptsize{(0.051)} \\
 &  & \textsf{H} & \textbf{6.899} \scriptsize{(0.202)} & \textbf{5.613}  \scriptsize{(0.177)} & \textbf{4.944}  \scriptsize{(0.170)} & \textbf{3.375}  \scriptsize{(0.092)} & \textbf{2.770}  \scriptsize{(0.080)} & \textbf{2.471}  \scriptsize{(0.076)} & \textbf{1.751}  \scriptsize{(0.045)} & \textbf{1.490}  \scriptsize{(0.046)} & \textbf{1.313}  \scriptsize{(0.043)} \\
 &  & \textsf{L} & 7.606  \scriptsize{(0.228)} & 6.274  \scriptsize{(0.201)} & 5.536  \scriptsize{(0.194)} & 3.704  \scriptsize{(0.103)} & 3.086  \scriptsize{(0.093)} & 2.757  \scriptsize{(0.088)} & 1.909  \scriptsize{(0.049)} & 1.682  \scriptsize{(0.054)} & 1.500  \scriptsize{(0.051)} \\
 &  & \textsf{P} & 30.011  \scriptsize{(4.593)} & 25.297  \scriptsize{(4.289)} & 16.091  \scriptsize{(1.303)} & 27.862  \scriptsize{(10.600)} & 14.815  \scriptsize{(3.984)} & 9.750  \scriptsize{(0.841)} & 8.040  \scriptsize{(0.738)} & 7.065  \scriptsize{(0.708)} & 5.213  \scriptsize{(0.354)} \\
 &  & \textsf{NP} & 24.755  \scriptsize{(2.288)} & 23.851  \scriptsize{(3.899)} & 15.767  \scriptsize{(1.281)} & 24.325  \scriptsize{(8.309)} & 14.120  \scriptsize{(3.617)} & 9.538  \scriptsize{(0.853)} & 7.691  \scriptsize{(0.746)} & 6.980  \scriptsize{(0.661)} & 5.201  \scriptsize{(0.368)} \\
 &  & \textsf{N$_1$} & 23.568  \scriptsize{(2.689)} & 23.028  \scriptsize{(3.960)} & 16.196  \scriptsize{(1.300)} & 19.292  \scriptsize{(4.618)} & 12.780  \scriptsize{(3.031)} & 9.671  \scriptsize{(0.868)} & 8.466  \scriptsize{(0.542)} & 6.103  \scriptsize{(0.586)} & 5.141  \scriptsize{(0.370)} \\
 &  & \textsf{N$_2$} & 8.907  \scriptsize{(0.206)} & 6.336  \scriptsize{(0.177)} & 6.057  \scriptsize{(0.156)} & 5.539  \scriptsize{(0.101)} & 2.919  \scriptsize{(0.068)} & 3.007  \scriptsize{(0.075)} & 3.971  \scriptsize{(0.059)} & 1.647  \scriptsize{(0.039)} & 1.511  \scriptsize{(0.038)} \\ \cdashline{2-12}
 & \parbox[t]{2mm}{\multirow{8}{*}{\rotatebox[origin=c]{90}{\textbf{Model$^{1}$~3}}}} & \textsf{S} & $>1000$ & $>1000$ & $>1000$ & $>1000$ & $>1000$ & $>1000$ & $>1000$ & $>1000$ & $>1000$ \\
 &  & \textsf{A} & 35.324  \scriptsize{(2.828)} & 20.878  \scriptsize{(0.705)} & 21.301  \scriptsize{(2.555)} & 17.690  \scriptsize{(2.063)} & 12.434  \scriptsize{(0.417)} & 10.227  \scriptsize{(0.335)} & 10.012  \scriptsize{(0.310)} & 6.953  \scriptsize{(0.181)} & 5.748  \scriptsize{(0.167)} \\
 &  & \textsf{H} & \textbf{28.904} \scriptsize{(0.825)} & \textbf{20.086}  \scriptsize{(0.624)} & \textbf{17.804}  \scriptsize{(0.574)} & \textbf{14.099}  \scriptsize{(0.395)} & \textbf{10.645}  \scriptsize{(0.313)} & \textbf{9.350}  \scriptsize{(0.284)} & \textbf{7.659}  \scriptsize{(0.210)} & \textbf{5.355}  \scriptsize{(0.144)} & \textbf{4.630}  \scriptsize{(0.143)} \\
 &  & \textsf{L} & 33.055  \scriptsize{(0.937)} & 23.276  \scriptsize{(0.720)} & 20.557  \scriptsize{(0.657)} & 16.095  \scriptsize{(0.416)} & 12.482  \scriptsize{(0.356)} & 10.749  \scriptsize{(0.321)} & 8.780  \scriptsize{(0.235)} & 6.269  \scriptsize{(0.173)} & 5.353  \scriptsize{(0.161)} \\
 &  & \textsf{P} & $>1000$ & $>1000$ & $>1000$ & $>1000$ & $>1000$ & $>1000$ & $>1000$ & $>1000$ & $>1000$ \\
 &  & \textsf{NP} & $>1000$ & $>1000$ & $>1000$ & $>1000$ & $>1000$ & $>1000$ & $>1000$ & $>1000$ & $>1000$ \\
 &  & \textsf{N$_1$} & $>1000$ & $>1000$ & $>1000$ & $>1000$ & $>1000$ & $>1000$ & $>1000$ & $>1000$ & $>1000$ \\
 &  & \textsf{N$_2$} & 29.754  \scriptsize{(0.707)} & 26.409  \scriptsize{(0.600)} & 29.575  \scriptsize{(0.627)} & 15.073  \scriptsize{(0.288)} & 12.995  \scriptsize{(0.273)} & 13.917  \scriptsize{(0.274)} & 8.817  \scriptsize{(0.161)} & 6.222  \scriptsize{(0.106)} & 6.590  \scriptsize{(0.119)} \\ \hline
\end{tabular}}
\end{table}

% latex table generated in R 4.5.2 by xtable 1.8-4 package
% Sun Oct 11 17:18:15 2026
\begin{table}[htpb]
\centering
\caption{Two-dimensional case, means and standard errors (in brackets) of the squared estimation errors ($\times 100$) for \textbf{Models$^2$ 1--3}; 500 runs of the simulation study. Values exceeding 1000 are truncated. The mean of the best-performing method in each setup is in bold.} 
\label{Tab:model_d2_regular}
\resizebox{\textwidth}{!}{\begin{tabular}{cc|cccc|cccc|cccc}
   \hline
 & & \multicolumn{4}{c}{$n=10$} & \multicolumn{4}{c}{$n=25$} & \multicolumn{4}{c}{$n=50$} \\
 & & $p=10$  & $p=50$  & $p=100$  & $p=200$  & $p=10$  & $p=50$  & $p=100$  & $p=200$  & $p=10$  & $p=50$  & $p=100$  & $p=200$ \\
   \hline
\parbox[t]{2mm}{\multirow{7}{*}{\rotatebox[origin=c]{90}{\hspace{-2.5cm}\textbf{Model$^2$ 1}}}} & \textsf{S}     & \shortstack{\rule{0pt}{3ex}\textbf{36.535} \\[-0.2ex] {\tiny(1.830)}\rule[-1.5ex]{0pt}{1.0ex}}    & \shortstack{\rule{0pt}{3ex}3.516 \\[-0.2ex] {\tiny(0.127)}\rule[-1.0ex]{0pt}{1.0ex}}     & \shortstack{\rule{0pt}{3ex}2.450 \\[-0.2ex] {\tiny(0.098)}\rule[-1.0ex]{0pt}{1.0ex}}     & \shortstack{\rule{0pt}{3ex}2.464 \\[-0.2ex] {\tiny(0.107)}\rule[-1.0ex]{0pt}{1.0ex}}     & \shortstack{\rule{0pt}{3ex}\textbf{31.757} \\[-0.2ex] {\tiny(1.506)}\rule[-1.0ex]{0pt}{1.0ex}}    & \shortstack{\rule{0pt}{3ex}1.934 \\[-0.2ex] {\tiny(0.072)}\rule[-1.0ex]{0pt}{1.0ex}}     & \shortstack{\rule{0pt}{3ex}1.216 \\[-0.2ex] {\tiny(0.045)}\rule[-1.0ex]{0pt}{1.0ex}}     & \shortstack{\rule{0pt}{3ex}0.942 \\[-0.2ex] {\tiny(0.043)}\rule[-1.0ex]{0pt}{1.0ex}}     & \shortstack{\rule{0pt}{3ex}\textbf{31.582} \\[-0.2ex] {\tiny(1.497)}\rule[-1.0ex]{0pt}{1.0ex}}    & \shortstack{\rule{0pt}{3ex}1.257 \\[-0.2ex] {\tiny(0.039)}\rule[-1.0ex]{0pt}{1.0ex}}     & \shortstack{\rule{0pt}{3ex}0.657 \\[-0.2ex] {\tiny(0.022)}\rule[-1.0ex]{0pt}{1.0ex}}     & \shortstack{\rule{0pt}{3ex}0.522 \\[-0.2ex] {\tiny(0.020)}\rule[-1.0ex]{0pt}{1.0ex}}     \\ 
                                                                                   & \textsf{A}     & \shortstack{\rule{0pt}{3ex}46.314 \\[-0.2ex] {\tiny(2.175)}\rule[-1.0ex]{0pt}{1.0ex}}    & \shortstack{\rule{0pt}{3ex}5.582 \\[-0.2ex] {\tiny(0.177)}\rule[-1.0ex]{0pt}{1.0ex}}     & \shortstack{\rule{0pt}{3ex}3.716 \\[-0.2ex] {\tiny(0.134)}\rule[-1.0ex]{0pt}{1.0ex}}     & \shortstack{\rule{0pt}{3ex}3.204 \\[-0.2ex] {\tiny(0.127)}\rule[-1.0ex]{0pt}{1.0ex}}     & \shortstack{\rule{0pt}{3ex}40.119 \\[-0.2ex] {\tiny(1.715)}\rule[-1.0ex]{0pt}{1.0ex}}    & \shortstack{\rule{0pt}{3ex}3.737 \\[-0.2ex] {\tiny(0.114)}\rule[-1.0ex]{0pt}{1.0ex}}     & \shortstack{\rule{0pt}{3ex}2.153 \\[-0.2ex] {\tiny(0.057)}\rule[-1.0ex]{0pt}{1.0ex}}     & \shortstack{\rule{0pt}{3ex}1.453 \\[-0.2ex] {\tiny(0.051)}\rule[-1.0ex]{0pt}{1.0ex}}     & \shortstack{\rule{0pt}{3ex}32.090 \\[-0.2ex] {\tiny(1.503)}\rule[-1.0ex]{0pt}{1.0ex}}    & \shortstack{\rule{0pt}{3ex}3.049 \\[-0.2ex] {\tiny(0.077)}\rule[-1.0ex]{0pt}{1.0ex}}     & \shortstack{\rule{0pt}{3ex}1.476 \\[-0.2ex] {\tiny(0.031)}\rule[-1.0ex]{0pt}{1.0ex}}     & \shortstack{\rule{0pt}{3ex}0.907 \\[-0.2ex] {\tiny(0.052)}\rule[-1.0ex]{0pt}{1.0ex}}     \\ 
                                                                                   & \textsf{H}     & \shortstack{\rule{0pt}{3ex}35.899 \\[-0.2ex] {\tiny(1.811)}\rule[-1.0ex]{0pt}{1.0ex}}    & \shortstack{\rule{0pt}{3ex}3.578 \\[-0.2ex] {\tiny(0.128)}\rule[-1.0ex]{0pt}{1.0ex}}     & \shortstack{\rule{0pt}{3ex}2.513 \\[-0.2ex] {\tiny(0.100)}\rule[-1.0ex]{0pt}{1.0ex}}     & \shortstack{\rule{0pt}{3ex}2.507 \\[-0.2ex] {\tiny(0.107)}\rule[-1.0ex]{0pt}{1.0ex}}     & \shortstack{\rule{0pt}{3ex}31.770 \\[-0.2ex] {\tiny(1.491)}\rule[-1.0ex]{0pt}{1.0ex}}    & \shortstack{\rule{0pt}{3ex}1.981 \\[-0.2ex] {\tiny(0.072)}\rule[-1.0ex]{0pt}{1.0ex}}     & \shortstack{\rule{0pt}{3ex}1.265 \\[-0.2ex] {\tiny(0.046)}\rule[-1.0ex]{0pt}{1.0ex}}     & \shortstack{\rule{0pt}{3ex}0.976 \\[-0.2ex] {\tiny(0.043)}\rule[-1.0ex]{0pt}{1.0ex}}     & \shortstack{\rule{0pt}{3ex}31.894 \\[-0.2ex] {\tiny(1.485)}\rule[-1.0ex]{0pt}{1.0ex}}    & \shortstack{\rule{0pt}{3ex}1.340 \\[-0.2ex] {\tiny(0.040)}\rule[-1.0ex]{0pt}{1.0ex}}     & \shortstack{\rule{0pt}{3ex}0.712 \\[-0.2ex] {\tiny(0.022)}\rule[-1.0ex]{0pt}{1.0ex}}     & \shortstack{\rule{0pt}{3ex}0.573 \\[-0.2ex] {\tiny(0.020)}\rule[-1.0ex]{0pt}{1.0ex}}     \\ 
                                                                                   & \textsf{L}     & \shortstack{\rule{0pt}{3ex}36.679 \\[-0.2ex] {\tiny(1.832)}\rule[-1.0ex]{0pt}{1.0ex}}    & \shortstack{\rule{0pt}{3ex}3.526 \\[-0.2ex] {\tiny(0.128)}\rule[-1.0ex]{0pt}{1.0ex}}     & \shortstack{\rule{0pt}{3ex}2.453 \\[-0.2ex] {\tiny(0.099)}\rule[-1.0ex]{0pt}{1.0ex}}     & \shortstack{\rule{0pt}{3ex}2.475 \\[-0.2ex] {\tiny(0.107)}\rule[-1.0ex]{0pt}{1.0ex}}     & \shortstack{\rule{0pt}{3ex}31.967 \\[-0.2ex] {\tiny(1.493)}\rule[-1.0ex]{0pt}{1.0ex}}    & \shortstack{\rule{0pt}{3ex}1.956 \\[-0.2ex] {\tiny(0.073)}\rule[-1.0ex]{0pt}{1.0ex}}     & \shortstack{\rule{0pt}{3ex}1.244 \\[-0.2ex] {\tiny(0.046)}\rule[-1.0ex]{0pt}{1.0ex}}     & \shortstack{\rule{0pt}{3ex}0.967 \\[-0.2ex] {\tiny(0.043)}\rule[-1.0ex]{0pt}{1.0ex}}     & \shortstack{\rule{0pt}{3ex}32.055 \\[-0.2ex] {\tiny(1.502)}\rule[-1.0ex]{0pt}{1.0ex}}    & \shortstack{\rule{0pt}{3ex}1.342 \\[-0.2ex] {\tiny(0.041)}\rule[-1.0ex]{0pt}{1.0ex}}     & \shortstack{\rule{0pt}{3ex}0.709 \\[-0.2ex] {\tiny(0.022)}\rule[-1.0ex]{0pt}{1.0ex}}     & \shortstack{\rule{0pt}{3ex}0.572 \\[-0.2ex] {\tiny(0.020)}\rule[-1.0ex]{0pt}{1.0ex}}     \\ 
                                                                                   & \textsf{NP}    & \shortstack{\rule{0pt}{3ex}$>1000$ \\[-0.2ex] {\tiny}\rule[-1.0ex]{0pt}{1.0ex}}     & \shortstack{\rule{0pt}{3ex}\textbf{2.607} \\[-0.2ex] {\tiny(0.119)}\rule[-1.0ex]{0pt}{1.0ex}}     & \shortstack{\rule{0pt}{3ex}\textbf{2.084} \\[-0.2ex] {\tiny(0.101)}\rule[-1.0ex]{0pt}{1.0ex}}     & \shortstack{\rule{0pt}{3ex}\textbf{2.220} \\[-0.2ex] {\tiny(0.103)}\rule[-1.0ex]{0pt}{1.0ex}}     & \shortstack{\rule{0pt}{3ex}$>1000$ \\[-0.2ex] {\tiny}\rule[-1.0ex]{0pt}{1.0ex}} & \shortstack{\rule{0pt}{3ex}\textbf{1.201} \\[-0.2ex] {\tiny(0.055)}\rule[-1.0ex]{0pt}{1.0ex}}     & \shortstack{\rule{0pt}{3ex}\textbf{0.996} \\[-0.2ex] {\tiny(0.045)}\rule[-1.0ex]{0pt}{1.0ex}}     & \shortstack{\rule{0pt}{3ex}\textbf{0.824} \\[-0.2ex] {\tiny(0.043)}\rule[-1.0ex]{0pt}{1.0ex}}     & \shortstack{\rule{0pt}{3ex}$>1000$ \\[-0.2ex] {\tiny}\rule[-1.0ex]{0pt}{1.0ex}} & \shortstack{\rule{0pt}{3ex}\textbf{0.599} \\[-0.2ex] {\tiny(0.024)}\rule[-1.0ex]{0pt}{1.0ex}}     & \shortstack{\rule{0pt}{3ex}\textbf{0.451} \\[-0.2ex] {\tiny(0.037)}\rule[-1.0ex]{0pt}{1.0ex}}     & \shortstack{\rule{0pt}{3ex}\textbf{0.447} \\[-0.2ex] {\tiny(0.046)}\rule[-1.0ex]{0pt}{1.0ex}}     \\ 
                                                                                   & \textsf{N$_1$} & \shortstack{\rule{0pt}{3ex}56.210 \\[-0.2ex] {\tiny(1.334)}\rule[-1.0ex]{0pt}{1.0ex}}    & \shortstack{\rule{0pt}{3ex}8.845 \\[-0.2ex] {\tiny(0.224)}\rule[-1.0ex]{0pt}{1.0ex}}     & \shortstack{\rule{0pt}{3ex}4.161 \\[-0.2ex] {\tiny(0.117)}\rule[-1.0ex]{0pt}{1.0ex}}     & \shortstack{\rule{0pt}{3ex}2.882 \\[-0.2ex] {\tiny(0.110)}\rule[-1.0ex]{0pt}{1.0ex}}     & \shortstack{\rule{0pt}{3ex}54.986 \\[-0.2ex] {\tiny(1.355)}\rule[-1.0ex]{0pt}{1.0ex}}    & \shortstack{\rule{0pt}{3ex}7.468 \\[-0.2ex] {\tiny(0.174)}\rule[-1.0ex]{0pt}{1.0ex}}     & \shortstack{\rule{0pt}{3ex}3.080 \\[-0.2ex] {\tiny(0.077)}\rule[-1.0ex]{0pt}{1.0ex}}     & \shortstack{\rule{0pt}{3ex}1.331 \\[-0.2ex] {\tiny(0.045)}\rule[-1.0ex]{0pt}{1.0ex}}     & \shortstack{\rule{0pt}{3ex}55.376 \\[-0.2ex] {\tiny(1.328)}\rule[-1.0ex]{0pt}{1.0ex}}    & \shortstack{\rule{0pt}{3ex}7.172 \\[-0.2ex] {\tiny(0.158)}\rule[-1.0ex]{0pt}{1.0ex}}     & \shortstack{\rule{0pt}{3ex}2.501 \\[-0.2ex] {\tiny(0.067)}\rule[-1.0ex]{0pt}{1.0ex}}     & \shortstack{\rule{0pt}{3ex}0.939 \\[-0.2ex] {\tiny(0.028)}\rule[-1.0ex]{0pt}{1.0ex}}     \\ 
                                                                                   & \textsf{N$_2$} & \shortstack{\rule{0pt}{3ex}56.932 \\[-0.2ex] {\tiny(1.351)}\rule[-1.0ex]{0pt}{1.0ex}}    & \shortstack{\rule{0pt}{3ex}9.588 \\[-0.2ex] {\tiny(0.230)}\rule[-1.0ex]{0pt}{1.0ex}}     & \shortstack{\rule{0pt}{3ex}5.049 \\[-0.2ex] {\tiny(0.136)}\rule[-1.0ex]{0pt}{1.0ex}}     & \shortstack{\rule{0pt}{3ex}3.768 \\[-0.2ex] {\tiny(0.127)}\rule[-1.0ex]{0pt}{1.0ex}}     & \shortstack{\rule{0pt}{3ex}55.459 \\[-0.2ex] {\tiny(1.350)}\rule[-1.0ex]{0pt}{1.0ex}}    & \shortstack{\rule{0pt}{3ex}7.889 \\[-0.2ex] {\tiny(0.178)}\rule[-1.0ex]{0pt}{1.0ex}}     & \shortstack{\rule{0pt}{3ex}3.558 \\[-0.2ex] {\tiny(0.082)}\rule[-1.0ex]{0pt}{1.0ex}}     & \shortstack{\rule{0pt}{3ex}1.815 \\[-0.2ex] {\tiny(0.053)}\rule[-1.0ex]{0pt}{1.0ex}}     & \shortstack{\rule{0pt}{3ex}55.334 \\[-0.2ex] {\tiny(1.318)}\rule[-1.0ex]{0pt}{1.0ex}}    & \shortstack{\rule{0pt}{3ex}7.396 \\[-0.2ex] {\tiny(0.159)}\rule[-1.0ex]{0pt}{1.0ex}}     & \shortstack{\rule{0pt}{3ex}2.727 \\[-0.2ex] {\tiny(0.069)}\rule[-1.0ex]{0pt}{1.0ex}}     & \shortstack{\rule{0pt}{3ex}1.170 \\[-0.2ex] {\tiny(0.030)}\rule[-1.0ex]{0pt}{1.0ex}}     \\   \cdashline{1-14}
  \parbox[t]{2mm}{\multirow{7}{*}{\rotatebox[origin=c]{90}{\hspace{-2.5cm}\textbf{Model$^2$ 2}}}} & \textsf{S}     & \shortstack{\rule{0pt}{3ex}64.174 \\[-0.2ex] {\tiny(8.057)}\rule[-1.0ex]{0pt}{1.0ex}}    & \shortstack{\rule{0pt}{3ex}19.511 \\[-0.2ex] {\tiny(1.729)}\rule[-1.0ex]{0pt}{1.0ex}}    & \shortstack{\rule{0pt}{3ex}49.658 \\[-0.2ex] {\tiny(31.172)}\rule[-1.0ex]{0pt}{1.0ex}}   & \shortstack{\rule{0pt}{3ex}24.520 \\[-0.2ex] {\tiny(6.089)}\rule[-1.0ex]{0pt}{1.0ex}}    & \shortstack{\rule{0pt}{3ex}44.692 \\[-0.2ex] {\tiny(2.274)}\rule[-1.0ex]{0pt}{1.0ex}}    & \shortstack{\rule{0pt}{3ex}8.877 \\[-0.2ex] {\tiny(0.508)}\rule[-1.0ex]{0pt}{1.0ex}}     & \shortstack{\rule{0pt}{3ex}11.952 \\[-0.2ex] {\tiny(2.017)}\rule[-1.0ex]{0pt}{1.0ex}}    & \shortstack{\rule{0pt}{3ex}8.235 \\[-0.2ex] {\tiny(0.938)}\rule[-1.0ex]{0pt}{1.0ex}}     & \shortstack{\rule{0pt}{3ex}41.376 \\[-0.2ex] {\tiny(1.989)}\rule[-1.0ex]{0pt}{1.0ex}}    & \shortstack{\rule{0pt}{3ex}5.206 \\[-0.2ex] {\tiny(0.272)}\rule[-1.0ex]{0pt}{1.0ex}}     & \shortstack{\rule{0pt}{3ex}6.999 \\[-0.2ex] {\tiny(1.275)}\rule[-1.0ex]{0pt}{1.0ex}}     & \shortstack{\rule{0pt}{3ex}6.415 \\[-0.2ex] {\tiny(1.952)}\rule[-1.0ex]{0pt}{1.0ex}}     \\ 
                                                                                   & \textsf{A}     & \shortstack{\rule{0pt}{3ex}55.000 \\[-0.2ex] {\tiny(2.214)}\rule[-1.0ex]{0pt}{1.0ex}}    & \shortstack{\rule{0pt}{3ex}10.854 \\[-0.2ex] {\tiny(0.463)}\rule[-1.0ex]{0pt}{1.0ex}}    & \shortstack{\rule{0pt}{3ex}8.088 \\[-0.2ex] {\tiny(0.427)}\rule[-1.0ex]{0pt}{1.0ex}}     & \shortstack{\rule{0pt}{3ex}\textbf{6.653} \\[-0.2ex] {\tiny(0.346)}\rule[-1.0ex]{0pt}{1.0ex}}     & \shortstack{\rule{0pt}{3ex}46.215 \\[-0.2ex] {\tiny(1.914)}\rule[-1.0ex]{0pt}{1.0ex}}    & \shortstack{\rule{0pt}{3ex}6.048 \\[-0.2ex] {\tiny(0.211)}\rule[-1.0ex]{0pt}{1.0ex}}     & \shortstack{\rule{0pt}{3ex}3.823 \\[-0.2ex] {\tiny(0.157)}\rule[-1.0ex]{0pt}{1.0ex}}     & \shortstack{\rule{0pt}{3ex}3.134 \\[-0.2ex] {\tiny(0.145)}\rule[-1.0ex]{0pt}{1.0ex}}     & \shortstack{\rule{0pt}{3ex}39.531 \\[-0.2ex] {\tiny(2.437)}\rule[-1.0ex]{0pt}{1.0ex}}    & \shortstack{\rule{0pt}{3ex}4.455 \\[-0.2ex] {\tiny(0.187)}\rule[-1.0ex]{0pt}{1.0ex}}     & \shortstack{\rule{0pt}{3ex}2.671 \\[-0.2ex] {\tiny(0.091)}\rule[-1.0ex]{0pt}{1.0ex}}     & \shortstack{\rule{0pt}{3ex}1.876 \\[-0.2ex] {\tiny(0.050)}\rule[-1.0ex]{0pt}{1.0ex}}     \\ 
                                                                                   & \textsf{H}     & \shortstack{\rule{0pt}{3ex}\textbf{44.195} \\[-0.2ex] {\tiny(2.018)}\rule[-1.0ex]{0pt}{1.0ex}}    & \shortstack{\rule{0pt}{3ex}\textbf{9.293} \\[-0.2ex] {\tiny(0.369)}\rule[-1.0ex]{0pt}{1.0ex}}     & \shortstack{\rule{0pt}{3ex}\textbf{7.743} \\[-0.2ex] {\tiny(0.336)}\rule[-1.0ex]{0pt}{1.0ex}}     & \shortstack{\rule{0pt}{3ex}7.226 \\[-0.2ex] {\tiny(0.359)}\rule[-1.0ex]{0pt}{1.0ex}}     & \shortstack{\rule{0pt}{3ex}\textbf{37.050} \\[-0.2ex] {\tiny(1.678)}\rule[-1.0ex]{0pt}{1.0ex}}    & \shortstack{\rule{0pt}{3ex}\textbf{4.419} \\[-0.2ex] {\tiny(0.130)}\rule[-1.0ex]{0pt}{1.0ex}}     & \shortstack{\rule{0pt}{3ex}\textbf{3.529} \\[-0.2ex] {\tiny(0.132)}\rule[-1.0ex]{0pt}{1.0ex}}     & \shortstack{\rule{0pt}{3ex}\textbf{3.078} \\[-0.2ex] {\tiny(0.122)}\rule[-1.0ex]{0pt}{1.0ex}}     & \shortstack{\rule{0pt}{3ex}\textbf{37.239} \\[-0.2ex] {\tiny(1.781)}\rule[-1.0ex]{0pt}{1.0ex}}    & \shortstack{\rule{0pt}{3ex}\textbf{2.617} \\[-0.2ex] {\tiny(0.071)}\rule[-1.0ex]{0pt}{1.0ex}}     & \shortstack{\rule{0pt}{3ex}\textbf{1.830} \\[-0.2ex] {\tiny(0.059)}\rule[-1.0ex]{0pt}{1.0ex}}     & \shortstack{\rule{0pt}{3ex}\textbf{1.533} \\[-0.2ex] {\tiny(0.047)}\rule[-1.0ex]{0pt}{1.0ex}}     \\ 
                                                                                   & \textsf{L}     & \shortstack{\rule{0pt}{3ex}45.839 \\[-0.2ex] {\tiny(2.036)}\rule[-1.0ex]{0pt}{1.0ex}}    & \shortstack{\rule{0pt}{3ex}10.056 \\[-0.2ex] {\tiny(0.420)}\rule[-1.0ex]{0pt}{1.0ex}}    & \shortstack{\rule{0pt}{3ex}8.425 \\[-0.2ex] {\tiny(0.377)}\rule[-1.0ex]{0pt}{1.0ex}}     & \shortstack{\rule{0pt}{3ex}8.103 \\[-0.2ex] {\tiny(0.434)}\rule[-1.0ex]{0pt}{1.0ex}}     & \shortstack{\rule{0pt}{3ex}37.799 \\[-0.2ex] {\tiny(1.688)}\rule[-1.0ex]{0pt}{1.0ex}}    & \shortstack{\rule{0pt}{3ex}4.723 \\[-0.2ex] {\tiny(0.147)}\rule[-1.0ex]{0pt}{1.0ex}}     & \shortstack{\rule{0pt}{3ex}3.858 \\[-0.2ex] {\tiny(0.155)}\rule[-1.0ex]{0pt}{1.0ex}}     & \shortstack{\rule{0pt}{3ex}3.380 \\[-0.2ex] {\tiny(0.139)}\rule[-1.0ex]{0pt}{1.0ex}}     & \shortstack{\rule{0pt}{3ex}37.666 \\[-0.2ex] {\tiny(1.764)}\rule[-1.0ex]{0pt}{1.0ex}}    & \shortstack{\rule{0pt}{3ex}2.780 \\[-0.2ex] {\tiny(0.076)}\rule[-1.0ex]{0pt}{1.0ex}}     & \shortstack{\rule{0pt}{3ex}1.987 \\[-0.2ex] {\tiny(0.067)}\rule[-1.0ex]{0pt}{1.0ex}}     & \shortstack{\rule{0pt}{3ex}1.684 \\[-0.2ex] {\tiny(0.054)}\rule[-1.0ex]{0pt}{1.0ex}}     \\ 
                                                                                   & \textsf{NP}    & \shortstack{\rule{0pt}{3ex}83.612 \\[-0.2ex] {\tiny(11.526)}\rule[-1.0ex]{0pt}{1.0ex}}   & \shortstack{\rule{0pt}{3ex}19.392 \\[-0.2ex] {\tiny(1.947)}\rule[-1.0ex]{0pt}{1.0ex}}    & \shortstack{\rule{0pt}{3ex}50.023 \\[-0.2ex] {\tiny(31.162)}\rule[-1.0ex]{0pt}{1.0ex}}   & \shortstack{\rule{0pt}{3ex}23.643 \\[-0.2ex] {\tiny(5.612)}\rule[-1.0ex]{0pt}{1.0ex}}    & \shortstack{\rule{0pt}{3ex}43.183 \\[-0.2ex] {\tiny(3.675)}\rule[-1.0ex]{0pt}{1.0ex}}    & \shortstack{\rule{0pt}{3ex}8.184 \\[-0.2ex] {\tiny(0.610)}\rule[-1.0ex]{0pt}{1.0ex}}     & \shortstack{\rule{0pt}{3ex}12.741 \\[-0.2ex] {\tiny(2.541)}\rule[-1.0ex]{0pt}{1.0ex}}    & \shortstack{\rule{0pt}{3ex}8.012 \\[-0.2ex] {\tiny(0.882)}\rule[-1.0ex]{0pt}{1.0ex}}     & \shortstack{\rule{0pt}{3ex}$>1000$ \\[-0.2ex] {\tiny}\rule[-1.0ex]{0pt}{1.0ex}} & \shortstack{\rule{0pt}{3ex}4.920 \\[-0.2ex] {\tiny(0.415)}\rule[-1.0ex]{0pt}{1.0ex}}     & \shortstack{\rule{0pt}{3ex}7.365 \\[-0.2ex] {\tiny(3.529)}\rule[-1.0ex]{0pt}{1.0ex}}     & \shortstack{\rule{0pt}{3ex}2.666 \\[-0.2ex] {\tiny(0.384)}\rule[-1.0ex]{0pt}{1.0ex}}     \\ 
                                                                                   & \textsf{N$_1$} & \shortstack{\rule{0pt}{3ex}74.533 \\[-0.2ex] {\tiny(7.442)}\rule[-1.0ex]{0pt}{1.0ex}}    & \shortstack{\rule{0pt}{3ex}24.714 \\[-0.2ex] {\tiny(1.808)}\rule[-1.0ex]{0pt}{1.0ex}}    & \shortstack{\rule{0pt}{3ex}52.904 \\[-0.2ex] {\tiny(31.198)}\rule[-1.0ex]{0pt}{1.0ex}}   & \shortstack{\rule{0pt}{3ex}26.914 \\[-0.2ex] {\tiny(6.208)}\rule[-1.0ex]{0pt}{1.0ex}}    & \shortstack{\rule{0pt}{3ex}60.935 \\[-0.2ex] {\tiny(1.815)}\rule[-1.0ex]{0pt}{1.0ex}}    & \shortstack{\rule{0pt}{3ex}14.107 \\[-0.2ex] {\tiny(0.644)}\rule[-1.0ex]{0pt}{1.0ex}}    & \shortstack{\rule{0pt}{3ex}15.334 \\[-0.2ex] {\tiny(2.615)}\rule[-1.0ex]{0pt}{1.0ex}}    & \shortstack{\rule{0pt}{3ex}9.512 \\[-0.2ex] {\tiny(0.981)}\rule[-1.0ex]{0pt}{1.0ex}}     & \shortstack{\rule{0pt}{3ex}57.122 \\[-0.2ex] {\tiny(1.265)}\rule[-1.0ex]{0pt}{1.0ex}}    & \shortstack{\rule{0pt}{3ex}10.600 \\[-0.2ex] {\tiny(0.355)}\rule[-1.0ex]{0pt}{1.0ex}}    & \shortstack{\rule{0pt}{3ex}8.841 \\[-0.2ex] {\tiny(1.256)}\rule[-1.0ex]{0pt}{1.0ex}}     & \shortstack{\rule{0pt}{3ex}7.056 \\[-0.2ex] {\tiny(1.991)}\rule[-1.0ex]{0pt}{1.0ex}}     \\ 
                                                                                   & \textsf{N$_2$} & \shortstack{\rule{0pt}{3ex}58.754 \\[-0.2ex] {\tiny(1.251)}\rule[-1.0ex]{0pt}{1.0ex}}    & \shortstack{\rule{0pt}{3ex}14.355 \\[-0.2ex] {\tiny(0.361)}\rule[-1.0ex]{0pt}{1.0ex}}    & \shortstack{\rule{0pt}{3ex}10.417 \\[-0.2ex] {\tiny(0.326)}\rule[-1.0ex]{0pt}{1.0ex}}    & \shortstack{\rule{0pt}{3ex}9.240 \\[-0.2ex] {\tiny(0.310)}\rule[-1.0ex]{0pt}{1.0ex}}     & \shortstack{\rule{0pt}{3ex}55.596 \\[-0.2ex] {\tiny(1.240)}\rule[-1.0ex]{0pt}{1.0ex}}    & \shortstack{\rule{0pt}{3ex}9.865 \\[-0.2ex] {\tiny(0.202)}\rule[-1.0ex]{0pt}{1.0ex}}     & \shortstack{\rule{0pt}{3ex}5.770 \\[-0.2ex] {\tiny(0.129)}\rule[-1.0ex]{0pt}{1.0ex}}     & \shortstack{\rule{0pt}{3ex}4.253 \\[-0.2ex] {\tiny(0.112)}\rule[-1.0ex]{0pt}{1.0ex}}     & \shortstack{\rule{0pt}{3ex}54.586 \\[-0.2ex] {\tiny(1.136)}\rule[-1.0ex]{0pt}{1.0ex}}    & \shortstack{\rule{0pt}{3ex}8.073 \\[-0.2ex] {\tiny(0.168)}\rule[-1.0ex]{0pt}{1.0ex}}     & \shortstack{\rule{0pt}{3ex}3.893 \\[-0.2ex] {\tiny(0.090)}\rule[-1.0ex]{0pt}{1.0ex}}     & \shortstack{\rule{0pt}{3ex}2.251 \\[-0.2ex] {\tiny(0.048)}\rule[-1.0ex]{0pt}{1.0ex}}     \\   \cdashline{1-14}
  \parbox[t]{2mm}{\multirow{7}{*}{\rotatebox[origin=c]{90}{\hspace{-2.5cm}\textbf{Model$^2$ 3}}}} & \textsf{S}     & \shortstack{\rule{0pt}{3ex}$>1000$ \\[-0.2ex] {\tiny}\rule[-1.0ex]{0pt}{1.0ex}} & \shortstack{\rule{0pt}{3ex}$>1000$ \\[-0.2ex] {\tiny}\rule[-1.0ex]{0pt}{1.0ex}} & \shortstack{\rule{0pt}{3ex}$>1000$ \\[-0.2ex] {\tiny}\rule[-1.0ex]{0pt}{1.0ex}} & \shortstack{\rule{0pt}{3ex}$>1000$ \\[-0.2ex] {\tiny}\rule[-1.0ex]{0pt}{1.0ex}} & \shortstack{\rule{0pt}{3ex}$>1000$ \\[-0.2ex] {\tiny}\rule[-1.0ex]{0pt}{1.0ex}} & \shortstack{\rule{0pt}{3ex}$>1000$ \\[-0.2ex] {\tiny}\rule[-1.0ex]{0pt}{1.0ex}} & \shortstack{\rule{0pt}{3ex}$>1000$ \\[-0.2ex] {\tiny}\rule[-1.0ex]{0pt}{1.0ex}} & \shortstack{\rule{0pt}{3ex}$>1000$ \\[-0.2ex] {\tiny}\rule[-1.0ex]{0pt}{1.0ex}} & \shortstack{\rule{0pt}{3ex}$>1000$ \\[-0.2ex] {\tiny}\rule[-1.0ex]{0pt}{1.0ex}} & \shortstack{\rule{0pt}{3ex}$>1000$ \\[-0.2ex] {\tiny}\rule[-1.0ex]{0pt}{1.0ex}} & \shortstack{\rule{0pt}{3ex}$>1000$ \\[-0.2ex] {\tiny}\rule[-1.0ex]{0pt}{1.0ex}} & \shortstack{\rule{0pt}{3ex}$>1000$ \\[-0.2ex] {\tiny}\rule[-1.0ex]{0pt}{1.0ex}} \\ 
                                                                                   & \textsf{A}     & \shortstack{\rule{0pt}{3ex}163.190 \\[-0.2ex] {\tiny(6.890)}\rule[-1.0ex]{0pt}{1.0ex}}   & \shortstack{\rule{0pt}{3ex}64.906 \\[-0.2ex] {\tiny(2.426)}\rule[-1.0ex]{0pt}{1.0ex}}    & \shortstack{\rule{0pt}{3ex}54.408 \\[-0.2ex] {\tiny(3.065)}\rule[-1.0ex]{0pt}{1.0ex}}    & \shortstack{\rule{0pt}{3ex}\textbf{50.649} \\[-0.2ex] {\tiny(3.448)}\rule[-1.0ex]{0pt}{1.0ex}}    & \shortstack{\rule{0pt}{3ex}105.937 \\[-0.2ex] {\tiny(3.146)}\rule[-1.0ex]{0pt}{1.0ex}}   & \shortstack{\rule{0pt}{3ex}30.958 \\[-0.2ex] {\tiny(0.939)}\rule[-1.0ex]{0pt}{1.0ex}}    & \shortstack{\rule{0pt}{3ex}23.008 \\[-0.2ex] {\tiny(0.908)}\rule[-1.0ex]{0pt}{1.0ex}}    & \shortstack{\rule{0pt}{3ex}\textbf{17.259} \\[-0.2ex] {\tiny(0.637)}\rule[-1.0ex]{0pt}{1.0ex}}    & \shortstack{\rule{0pt}{3ex}88.313 \\[-0.2ex] {\tiny(2.991)}\rule[-1.0ex]{0pt}{1.0ex}}    & \shortstack{\rule{0pt}{3ex}20.785 \\[-0.2ex] {\tiny(0.663)}\rule[-1.0ex]{0pt}{1.0ex}}    & \shortstack{\rule{0pt}{3ex}13.002 \\[-0.2ex] {\tiny(0.444)}\rule[-1.0ex]{0pt}{1.0ex}}    & \shortstack{\rule{0pt}{3ex}9.391 \\[-0.2ex] {\tiny(0.258)}\rule[-1.0ex]{0pt}{1.0ex}}     \\ 
                                                                                   & \textsf{H}     & \shortstack{\rule{0pt}{3ex}154.331 \\[-0.2ex] {\tiny(6.337)}\rule[-1.0ex]{0pt}{1.0ex}}   & \shortstack{\rule{0pt}{3ex}\textbf{56.710} \\[-0.2ex] {\tiny(2.168)}\rule[-1.0ex]{0pt}{1.0ex}}    & \shortstack{\rule{0pt}{3ex}\textbf{52.565} \\[-0.2ex] {\tiny(2.661)}\rule[-1.0ex]{0pt}{1.0ex}}    & \shortstack{\rule{0pt}{3ex}51.742 \\[-0.2ex] {\tiny(3.063)}\rule[-1.0ex]{0pt}{1.0ex}}    & \shortstack{\rule{0pt}{3ex}85.538 \\[-0.2ex] {\tiny(2.907)}\rule[-1.0ex]{0pt}{1.0ex}}    & \shortstack{\rule{0pt}{3ex}\textbf{24.621} \\[-0.2ex] {\tiny(0.856)}\rule[-1.0ex]{0pt}{1.0ex}}    & \shortstack{\rule{0pt}{3ex}\textbf{19.685} \\[-0.2ex] {\tiny(0.711)}\rule[-1.0ex]{0pt}{1.0ex}}    & \shortstack{\rule{0pt}{3ex}17.553 \\[-0.2ex] {\tiny(0.504)}\rule[-1.0ex]{0pt}{1.0ex}}    & \shortstack{\rule{0pt}{3ex}\textbf{66.009} \\[-0.2ex] {\tiny(2.657)}\rule[-1.0ex]{0pt}{1.0ex}}    & \shortstack{\rule{0pt}{3ex}\textbf{13.939} \\[-0.2ex] {\tiny(0.391)}\rule[-1.0ex]{0pt}{1.0ex}}    & \shortstack{\rule{0pt}{3ex}\textbf{10.036} \\[-0.2ex] {\tiny(0.290)}\rule[-1.0ex]{0pt}{1.0ex}}    & \shortstack{\rule{0pt}{3ex}\textbf{8.729} \\[-0.2ex] {\tiny(0.257)}\rule[-1.0ex]{0pt}{1.0ex}}     \\ 
                                                                                   & \textsf{L}     & \shortstack{\rule{0pt}{3ex}158.978 \\[-0.2ex] {\tiny(6.335)}\rule[-1.0ex]{0pt}{1.0ex}}   & \shortstack{\rule{0pt}{3ex}62.044 \\[-0.2ex] {\tiny(2.376)}\rule[-1.0ex]{0pt}{1.0ex}}    & \shortstack{\rule{0pt}{3ex}59.057 \\[-0.2ex] {\tiny(2.817)}\rule[-1.0ex]{0pt}{1.0ex}}    & \shortstack{\rule{0pt}{3ex}59.011 \\[-0.2ex] {\tiny(3.302)}\rule[-1.0ex]{0pt}{1.0ex}}    & \shortstack{\rule{0pt}{3ex}90.683 \\[-0.2ex] {\tiny(3.214)}\rule[-1.0ex]{0pt}{1.0ex}}    & \shortstack{\rule{0pt}{3ex}27.366 \\[-0.2ex] {\tiny(0.905)}\rule[-1.0ex]{0pt}{1.0ex}}    & \shortstack{\rule{0pt}{3ex}22.185 \\[-0.2ex] {\tiny(0.770)}\rule[-1.0ex]{0pt}{1.0ex}}    & \shortstack{\rule{0pt}{3ex}20.436 \\[-0.2ex] {\tiny(0.565)}\rule[-1.0ex]{0pt}{1.0ex}}    & \shortstack{\rule{0pt}{3ex}68.057 \\[-0.2ex] {\tiny(2.754)}\rule[-1.0ex]{0pt}{1.0ex}}    & \shortstack{\rule{0pt}{3ex}15.253 \\[-0.2ex] {\tiny(0.409)}\rule[-1.0ex]{0pt}{1.0ex}}    & \shortstack{\rule{0pt}{3ex}11.361 \\[-0.2ex] {\tiny(0.324)}\rule[-1.0ex]{0pt}{1.0ex}}    & \shortstack{\rule{0pt}{3ex}10.144 \\[-0.2ex] {\tiny(0.293)}\rule[-1.0ex]{0pt}{1.0ex}}    \\ 
                                                                                   & \textsf{NP}    & \shortstack{\rule{0pt}{3ex}$>1000$ \\[-0.2ex] {\tiny}\rule[-1.0ex]{0pt}{1.0ex}} & \shortstack{\rule{0pt}{3ex}$>1000$ \\[-0.2ex] {\tiny}\rule[-1.0ex]{0pt}{1.0ex}} & \shortstack{\rule{0pt}{3ex}$>1000$ \\[-0.2ex] {\tiny}\rule[-1.0ex]{0pt}{1.0ex}} & \shortstack{\rule{0pt}{3ex}$>1000$ \\[-0.2ex] {\tiny}\rule[-1.0ex]{0pt}{1.0ex}} & \shortstack{\rule{0pt}{3ex}$>1000$ \\[-0.2ex] {\tiny}\rule[-1.0ex]{0pt}{1.0ex}} & \shortstack{\rule{0pt}{3ex}$>1000$ \\[-0.2ex] {\tiny}\rule[-1.0ex]{0pt}{1.0ex}} & \shortstack{\rule{0pt}{3ex}$>1000$ \\[-0.2ex] {\tiny}\rule[-1.0ex]{0pt}{1.0ex}} & \shortstack{\rule{0pt}{3ex}$>1000$ \\[-0.2ex] {\tiny}\rule[-1.0ex]{0pt}{1.0ex}} & \shortstack{\rule{0pt}{3ex}$>1000$ \\[-0.2ex] {\tiny}\rule[-1.0ex]{0pt}{1.0ex}} & \shortstack{\rule{0pt}{3ex}$>1000$ \\[-0.2ex] {\tiny}\rule[-1.0ex]{0pt}{1.0ex}} & \shortstack{\rule{0pt}{3ex}$>1000$ \\[-0.2ex] {\tiny}\rule[-1.0ex]{0pt}{1.0ex}} & \shortstack{\rule{0pt}{3ex}$>1000$ \\[-0.2ex] {\tiny}\rule[-1.0ex]{0pt}{1.0ex}} \\ 
                                                                                   & \textsf{N$_1$} & \shortstack{\rule{0pt}{3ex}$>1000$ \\[-0.2ex] {\tiny}\rule[-1.0ex]{0pt}{1.0ex}} & \shortstack{\rule{0pt}{3ex}$>1000$ \\[-0.2ex] {\tiny}\rule[-1.0ex]{0pt}{1.0ex}} & \shortstack{\rule{0pt}{3ex}$>1000$ \\[-0.2ex] {\tiny}\rule[-1.0ex]{0pt}{1.0ex}} & \shortstack{\rule{0pt}{3ex}$>1000$ \\[-0.2ex] {\tiny}\rule[-1.0ex]{0pt}{1.0ex}} & \shortstack{\rule{0pt}{3ex}$>1000$ \\[-0.2ex] {\tiny}\rule[-1.0ex]{0pt}{1.0ex}} & \shortstack{\rule{0pt}{3ex}$>1000$ \\[-0.2ex] {\tiny}\rule[-1.0ex]{0pt}{1.0ex}} & \shortstack{\rule{0pt}{3ex}$>1000$ \\[-0.2ex] {\tiny}\rule[-1.0ex]{0pt}{1.0ex}} & \shortstack{\rule{0pt}{3ex}$>1000$ \\[-0.2ex] {\tiny}\rule[-1.0ex]{0pt}{1.0ex}} & \shortstack{\rule{0pt}{3ex}$>1000$ \\[-0.2ex] {\tiny}\rule[-1.0ex]{0pt}{1.0ex}} & \shortstack{\rule{0pt}{3ex}$>1000$ \\[-0.2ex] {\tiny}\rule[-1.0ex]{0pt}{1.0ex}} & \shortstack{\rule{0pt}{3ex}$>1000$ \\[-0.2ex] {\tiny}\rule[-1.0ex]{0pt}{1.0ex}} & \shortstack{\rule{0pt}{3ex}$>1000$ \\[-0.2ex] {\tiny}\rule[-1.0ex]{0pt}{1.0ex}} \\ 
                                                                                   & \textsf{N$_2$} & \shortstack{\rule{0pt}{3ex}\textbf{138.565} \\[-0.2ex] {\tiny(4.212)}\rule[-1.0ex]{0pt}{1.0ex}}   & \shortstack{\rule{0pt}{3ex}96.045 \\[-0.2ex] {\tiny(3.248)}\rule[-1.0ex]{0pt}{1.0ex}}    & \shortstack{\rule{0pt}{3ex}105.830 \\[-0.2ex] {\tiny(3.826)}\rule[-1.0ex]{0pt}{1.0ex}}   & \shortstack{\rule{0pt}{3ex}122.027 \\[-0.2ex] {\tiny(6.356)}\rule[-1.0ex]{0pt}{1.0ex}}   & \shortstack{\rule{0pt}{3ex}\textbf{84.979} \\[-0.2ex] {\tiny(1.846)}\rule[-1.0ex]{0pt}{1.0ex}}    & \shortstack{\rule{0pt}{3ex}37.252 \\[-0.2ex] {\tiny(0.773)}\rule[-1.0ex]{0pt}{1.0ex}}    & \shortstack{\rule{0pt}{3ex}35.195 \\[-0.2ex] {\tiny(0.867)}\rule[-1.0ex]{0pt}{1.0ex}}    & \shortstack{\rule{0pt}{3ex}36.163 \\[-0.2ex] {\tiny(0.679)}\rule[-1.0ex]{0pt}{1.0ex}}    & \shortstack{\rule{0pt}{3ex}67.825 \\[-0.2ex] {\tiny(1.514)}\rule[-1.0ex]{0pt}{1.0ex}}    & \shortstack{\rule{0pt}{3ex}21.406 \\[-0.2ex] {\tiny(0.400)}\rule[-1.0ex]{0pt}{1.0ex}}    & \shortstack{\rule{0pt}{3ex}17.421 \\[-0.2ex] {\tiny(0.277)}\rule[-1.0ex]{0pt}{1.0ex}}    & \shortstack{\rule{0pt}{3ex}17.046 \\[-0.2ex] {\tiny(0.278)}\rule[-1.0ex]{0pt}{1.0ex}}    \\ 
   \hline
\end{tabular}}
\end{table}

Tables~\ref{Tab:results_model_d1} and~\ref{Tab:model_d2_regular} report the means and standard errors of the SEE over $1000$ replications with $d=1$, and $500$ replications for $d=2$. As expected, estimation errors generally decrease with growing $n$ and $p$. The results clearly illustrate the efficiency/robustness trade-off. Under \textbf{Model$^{d}$~1}, the non-robust method \textsf{S} is among the most accurate, but \textsf{H} and \textsf{L} achieve similar accuracy, indicating virtually no efficiency loss in this light-tailed setting. The $L_1$ estimator \textsf{A} and the coordinate-wise median \textsf{N$_2$} are less accurate, while \textsf{S}, \textsf{H}, and \textsf{L} are generally at least as good as the competitors \textsf{P}, \textsf{NP}, \textsf{N$_1$}, and \textsf{N$_2$}, with the largest differences for small $p$. Under \textbf{Model$^{d}$~2}, the non-robust methods \textsf{S}, \textsf{P}, \textsf{NP}, and \textsf{N$_1$} start to struggle, whereas the robust methods remain stable; among them, \textsf{H} and \textsf{L} are the best. Under the extreme \textbf{Model$^{d}$~3}, the non-robust methods fold completely, whereas \textsf{A}, \textsf{H}, \textsf{L}, and \textsf{N$_2$} remain stable and improve with growing $n$ and $p$. The robust spline estimators generally outperform \textsf{N$_2$}, especially for larger $p$.

The unbalanced design yields the same qualitative conclusions, with an even clearer decrease in error as $n$ grows, particularly for the spline estimators. The differences are also more pronounced. The specific results for $d=2$ also exhibit a sharper separation between robust and non-robust methods as the error distribution becomes heavier. Estimation errors generally decrease with $n$ and $p$, although the very sparse setting $p=10$ remains challenging for all methods, with little improvement from increasing $n$. The two-dimensional results also highlight a marked improvement of \textsf{A} with increasing $p$. Under \textbf{Model$^2$~2}, \textsf{A} becomes comparable to \textsf{H} and \textsf{L} at $p=200$, while under \textbf{Model$^2$~3} it is essentially as accurate as them for $p=200$, without requiring a scale estimate (as \textsf{H} and \textsf{L} do). The two-dimensional simulations provide strong evidence of the competitiveness of \textsf{A} at sufficiently dense designs.

% Overall, \textsf{H} and \textsf{L} provide an attractive combination of efficiency under light-tailed errors and robustness under heavy tails, whereas \textsf{A} is less efficient under \textbf{Model~1} but remains highly competitive under \textbf{Models~2} and~\textbf{3}. This comparison is, however, favorable to \textsf{H} and \textsf{L}, since their losses require the residuals to be standardized by a scale parameter and depend on tuning parameters (like $k$ for the Huber estimator \SN{Does L depend on tuning parameters? If not, isn't it better to say ``scale parameter and as such, they''? Is it then correct to say that they ``are fixed in this simulation'' (if the scale estimator is considered a ``tuning parameter'')?} \MC{Yes, I wasn't completely satisfied by that writing either. What do you think about the proposal below?}) which are controlled and fixed in this simulation. 
Overall, \textsf{H} and \textsf{L} provide an attractive combination of efficiency under light-tailed errors and robustness under heavy tails, whereas \textsf{A} is less efficient under \textbf{Model$^{d}$~1} but remains highly competitive under \textbf{Model$^{d}$~2} and~\textbf{3}. This comparison, however, favors \textsf{H} and \textsf{L}: their losses are not scale-equivariant, so their residuals must be standardized by a scale parameter. In our simulation, this scale is set to its true value; in practice, estimating such a scale reliably can be difficult, particularly with noisy or sparse observations, and misspecification may compromise both efficiency and robustness. In contrast, \textsf{A} requires no scale estimator due to its natural scale-equivariance. % neither a scale nor any tuning constant. 
Thus, one may argue that \textsf{A} is the most practical member of the proposed family, % and the one we advocate for applications, 
despite its moderate loss of efficiency in the light-tailed settings.

\section{Analysis of Temperature Anomalies in Southeast Asia}
\label{sec:application}

As a real data example, we consider historical surface temperature anomaly records from the NOAA Merged Land-Ocean Surface Temperature (MLOST) monthly dataset \citep{Lawrimore:2011,Smith:2008}. The MLOST dataset combines land surface air temperatures from the Global Historical Climatology Network Monthly (GHCN-M v3) with sea surface temperatures from the Extended Reconstructed Sea Surface Temperature (ERSST) analysis. The data are provided on a regular $5^\circ \times 5^\circ$ grid and consist of temperature anomalies (in $[^{\circ}C]$) computed with respect to the 30-year climatological average over the period 1961--1990.

In this analysis, we focus on the Southeast Asian region (latitudes $10^\circ\text{S}$--$55^\circ\text{N}$ and longitudes $60^\circ\text{E}$--$150^\circ\text{E}$) during May, covering the period from 1990 to 2016. We selected this region and period because May represents a climatologically important transition period in Southeast Asia, immediately preceding the monsoon season, and is characterized by substantial temperature variability. Moreover, the region is affected by large-scale climate variability associated with phenomena such as El Niño, which may contribute to anomalous temperature observations. We consider a total of $n = 26$ temperature functions measured at $m_{i} = m = 233$ locations, for $i = 1,\dots, 26$.

We use the proposed quantile-loss-based spline estimator to illustrate the spatial patterns of temperature anomalies across different quantile levels. A sensitivity analysis to data contamination is also reported in Appendix~\ref{sec:appendix_C}, showing that the quantile and absolute loss functions are more effective in handling these data than the Huber and least-squares loss functions.

\begin{sidewaysfigure}
    \centering
    \includegraphics[width=0.85\textheight]{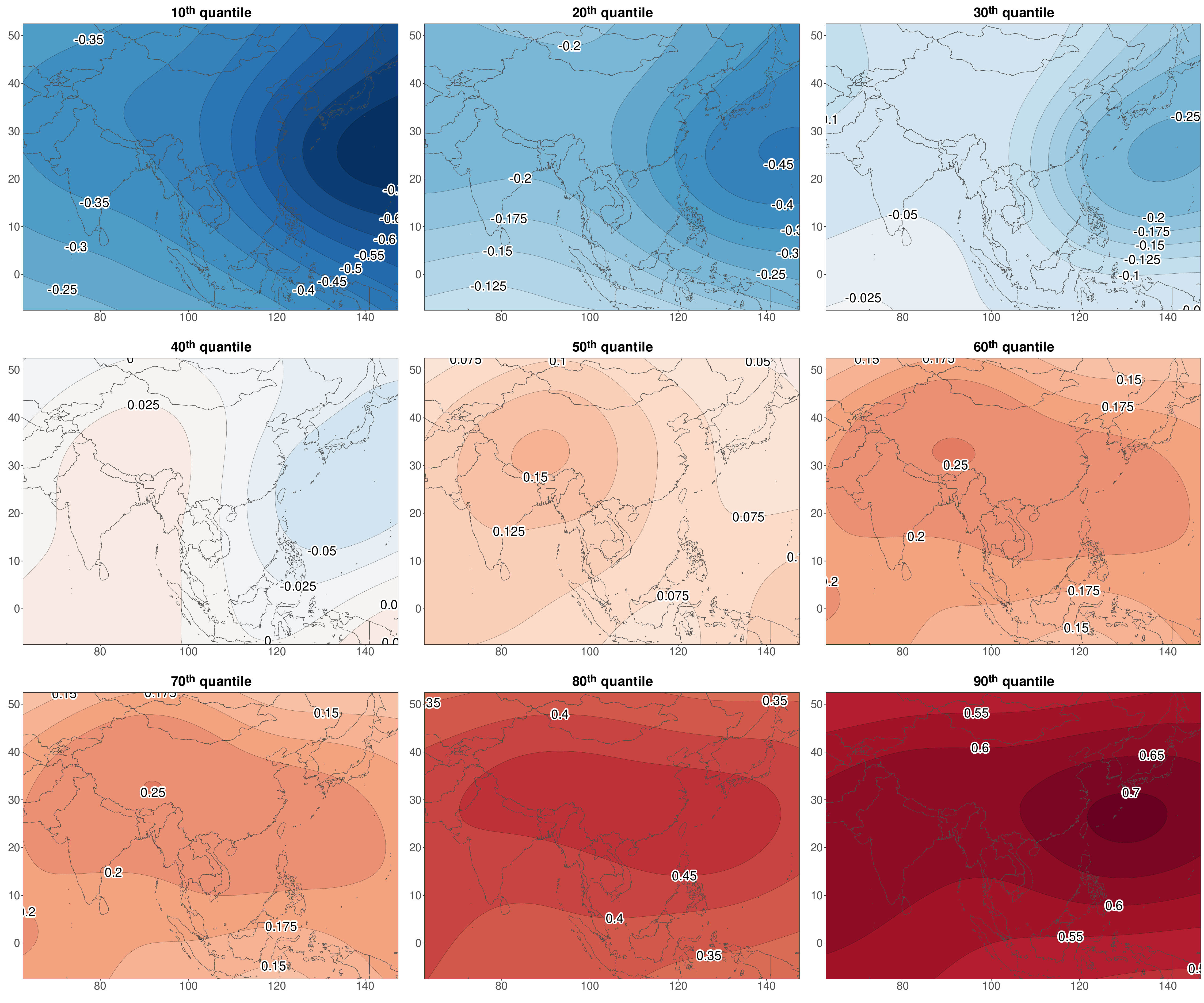}%
    \hspace{0.5cm}
    \raisebox{0.3\height}{%
        \includegraphics[width=0.055\textheight]{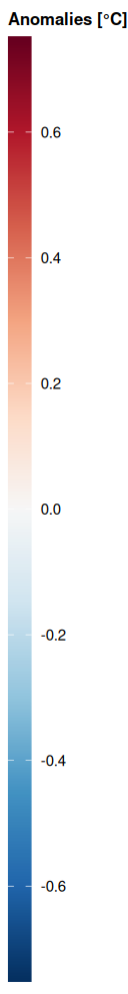}%
    }
    \caption{Contour plots of selected quantiles, ranging from the 10$^{\text{th}}$ to the 90$^{\text{th}}$ percentile, of the temperature anomalies (in [$^{\circ}$C]) in the period 1990--2016. Blue colors indicate negative values while red colors indicate positive values.}
    \label{fig:may_anomalies_1990_2016_contours}
\end{sidewaysfigure}
Figure \ref{fig:may_anomalies_1990_2016_contours} shows the contour plots of nine quantile functions of the temperature anomalies, corresponding to the loss~\eqref{eq::quantile_loss} with $\tau = 0.1, 0.2, \dots, 0.9$. The $50^{\textrm{th}}$ quantile (level 0.5) is the location estimator based on the Absolute loss. These quantile functions summarize the spatial distribution of temperature anomalies using a loss function that is robust to outliers.

We notice that the transition from negative to positive anomalies already occurs at the $40^{\textrm{th}}$ quantile, rather than at the $50^{\textrm{th}}$ quantile. This indicates a slight asymmetry of the distribution towards higher values. Such a pattern is consistent with the expected effects of global warming, although this interpretation should be treated with caution given the short time span considered here ($n = 26$ years). The quantile functions also reveal a clear spatial pattern. Negative anomalies tend to be more pronounced over the oceanic areas in the eastern part of the domain, particularly in the $10^{\textrm{th}}$--$40^{\textrm{th}}$ quantiles. In contrast, positive anomalies are more concentrated over the continental areas covering India and western China in the $40^{\textrm{th}}$--$50^{\textrm{th}}$ quantiles, while they become substantially more homogeneous over the entire domain in the $60^{\textrm{th}}$--$90^{\textrm{th}}$ quantiles.

An explanation for this spatial heterogeneity may lie in the different seasonal behavior of the oceanic and continental areas in the region. Oceanic areas are strongly affected by monsoon-related variability during this period. An earlier start of the monsoon season increases the frequency of storms, winds, and precipitation, potentially leading to negative temperature anomalies relative to the reference climate. On the other hand, a delayed onset may be associated with more positive temperature anomalies. The continental areas considered here are less directly affected by this variability and therefore tend to exhibit less extreme anomalies, as reflected in the $10^{\textrm{th}}$--$30^{\textrm{th}}$ and $90^{\textrm{th}}$ quantiles. The reduced magnitude of the anomalies over these areas also makes the shift towards predominantly positive anomalies more apparent in the central and upper quantiles. We conclude that the proposed quantile framework provides a rich spatial description of the temperature anomaly distribution, while preserving robustness to anomalous observations.

An additional real data analysis exploring the robustness properties of the proposed estimators can be found in Appendix~\ref{sec:appendix_C}.

\section{Discussion}

The results of this paper provide theoretical and practical support for the use of thin-plate spline estimators for multidimensional location estimation from discretely sampled functional data. These estimators emerge as natural generalizations of one-dimensional smoothing splines and are theoretically attractive under mild conditions covering all commonly used loss functions. Our theoretical analysis draws on approximation theory, in particular the theory of conditionally positive definite kernels, and empirical process theory; we show that these two fields may be effectively combined to yield strong theoretical results for a wide class of estimators. Computationally, these estimators may be implemented using either traditional iteratively reweighted least-squares algorithms or modern quadratic programming routines, the latter being substantially faster and more accurate when the loss is not smooth, e.g., the quantile loss. 

Several research directions remain open, most immediately the robust estimation of dispersion from discretely sampled functional data on $\mathbbm{R}^d$. Robust estimation of dispersion is important in its own right, but it is also crucial for downstream procedures, such as robust functional principal component and canonical correlation analysis. We anticipate that combining M-estimation with tensor thin-plate splines will provide a fruitful framework for this problem and we aim to investigate this intriguing combination in future work. Our results could also be extended through a non-asymptotic analysis of thin-plate spline estimators. Such an analysis would be valuable to both theoreticians and practitioners, as it would shed light on critical aspects of the theory, including the phase transition.

\appendix

\section{Appendix: Proofs}

To lighten the notation in the course of our proofs we use $c_0$ to denote generic positive constants whose precise value is unimportant. Thus, the value of $c_0$ may change from appearance to appearance. Given two sequences $\{a_n\}_{n=1}^\infty$ and $\{b_n\}_{n=1}^\infty$ we write $a_n \sim b_n$ if $\lim_{n \to \infty} a_n/b_n = 1$ and $a_n \lesssim b_n$ if $a_n \leq C b_n$ for some universal constant $C$. In our proofs, we occasionally use inequalities of the form $A \leq B$ for $A, B$ random variables defined on the same probability space $(\Omega,\mathcal{A},\mathbb{P})$. These inequalities are to be interpreted in the almost sure sense, i.e., $A(\omega) \leq B(\omega)$ for $\mathbb{P}$-almost all $\omega \in \Omega$.

\begin{proof}[Proof of Proposition~\ref{prop:1}]
Since $2r>d$, Theorem~10.43 in \citet{Wend:2005} identifies $\BL_r(\mathbbm{R}^d)$ with the native space $\mathcal{N}_{\Phi_{dr}}(\mathbbm{R}^d)$ of the kernel $\Phi_{dr} := \eta_{dr}(\|\cdot\|_{\mathbbm{R}^d})$, which is conditionally positive definite with respect to $\mathcal{P}_r$, and the semi-inner product $(\cdot,\cdot)_{\BL_r(\mathbbm{R}^d)}$ agrees with the semi-inner product $(\cdot,\cdot)_{\mathcal{N}_{\Phi_{dr}}}$ of \citet[Definition~10.16]{Wend:2005}. Since $\mathcal{T}$ is $\mathcal{P}_r$-unisolvent, it contains a set $\Xi = \{\boldsymbol{\xi}_1,\ldots,\boldsymbol{\xi}_M\}$ of $M=\dim\mathcal{P}_r$ points that is unisolvent in the minimal sense of \citet{Wend:2005}. Theorem~10.20 in \citet{Wend:2005} then shows that
\begin{align*}
\langle f, g \rangle_{\Xi} := (f,g)_{\BL_r(\mathbbm{R}^d)} + \sum_{k=1}^M f(\boldsymbol{\xi}_k) g(\boldsymbol{\xi}_k)
\end{align*}
is an inner product under which $\BL_r(\mathbbm{R}^d)$ is an RKHS. In particular, with $\|f\|_{\Xi}^2 = \langle f,f\rangle_{\Xi}$, this space is complete and, for each $\mathbf{x}\in\mathbbm{R}^d$, there is a constant $C_{\mathbf{x}}$ with $|f(\mathbf{x})|\le C_{\mathbf{x}}\|f\|_{\Xi}$ for all $f$, as the evaluation functionals are continuous. Let $\|f\|_n^2 = \langle f,f\rangle_n$. Because $\Xi$ consists of points of $\mathcal{T}$, we have $\|f\|_{\Xi}\le\|f\|_n$. Conversely, writing $\mathcal{J}$ for the index pairs $(i,j)$ not used by $\Xi$,
\begin{align*}
\|f\|_n^2 = \|f\|_{\Xi}^2 + \sum_{(i,j)\in\mathcal{J}} f(\mathbf{T}_{ij})^2 \le \Bigl(1 + \sum_{(i,j)\in\mathcal{J}} C_{\mathbf{T}_{ij}}^2\Bigr)\|f\|_{\Xi}^2 .
\end{align*}
It follows $\|\cdot\|_n$ and $\|\cdot\|_{\Xi}$ are equivalent. Hence, $\langle\cdot,\cdot\rangle_n$ is an inner product, $\BL_r(\mathbbm{R}^d)$ is complete under it, and $|f(\mathbf{x})|\leq C_{\mathbf{x}}\|f\|_{\Xi}\leq C_{\mathbf{x}}\|f\|_n$, so every evaluation functional $\delta_{\mathbf{x}}$ is continuous. Thus $\BL_r(\mathbbm{R}^d)$ equipped with $\langle\cdot,\cdot\rangle_n$ is an RKHS.
\end{proof}

\begin{proof}[Proof of Theorem~\ref{thm:1}]
By assumption, $G(\cdot,y)$ is bounded from below for every $y$, and there are finitely many $Y_{ij}$, so $L_n$ is bounded from below by a constant $c>-\infty$. In particular, $L_n$ is bounded from below, so there exists a minimising sequence $\{f_k\}_k \in \BL_{r}(\mathbbm{R}^d)$ with the property
\begin{align*}
\lim_{k \to \infty} L_n(f_k) = \inf_{f \in \BL_{r}(\mathbbm{R}^d)} L_n(f).
\end{align*}
We may assume without loss of generality that $L_n(f_k) \leq L_n(f_1)$ for all $k \in \mathbbm{N}$ switching to a subsequence if necessary. We claim that $\{f_k\}_k$ is bounded in $\|\cdot\|_n$. Indeed, we have $\sup_{k} |f_k|_{\BL_r(\mathbbm{R}^d)}^2 \leq \lambda^{-1} (L_n(f_1) - c) < \infty$, so if a subsequence of $\{f_k\}_k$, $\{f_{n_k}\}_{k}$, were to diverge there would exist a $\mathbf{T}_{ij} \in \mathbbm{R}^d$, for which $|f_{n_k}(\mathbf{T}_{ij})| \to \infty$ as $k \to \infty$ (if no such subsequence existed, we would have $\sup_{k} \|f_k\|_n < \infty$ and $\{f_k \}_k$ would be bounded). But, by assumption, $G(x,y) \to \infty$ as $|x| \to \infty$ for all $y \in \mathbbm{R}$, so it must be that
\begin{align*}
\frac{1}{n} \sum_{i=1}^n \frac{1}{m_i} \sum_{j=1}^{m_i} G(f_k(\mathbf{T}_{ij}),Y_{ij}) \to \infty,
\end{align*}
as $k \to \infty$. But this implies that $L_n(f_{n_k}) > L_n(f_1)$ for all large $k$, a contradiction. Hence, $\{f_k\}_k$ is bounded. Now, since $\{f_k\}_k$ is bounded, by a well-known property of reflexive Banach spaces, there exists a subsequence $\{f_{n_k}\}_k$, with the property that $f_{n_k} \to f_0$ weakly for some $f_0 \in \BL_r(\mathbbm{R}^d)$. By Proposition~\ref{prop:1}, $\BL_r(\mathbbm{R}^d)$ equipped with $(\cdot, \cdot)_n$ is an RKHS, hence weak convergence in  $\BL_r(\mathbbm{R}^d)$ also implies pointwise convergence as, by the reproducing property,
\begin{align*}
\lim_{n \to \infty} f_{n_k}(\mathbf{x}) = \lim_{k \to \infty}\langle f_{n_k} , K(\mathbf{x},\cdot)  \rangle_{n} = \langle f_0 , K(\mathbf{x},\cdot)  \rangle_{n} = f_0(\mathbf{x}),
\end{align*}
for every $\mathbf{x} \in \mathbbm{R}^d$, where $K:\mathbbm{R}^d \times \mathbbm{R}^d \to \mathbbm{R}$ denotes the reproducing kernel. By continuity, $G(f_{n_k}(\mathbf{T}_{ij}), Y_{ij}) \to G(f_{0}(\mathbf{T}_{ij}), Y_{ij})$ for all $(i,j)$, as $k \to \infty$. Since $|\cdot|_{\BL_r(\mathbbm{R}^d)}$ is convex and continuous with respect to $\|\cdot\|_n$ (as $|f|_{\BL_r(\mathbbm{R}^d)} \le \|f\|_n$), its lower level sets are convex and closed, hence weakly closed by Mazur's theorem \citep[Theorem 27.7]{Robinson2020}, so $|\cdot|_{\BL_r(\mathbbm{R}^d)}$ is weakly lower semicontinuous, i.e., $|f_0|_{\BL_r(\mathbbm{R}^d)} \leq \liminf_{k \to \infty} |f_{n_k}|_{\BL_r(\mathbbm{R}^d)}$. Combining these two facts along with the definition of the minimising sequence $\{f_k\}_k$ we obtain
\begin{align*}
\inf_{f \in \BL_{r}(\mathbbm{R}^d)} L_n(f) \leq L_n(f_0) \leq \liminf_{k \to \infty} L_n(f_{n_k}) = \inf_{f \in \BL_{r}(\mathbbm{R}^d)} L_n(f),
\end{align*}
which implies that equality is attained, i.e., $\inf_{f \in \BL_{r}(\mathbbm{R}^d)} L_n(f) = L_n(f_0)$, which proves the existence of a minimiser in $\BL_r(\mathbbm{R}^d)$. 

Let $\mathcal{M}$ denote the set of minimisers of \eqref{eq:objf} in $\BL_r(\mathbbm{R}^d)$. We have just proven that $\mathcal{M}$ is non-empty, so let $g_0 \in \mathcal{M}$ be arbitrary. Since $\mathcal{T}$ is $\mathcal{P}_r$-unisolvent, by \citet[Theorem~4 bis]{Duchon:1977}, there exists a unique natural thin plate spline of order $r$, $f_0 \in \BL_r(\mathbbm{R}^d)$, such that $f_0(\mathbf{T}_{ij}) = g_0(\mathbf{T}_{ij})$ for all $j=1, \ldots, m_i, \ i=1, \ldots, n$ with the property that $|f_0|_{\BL_r(\mathbbm{R}^d)} \leq |g_0|_{\BL_r(\mathbbm{R}^d)}$. Moreover, by the Pythagorean relation in \citet[Corollary 10.25]{Wend:2005},
\begin{align*}
|f_0-g_0|^2_{\BL_r(\mathbbm{R}^d)} + |f_0|^2_{\BL_r(\mathbbm{R}^d)} = |g_0|_{\BL_r(\mathbbm{R}^d)}^2.
\end{align*}
By definition of $g_0$, we also have $L_n(g_0) \leq L_n(f_0)$ and so, by the interpolation and because $\lambda>0$, $|g_0|_{\BL_r(\mathbbm{R}^d)}=|f_0|_{\BL_r(\mathbbm{R}^d)}$. The above display then yields $|f_0-g_0|_{\BL_r(\mathbbm{R}^d)}=0$ and so $f_0 = g_0$ as elements of $\BL_r(\mathbbm{R}^d)$, because $\|f_0-g_0\|_n = 0$. 

To show part (b) note that if $x \mapsto G(x,y)$ is strictly convex for every $y \in \mathbbm{R}$,  then $L_n$ is strictly convex because if $f \neq g \in \BL_r(\mathbb{R}^d)$ then it must be that either $|f-g|_{\BL_r(\mathbbm{R}^d)}>0$ or $\max_{i,j} |f(\mathbf{T}_{ij})-g(\mathbf{T}_{ij})|>0$. If $|f-g|_{\BL_r(\mathbbm{R}^d)}>0$ then $f \mapsto |f|_{\BL_r(\mathbbm{R}^d)}^2$ is strictly convex in the linear segment between $f$ and $g$, because, for every $t \in (0,1)$,
\begin{align*}
|tf+(1-t)g|_{\BL_r(\mathbbm{R}^d)}^2 - t|f|_{\BL_r(\mathbbm{R}^d)}^2 - (1-t)|g|_{\BL_r(\mathbbm{R}^d)}^2 = - t(1-t) |f-g|_{\BL_r(\mathbbm{R}^d)}^2<0.
\end{align*}
Hence, $L_n$ is itself strictly convex. On the other hand, if there exists a $\mathbf{T}_{ij}$ such that $f(\mathbf{T}_{ij}) \neq g(\mathbf{T}_{ij})$ then the strict convexity of $x \mapsto G(x,y)$ yields again the strict convexity of $L_n$. The uniqueness of the minimiser now follows from this strict convexity.

To establish the measurable selection in part (c), we observe that, for a natural thin-plate spline $f$ with knots at the distinct $\mathbf{T}_{ij}$, $|f|_{\BL_r(\mathbbm{R}^d)} < \infty$ if, and only if, its coefficients $(\boldsymbol{\gamma}, \boldsymbol{\delta}) \in \mathbbm{R}^{m_0} \times \mathbbm{R}^M$ are members of $\mathcal{F}(\omega)$, where the linear space $\mathcal{F}(\omega)$ is defined as
\begin{align*}
    \mathcal{F}(\omega) := \{(\boldsymbol{\alpha},\boldsymbol{\beta}) \in \mathbbm{R}^{m_0} \times \mathbbm{R}^M: \boldsymbol{\Phi}^{\top}(\omega) \boldsymbol{\alpha}= \mathbf{0} \},
\end{align*}
see, e.g., \citet{Wend:2005}. Here, $m_0$ is the number of distinct $\mathbf{T}_{ij}$ and $\boldsymbol{\Phi}(\omega) \in \mathbbm{R}^{m_0 \times M}$ is the matrix of discretised polynomials of order $r$ at the $\mathbf{T}_{ij}$. This space depends on $\omega$ through the $\mathbf{T}_{ij}$. Moreover, for every $\boldsymbol{\gamma}$ such that $\boldsymbol{\Phi}^{\top}(\omega) \boldsymbol{\gamma} = \mathbf{0}$, $|f|_{\BL_r(\mathbbm{R}^d)}^2 = \boldsymbol{\gamma}^{\top} \mathbf{K} \boldsymbol{\gamma}$ with $\mathbf{K} \in \mathbbm{R}^{m_0 \times m_0}$ the kernel matrix \citep[see, e.g.][]{Green:1994}. There exists a projection matrix $\boldsymbol{\Pi}(\omega) \in \mathbbm{R}^{m_0 \times m_0}$ onto the null space of $\boldsymbol{\Phi}^{\top}$ that is continuous in the $\mathbf{T}_{ij}$. In particular, $\boldsymbol{\Pi}(\omega) = \mathbf{I} - \boldsymbol{\Phi}(\omega) \left(\boldsymbol{\Phi}^{\top} (\omega)\boldsymbol{\Phi}(\omega) \right)^{-1} \boldsymbol{\Phi}^{\top}(\omega)$, where $\mathbf{I} \in \mathbbm{R}^{m_0 \times m_0}$ is the identity matrix and $\boldsymbol{\Phi}^{\top} (\omega)\boldsymbol{\Phi}(\omega)$ is invertible, because $\boldsymbol{\Phi}$ has full column rank thanks to the unisolvency of the $\mathbf{T}_{ij}$. Thus, the fixed space $\mathbbm{R}^{m_0} \times \mathbbm{R}^M$ is mapped continuously onto $\mathcal{F}(\omega)$ by $(\boldsymbol{\alpha}, \boldsymbol{\beta}) \mapsto (\boldsymbol{\Pi}(\omega)\boldsymbol{\alpha},\boldsymbol{\beta})$, and it suffices to minimise
\begin{align*}
    \mathcal{L}(\boldsymbol{\alpha},\boldsymbol{\beta}, \omega) = \frac{1}{n} \sum_{i=1}^n \frac{1}{m_i} \sum_{j=1}^{m_i} G\Bigl(\bigl(\mathbf{K}(\omega)\boldsymbol{\Pi}(\omega)\boldsymbol{\alpha} + \boldsymbol{\Phi}(\omega)\boldsymbol{\beta}\bigr)_{ij}, Y_{ij}(\omega)\Bigr) + \lambda \boldsymbol{\alpha}^{\top} \boldsymbol{\Pi}^{\top}(\omega) \mathbf{K}(\omega) \boldsymbol{\Pi}(\omega)\boldsymbol{\alpha}.
\end{align*}
This is a Carath\'eodory function: continuous in $(\boldsymbol{\alpha},\boldsymbol{\beta})$ for every $\omega \in \Omega$ and measurable in $\omega$ for every $(\boldsymbol{\alpha},\boldsymbol{\beta})$. By \citet[Theorem~14.37]{Rockafellar_Wets1998} it follows that the minimiser can be chosen measurable. \qedhere

\end{proof}

\begin{proof}[Proof of Proposition~\ref{prop:3}] Since $\mathcal{I}$ is open and bounded with a Lipschitz boundary and $2r>d$, the conditions in \citet[Theorem~3.9]{Agmon:1965} are satisfied. Hence that theorem yields the existence of a positive $\gamma = \gamma(\mathcal{I}, r)$ such that for all $s \geq 1$,
\begin{align*}
\sup_{ \mathbf{t} \in \mathcal{I}} \left|f(\mathbf{t}) \right| & \leq \gamma s^{-r+d/2} \left( s^{r} \left\| f \right\|_{L_2(\mathcal{I})}+ \left\{\sum_{|\alpha|= r} \frac{r!}{\alpha!} \|\partial^{\alpha} f \|_{L_2(\mathcal{I})}^2  \right\}^{1/2}\right) \\ &  = \gamma s^{d/2} \left( \left\|f\right\|_{L_2(\mathcal{I})}+ s^{-r}\left\{\sum_{|\alpha|= r} \frac{r!}{\alpha!} \|\partial^{\alpha} f \|_{L_2(\mathcal{I})}^2  \right\}^{1/2} \right).
\end{align*}
The inequality is valid for all $s \geq 1$, so set $s = \lambda^{-1/(2r)}$ so that $s^{d/2} = \lambda^{-d/(4r)}$ and $s^{-r} = \lambda^{1/2}$ and observe that, for $\lambda \in (0,1]$ we have $s \geq 1$. Therefore, 
\begin{align*}
\sup_{ \mathbf{t} \in \mathcal{I}} \left|f(\mathbf{t}) \right| & \leq \gamma \lambda^{-d/(4r)} \left( \left\| f\right\|_{L_2(\mathcal{I})} + \lambda^{1/2} \left\{\sum_{|\alpha|= r} \frac{r!}{\alpha!} \|\partial^{\alpha} f \|_{L_2(\mathcal{I})}^2 \right\}^{1/2} \right)
\\ & \leq 2^{1/2} \gamma \lambda^{-d/(4r)} \left\{ \left\| f\right\|_{L_2(\mathcal{I})}^2 + \lambda\sum_{|\alpha|= r} \frac{r!}{\alpha!} \|\partial^{\alpha} f \|_{L_2(\mathcal{I})}^2 \right\}^{1/2} =  c_0 \lambda^{-d/(4r)} \|f\|_{\mathcal{H}^r(\mathcal{I}),\lambda},
\end{align*}
as $x+y \leq 2^{1/2}(x^2+y^2)^{1/2}$. The proof is complete.
\end{proof}
\iffalse
Before giving the proof of Theorem~\ref{thm:2}, we note an important implication of Proposition~\ref{prop:3}: this result implies that point evaluations $\mathcal{H}^r(\mathcal{I}) \to \mathbbm{R}: f \mapsto f(\mathbf{t}),\ \mathbf{t} \in \mathcal{I}$, are continuous. These properties now reveal that
\begin{align*}
\left\| \mathcal{R}_{r,\mathcal{I},\lambda}(\mathbf{x}, \cdot)  \right\|_{r, \mathcal{I}, \lambda}^2 = \langle \mathcal{R}_{r,\mathcal{I},\lambda}(\mathbf{x}, \cdot), \mathcal{R}_{r,\mathcal{I},\lambda}(\cdot, \mathbf{x}) \rangle_{r, \mathcal{I}, \lambda} = \mathcal{R}_{r, \mathcal{I}, \lambda}(\mathbf{x}, \mathbf{x}) \leq c_0 \lambda^{-d/(4r)} \left\| \mathcal{R}_{r,\mathcal{I},\lambda}(\mathbf{x}, \cdot)  \right\|_{r, \mathcal{I}, \lambda}.
\end{align*}
Assuming that $\left\| \mathcal{R}_{r,\mathcal{I},\lambda}(\mathbf{x}, \cdot)  \right\|_{r, \mathcal{I}, \lambda}>0$ dividing both sides by $\left\| \mathcal{R}_{r,\mathcal{I},\lambda}(\mathbf{x}, \cdot)  \right\|_{r, \mathcal{I}, \lambda}$ yields
\begin{align}
\label{eq:rkineq}
\left\| \mathcal{R}_{r,\mathcal{I},\lambda}(\mathbf{x}, \cdot)  \right\|_{r, \mathcal{I}, \lambda} \leq c_0 \lambda^{-d/(4r)},
\end{align}
and, since $\lambda>0$, this inequality clearly also holds whenever $\left\| \mathcal{R}_{r,\mathcal{I},\lambda}(\mathbf{x}, \cdot)  \right\|_{r, \mathcal{I}, \lambda} = 0$.
\fi
\begin{proof}[Proof of Theorem~\ref{thm:2}]

Observe that in the case of a fixed design the objective function reduces to
\begin{align*}
L_n(g) = \frac{1}{n m} \sum_{i=1}^n \sum_{j=1}^{m} \rho(Y_{ij}-g(\mathbf{T}_{j})) + \lambda |g|_{\BL_r(\mathbbm{R}^d)}^2.
\end{align*}
Let $C_n = n^{-1} + h_{\max}^{2r}$ and $\mathcal{S}_{D} = \{g \in \BL_r(\mathbbm{R}^d): \|g\|_{\BL_r(\mathbbm{R}^d), \lambda}^2 = D^2\}$. We will show that for every $\varepsilon>0$ there exists a $D = D_{\varepsilon} \geq 1$ such that
\begin{align}
\label{eq:res1}
\liminf_{n \to \infty} \mathbb{P} \left( \inf_{g \in \mathcal{S}_D} L_n(\mu_0+ C_n^{1/2} g) > L_n(\mu_0) \right) \geq 1-\varepsilon.
\end{align}
An application of Theorem~1 of \citet{K:2022a} will then imply that there exists a minimiser $\widehat{\mu}_n \in \BL_r(\mathbbm{R}^d)$ such that
\begin{align*}
\|\widehat{\mu}_n - \mu_0\|_{L_2(\mathcal{I})}^2 + \lambda |\widehat{\mu}_n-\mu_0|_{\BL_r(\mathbbm{R}^d)}^2 = O_{\mathbb{P}}(C_n),
\end{align*}
whence $\|\widehat{\mu}_n - \mu_0\|_{L_2(\mathcal{I})}^2 = O_{\mathbb{P}}(C_n)$, which is the result of the theorem. Therefore, we need to establish \eqref{eq:res1}. To that end, use the absolute continuity of $\rho$ in assumption \ref{A3} and the fundamental theorem of calculus to decompose $ L_n(\mu_0+ C_n^{1/2} g) - L_n(\mu_0)$ as
\begin{align*}
L_n(\mu_0+ C_n^{1/2} g) - L_n(\mu_0) & = \frac{1}{n m} \sum_{i=1}^n \sum_{j=1}^{m} \int_{0}^{C_n^{1/2}g(\mathbf{T}_j)} \{\psi(\epsilon_{ij}+u)-\psi(\epsilon_{ij})\} \dd u
\\ & \quad + \frac{C_n^{1/2}}{nm} \sum_{i=1}^n \sum_{j=1}^m \psi(\epsilon_{ij})g(\mathbf{T}_j) + 2 \lambda C_n^{1/2} (\mu_0, g)_{\BL_r(\mathbbm{R}^d)}  
\\ & \quad + \lambda C_n |g|^2_{\BL_r(\mathbbm{R}^d)},
\end{align*}
Rewrite the above as
\begin{align*}
L_n(\mu_0+ C_n^{1/2} g) - L_n(\mu_0) = \mathbb{E}\{I_1(g)\} + (I_1(g) - \mathbb{E}\{I_1(g)\}) + I_2(g) + I_3(g), 
\end{align*}
with
\begin{align*}
I_1(g) & = \frac{1}{n m} \sum_{i=1}^n \sum_{j=1}^{m} \int_{0}^{C_n^{1/2}g(\mathbf{T}_j)} \{\psi(\epsilon_{ij}+u)-\psi(\epsilon_{ij})\} \dd u + \lambda C_n |g|^2_{\BL_r(\mathbbm{R}^d)},
\\ I_2(g) & = \frac{C_n^{1/2}}{nm} \sum_{i=1}^n \sum_{j=1}^m \psi(\epsilon_{ij})g(\mathbf{T}_j),
\\ I_3(g) & = 2 \lambda C_n^{1/2} (\mu_0,g)_{\BL_r(\mathbbm{R}^d)}.
\end{align*}
Let us define $\mathcal{G}_{D} = \{g \in \BL_r(\mathbbm{R}^d): \|g\|_{\BL_r(\mathbbm{R}^d),\lambda}^2 \leq D^2\}$. Clearly, $\mathcal{S}_D \subset \mathcal{G}_D$ and, by the superadditivity of the infimum, we have
\begin{align*}
\inf_{g \in \mathcal{S}_D} \left[ L_n(\mu_0+ C_n^{1/2} g) - L_n(\mu_0) \right] &\geq \inf_{g \in \mathcal{S}_D} I_1(g)  + \inf_{g \in  \mathcal{S}_D} (I_1(g) - \mathbb{E}\{I_1(g)\})
\\ & \quad + \inf_{g \in \mathcal{S}_D} I_2(g)  +\inf_{g \in  \mathcal{S}_D} I_3(g).
\end{align*}
To prove \eqref{eq:res1}, it suffices to show the following:
\begin{align}
\inf_{g \in  \mathcal{S}_D} I_1(g) & \geq c_0 D^2 C_n \label{eq:boundi1},
\\ \sup_{g \in \mathcal{G}_D} |I_1(g) - \mathbb{E}\{I_1(g)\}| & = o_{\mathbb{P}}(1) C_n\label{eq:boundi1expect},
\\ \sup_{g \in \mathcal{G}_D} |I_2(g)| & = O_{\mathbb{P}}(1)C_n\label{eq:boundi2},
\\ \sup_{g \in \mathcal{G}_D} |I_3(g)| & = O(1) D C_n\label{eq:boundi3},
\end{align}
for some strictly positive constant $c_0$ in \eqref{eq:boundi1}. These bounds jointly yield
\begin{align*}
\inf_{g \in \mathcal{S}_D} L_n(\mu_0+ C_n^{1/2} g) - L_n(\mu_0) \geq c_0 D^2 C_n + O_{\mathbb{P}}(1) DC_n + o_{\mathbb{P}}(1) C_n,
\end{align*}
from which we see that by fixing $D$ large enough, so that the quadratic term dominates, and taking $n \to \infty$ we obtain \eqref{eq:res1}. We begin by establishing \eqref{eq:boundi2}. Changing the order of summation and applying the Schwarz inequality, we obtain
\begin{align}
\left|\frac{1}{n m} \sum_{i=1}^n \sum_{j=1}^m \psi(\epsilon_{ij}) g(\mathbf{T}_j)\right| &= \left|\frac{1}{m} \sum_{j=1}^m \left\{ \frac{1}{n}\sum_{i=1}^n   \psi(\epsilon_{ij})\right\}  g(\mathbf{T}_j) \right|  \nonumber
\\ & \leq \left\{ \frac{1}{m} \sum_{j=1}^m \left| \frac{1}{n}\sum_{i=1}^n   \psi(\epsilon_{ij})\right|^2\right\}^{1/2} \left\{ \frac{1}{m} \sum_{j=1}^m |g(\mathbf{T}_j)|^2\right\}^{1/2} \label{eq:i2b}.
\end{align}
Observe that, by independence (assumption \ref{A5}) and zero means $\mathbb{E}\{\psi(\epsilon_{ij})\} = 0$ and finite second moments (contained in assumption \ref{A7}), we have
\begin{align*}
\mathbb{E}\left\{ \left| \frac{1}{n}\sum_{i=1}^n   \psi(\epsilon_{ij})\right|^2 \right\} = \frac{1}{n^2} \sum_{i=1}^n \sum_{k=1}^n \mathbb{E}\left\{ \psi(\epsilon_{ij}) \psi(\epsilon_{kj})\right\} = \frac{1}{n^2} \sum_{i=1}^n \mathbb{E}\left\{\left| \psi(\epsilon_{1j})\right|^2\right\} = O\left(\frac{1}{n}\right).
\end{align*}
From Markov's inequality, it now follows that the first factor in \eqref{eq:i2b} is $O_{\mathbb{P}}(n^{-1/2})$.

To bound the second factor, observe first that for $g \in \BL_{r}(\mathbbm{R}^d)$ we have $g\vert_{\mathcal{I}}  \in \mathcal{H}^r(\mathcal{I})$. Furthermore, by our limit assumptions, $h_{\max} \to 0$ as $n \to \infty$. The conditions of Theorem~3.3 of \citet{Utr:1988} are satisfied and so there exists a constant $c_0$ depending only on $r,d,\mathcal{I}$ but not on $g \in \BL_r(\mathbbm{R}^d)$, and the constant $B$ of assumption \ref{A2}, such that
\begin{align}
\label{eq:utr1}
\frac{1}{m} \sum_{j=1}^m |g(\mathbf{T}_j)|^2 \leq c_0 \left( \|g\|_{L_2( \mathcal{I})}^2 + h_{\max}^{2r} \sum_{|\alpha|= r} \frac{r!}{\alpha!} \|\partial^{\alpha}g\|_{L_2(\mathcal{I})}^2 \right).
\end{align}
By assumption, $h_{\max}^{2r} \leq c_0 \lambda$ and, as $g \in S_{D} \subset \mathcal{G}_D$, extending the domain of integration of the second term from $\mathcal{I} \subset \mathbbm{R}^d$ to $\mathbbm{R}^d$ proper  leads to
\begin{align*}
\frac{1}{m} \sum_{j=1}^m |g(\mathbf{T}_j)|^2 &\leq c_0 \|g\|_{\BL_r(\mathbbm{R}^d),\lambda}^2 \leq c_0 D^2,
\end{align*}
uniformly for every $g \in \mathcal{G}_D$. Thus, plugging both our bounds into the right-hand side of \eqref{eq:i2b} and recalling the definition of $I_2(g)$, we obtain
\begin{align*}
\sup_{g\in \mathcal{G}_D} |I_2(g)| = O_{\mathbb{P}}(1)  D C_n^{1/2} n^{-1/2}  \leq O_{\mathbb{P}}(1) D C_n,
\end{align*}
where the last inequality follows from $n^{-1/2} \leq C_n^{1/2}$. This is precisely \eqref{eq:boundi2}.

To obtain \eqref{eq:boundi3}, using the integral Schwarz inequality we can see that
\begin{align*}
|I_3(g)| \leq 2 \lambda C_n^{1/2} |\mu_0|_{\BL_r(\mathbbm{R}^d)} |g|_{\BL_r(\mathbbm{R}^d)}  = 2|\mu_0|_{\BL_r(\mathbbm{R}^d)} \lambda^{1/2} C_n^{1/2}\{ \lambda^{1/2} |g|_{\BL_r(\mathbbm{R}^d)} \}.
\end{align*}
But for all $g \in \mathcal{G}_D$ we have $\lambda^{1/2}  |g|_{\BL_r(\mathbbm{R}^d)} \leq D$ and, by assumption, $\mu_0 \in \BL_r(\mathbbm{R}^d)$ so that $ |\mu_0|_{\BL_r(\mathbbm{R}^d)} <\infty$. Hence taking the supremum leads to
\begin{align*}
\sup_{g \in \mathcal{G}_D} |I_3(g)| \leq 2D I_r(\mu_0) \lambda^{1/2} C_n^{1/2} =  O(1) C_n D,
\end{align*}
as $\lambda^{1/2} \sim C^{1/2} h_{\max}^{r} \leq C^{1/2} C_n^{1/2}$. This is precisely \eqref{eq:boundi3}.

We now turn to the lower bound \eqref{eq:boundi1}. By \ref{A1}, $\mathbf{T}_j \in \mathcal{I}$ with $\mathcal{I}$ open and bounded with Lipschitz boundary. Furthermore, by assumption, $2r>d$ and $\lambda \sim C h_{\max}^{2r}$ with $h_{\max} \to 0$, so $\lambda \in (0,1]$ for all large $n$. The conditions of Proposition~\ref{prop:3} are fulfilled and we find
\begin{align*}
C_n^{1/2}\max_{1 \leq j \leq m}|g(\mathbf{T}_j)| \leq  C_n^{1/2} \sup_{\mathbf{t} \in \mathcal{I}} |g(\mathbf{t})| & \leq c_0 C_n^{1/2} \lambda^{-d/(4r)} \left\|g\right\|_{\mathcal{H}^r(\mathcal{I}),\lambda} \\&  \leq c_0 C_n^{1/2} \lambda^{-d/(4r)} \left( \left\|g\right\|_{L_2(\mathcal{I})}^2 + \lambda |g|^2_{\BL_r(\mathbbm{R}^d)} \right)^{1/2},
\end{align*}
where to obtain the last inequality we have extended the domain of integration from $\mathcal{I}$ to $\mathbbm{R}^d$. Our limit conditions entail $(n^{1/2} \lambda^{d/(4r)})^{-1} \to 0$ and, since $2r>d$, we also have $\lambda^{-d/(4r)} h_{\max}^{r} \sim C h_{\max}^{r-d/2} \to 0$, as $h_{\max} \to 0$. Therefore, uniformly in $g \in \mathcal{G}_D$,
\begin{align}
\label{eq:limi1}
C_n^{1/2} \max_{1 \leq j \leq m}|g(\mathbf{T}_j)| \leq c_0 D C_n^{1/2} \lambda^{-d/(4r)} \to 0.
\end{align}
Thus, the domain of integration shrinks to zero so, by Fubini's theorem, \ref{A5} and \ref{A7},
\begin{align*}
\mathbb{E}\{I_1(g)\} & = \frac{1}{n} \sum_{i=1}^n \frac{1}{m} \sum_{j=1}^{m} \int_{0}^{C_n^{1/2}g(\mathbf{T}_j)} \mathbb{E}\{\psi(\epsilon_{1j}+u) \}\dd u + C_n \lambda |g|_{\BL_r(\mathbbm{R}^d)}^2
\\ & = \frac{1}{m} \sum_{j=1}^{m} \int_{0}^{C_n^{1/2}g(\mathbf{T}_j)} \{\delta(\mathbf{T}_j)u+o(u) \}\dd u + C_n \lambda |g|_{\BL_r(\mathbbm{R}^d)}^2
\\ & \geq \inf_{\mathbf{t} \in \mathcal{I}} \delta(\mathbf{t})  \frac{C_n}{2m} \sum_{j=1}^{m} |g(\mathbf{T}_j)|^2(1+o(1))+ C_n \lambda |g|_{\BL_r(\mathbbm{R}^d)}^2,
\end{align*}
where the $o(1)$-term is uniform in $g \in \mathcal{G}_D$ and $\mathbf{T}_j \in \mathcal{I}$.  Approximating the sum from below using Theorem~3.4 of \citet{Utr:1988} and extending the domain of integration we obtain
\begin{align}
\label{eq:utr2}
\frac{1}{m} \sum_{j=1}^{m} |g(\mathbf{T}_j)|^2 &\geq c_0 \|g\|_{L_2(\mathcal{I})}^2 - h_{\max}^{2r}  \sum_{|\alpha|= r} \frac{r!}{\alpha!} \|\partial^{\alpha} g\|_{L_2(\mathcal{I})}^2 \nonumber
\\ & \geq c_0 \|g\|_{L_2(\mathcal{I})}^2 - h_{\max}^{2r}  |g|_{\BL_r(\mathbbm{R}^d)}^2.
\end{align}
Plugging \eqref{eq:utr2} into our lower bound for $\mathbb{E}\{I_1(g)\}$ using $\inf_{\mathbf{t}} \delta(\mathbf{t})>0$ contained in assumption \ref{A7} leads to
\begin{align*}
\mathbb{E}\{I_1(g)\}  & \geq c_0 C_n (\|g\|_{L_2(\mathcal{I})}^2 - h_{\max}^{2r}  |g|_{\BL_r(\mathbbm{R}^d)}^2)(1+o(1)) + C_n \lambda |g|_{\BL_r(\mathbbm{R}^d)}^2
\\ & = c_0 C_n  \|g\|_{L_2(\mathcal{I})}^2(1+o(1))  + C_n \lambda |g|_{\BL_r(\mathbbm{R}^d)}^2 \left(1- c_0\frac{h_{\max}^{2r}}{\lambda}(1+o(1)) \right)
\\ & \geq  c_0 C_n \left( \|g\|_{L_2(\mathcal{I})}^2 + \lambda |g|_{\BL_r(\mathbbm{R}^d)}^2\right) = c_0 C_n \|g\|_{\BL_r(\mathbbm{R}^d),\lambda}^2,
\end{align*}
as $\lambda \sim C h_{\max}^{2r}$ for a large $C$ so that $h_{\max}^{2r}/\lambda$ is small for all large $n$. This allows us to bound $1-c_0(h_{\max}^{2r}/\lambda)(1+o(1)) \geq 1/2$, say, for all large $n$. Taking the infimum over $g \in \mathcal{S}_D$, we now obtain $\inf_{g \in \mathcal{S}_D} \mathbb{E}\{I_1(g)\} \geq c_0 D^2 C_n$, which is precisely \eqref{eq:boundi1}.

To prove \eqref{eq:i2b}, define $f_{g}:\mathbbm{R}^m \to \mathbbm{R}$ by
\begin{align*}
f_{g}(\mathbf{x}) = \frac{1}{m} \sum_{j=1}^m \int_{0}^{C_n^{1/2}g(\mathbf{T}_j)} \{\psi(x_j+u)-\psi(x_j)\} \dd u
\end{align*}
with $\mathbf{x} = (x_1, \dots, x_m)^\top \in \mathbbm{R}^m$, so that
\begin{align*}
|I_1(g) - \mathbb{E}\{I_1(g)\}| = \left| \frac{1}{n} \sum_{i=1}^n \{f_{g}(\boldsymbol{\epsilon}_i) - \mathbb{E}\{f_{g}(\boldsymbol{\epsilon}_i) \}\} \right|
\end{align*}
where $\boldsymbol{\epsilon}_i = (\epsilon_{i1}, \ldots, \epsilon_{im})^\top \in \mathbbm{R}^m$, $i = 1, \dots, n$, are the well-defined projections of the error processes $\epsilon_1, \ldots, \epsilon_n$ onto the points $\mathbf{T}_1, \ldots, \mathbf{T}_m$.  We aim to apply Lemma~19.36 of \citet{VDV:1998} and will therefore need 
\begin{itemize}
\item[(i)] a bound on $\sup_{g \in \mathcal{G}_D}\|f_{g}\|_{\infty}$,
\item[(ii)] a bound on $\sup_{g \in \mathcal{G}_D}\mathbb{E}\{|f_g(\boldsymbol{\epsilon}_1)|^2\}$, and
\item[(iii)] a bound on the covering number in the supremum norm: $\mathcal{N}_{\infty}(\varepsilon,\{f_g,g \in \mathcal{G}_D\})$ for small $\varepsilon>0$.
\end{itemize}
Recall that the covering number $\mathcal{N}_{\infty}(\varepsilon,\{f_g,g \in \mathcal{G}_D\})$ is the smallest number of balls of radius $\varepsilon$ in the sup-norm needed to cover $\{f_{g}, g \in \mathcal{G}_D\}$, i.e., $ N= \mathcal{N}_{\infty}(\varepsilon,\{f_g,g \in \mathcal{G}_D\})$ iff there exist $g_1, \ldots, g_N$ with the property that
\begin{align*}
\sup_{g \in \mathcal{G}_D} \min_{j=1, \ldots, N} \|f_g-f_{g_j}\|_{\infty} < \varepsilon.
\end{align*}
To derive a uniform bound for (i), recall that the domain of integration shrinks to zero as $n \to \infty$ uniformly in $g \in \mathcal{G}_D$ and $\mathbf{T}_j \in \mathcal{I}$ so, by assumption \ref{A4} and the Schwarz inequality,
\begin{align}
\label{eq:unibound}
|f_{g}(\mathbf{x})|  \leq   m^{-1}\sum_{j=1}^m \left|\int_{0}^{C_n^{1/2}g(\mathbf{T}_j)}| \psi(x_j+u)-\psi(x_j)|  \dd u \right|  & \leq  M_1 C_n^{1/2} m^{-1} \sum_{j=1}^m |g(\mathbf{T}_j)| \nonumber
\\ & \leq M_1 C_n^{1/2}\left\{ m^{-1} \sum_{j=1}^m |g(\mathbf{T}_j)|^2\right\}^{1/2} \nonumber
\\ & \leq c_0 C_n^{1/2}
\end{align}
where the last inequality follows from \eqref{eq:utr1} after using  $h_{\max}^{2r} \leq c_0 \lambda$ and absorbing $D$ into $c_0$. This bound holds uniformly in $g \in \mathcal{G}_D$ and $\mathbf{x} \in \mathbbm{R}^m$ and therefore fulfills (i).

We next bound the second moments $\mathbb{E}\{|f_g(\boldsymbol{\epsilon}_1)|^2\}$ uniformly in $g \in \mathcal{G}_D$.  By the Schwarz inequality in both its sum and integral forms along with Proposition~\ref{prop:3}, we obtain
\begin{align*}
\mathbb{E}\{|f_g(\boldsymbol{\epsilon}_1)|^2\} & = m^{-2} \mathbb{E} \left\{ \left| \sum_{j=1}^m \int_{0}^{C_n^{1/2} g(\mathbf{T}_j)} \{\psi(\epsilon_{1j}+u)-\psi(\epsilon_{1j})\} \dd u \right|^2 \right\}
\\ & \leq m^{-1} \sum_{j=1}^m \mathbb{E} \left\{ \left| \int_{0}^{C_n^{1/2} g(\mathbf{T}_j)} \{\psi(\epsilon_{1j}+u)-\psi(\epsilon_{1j})\} \dd u \right|^2\right\}
\\ & \leq m^{-1} C_n^{1/2} \sum_{j=1}^m |g(\mathbf{T}_j)| \left|\int_{0}^{C_n^{1/2} g(\mathbf{T}_j)} \mathbb{E} \left\{ \left|\psi(\epsilon_{1j}+u)-\psi(\epsilon_{1j}) \right|^2 \right\} \dd u \right| 
\\  &\leq m^{-1} C_n^{1/2} \sup_{\mathbf{t} \in \mathcal{I}} |g(\mathbf{t})| \sum_{j=1}^m  \left|\int_{0}^{C_n^{1/2} g(\mathbf{T}_j)} \mathbb{E} \left\{ \left|\psi(\epsilon_{1j}+u)-\psi(\epsilon_{1j}) \right|^2 \right\} \dd u \right| 
\\ &\lesssim   m^{-1} C_n^{1/2} \lambda^{-d/(4r)} \sum_{j=1}^m  \left|\int_{0}^{C_n^{1/2} g(\mathbf{T}_j)} \mathbb{E} \left\{ \left|\psi(\epsilon_{1j}+u)-\psi(\epsilon_{1j}) \right|^2 \right\} \dd u \right|.
\end{align*}
Now, the domain of integration shrinks to zero, so, applying assumption \ref{A6}, we may further bound the above with 
\begin{align}
\label{eq:expectbound}
\mathbb{E}\{|f_g(\boldsymbol{\epsilon}_1)|^2\} & \lesssim C_n^{1/2} \lambda^{-d/(4r)} m^{-1}\sum_{j=1}^m  \left|\int_{0}^{C_n^{1/2} g(\mathbf{T}_j)} |u| \dd u \right| \nonumber \lesssim C_n^{3/2} \lambda^{-d/(4r)} m^{-1} \sum_{j=1}^m |g(\mathbf{T}_j)|^2 \nonumber
\\ & \lesssim C_n^{3/2} \lambda^{-d/(4r)},
\end{align}
where the second inequality follows from the antiderivative $\int |x| \dd x= x|x|/2 + C$ valid for all $x \neq 0$ and the third inequality again from \eqref{eq:utr1} using again $h_{\max}^{2r} \lesssim \lambda$ and extending the domain of integration from $\mathcal{I}$ to the whole of $\mathbbm{R}^d$. This inequality is uniform in $g \in \mathcal{G}_D$, hence \eqref{eq:expectbound} provides the desired bound for (ii).

We now bound the covering number $\mathcal{N}_{\infty}(\varepsilon,\{f_g,g \in \mathcal{G}_D\})$. Observe first that we may restrict attention to $\mathcal{I} \subset \mathbbm{R}^d$, as $f_{g}$ only depends on $g(\mathbf{T}_1), \ldots, g(\mathbf{T}_m)$ and $\mathbf{T}_j \in \mathcal{I}$ so that modifying each $g$ outside $\mathcal{I}$ does not change $f_g$. Moreover, $g \vert_{\mathcal{I}} \in \mathcal{H}^{r}(\mathcal{I})$ and $\lambda |g|_{\BL_r(\mathbbm{R}^d)}^2\leq D^2$ implies $\lambda \sum_{|\alpha|=r} (r!/\alpha!)|\partial^{\alpha} g|_{L_2(\mathcal{I})}^2\leq D^2$. From these observations it follows that we may, instead of $\mathcal{N}_{\infty}(\varepsilon,\{f_g,g \in \mathcal{G}\})$, bound $\mathcal{N}_{\infty}(\varepsilon,\{f_g,g \in \mathcal{G}^{\prime}_D\})$ with $\mathcal{G}^{\prime}_D$ is defined as
\begin{align*}
\mathcal{G}^{\prime}_D = \{g \in \mathcal{H}^{r}(\mathcal{I}): \|g\|_{\mathcal{H}^r(\mathcal{I}),\lambda}^2 \leq D^2\}.
\end{align*}
We first bound $\mathcal{N}_{\infty}(\varepsilon, \mathcal{G}^{\prime}_D)$ and then show how this bound can be used to produce a bound on $\mathcal{N}_{\infty}(\varepsilon,\{f_g,g \in \mathcal{G}^{\prime}_D\})$. By assumption, $\lambda \in (0,1]$ for all large $n$, hence $\|g\|_{\mathcal{H}^r(\mathcal{I}),\lambda}^2 \geq \lambda \|g\|_{\mathcal{H}^r(\mathcal{I}),1}^2$.
This implies $\mathcal{G}^{\prime}_D \subset \{g \in \mathcal{H}^{r}(\mathcal{I}): \|g\|_{\mathcal{H}^r(\mathcal{I}),1}^2 \leq \lambda^{-1} D^2\}$, so $\mathcal{N}_{\infty}(\varepsilon, \mathcal{G}^{\prime}_D) \leq \mathcal{N}_{\infty}(\varepsilon, \{g \in \mathcal{H}^{r}(\mathcal{I}):\|g\|_{\mathcal{H}^r(\mathcal{I}),1}^2 \leq \lambda^{-1} D^2\}\}$. By Proposition~6 in \citet{Cucker:2001}, the latter may be bounded by
\begin{align*}
\log \mathcal{N}_{\infty}(\varepsilon, \{g \in \mathcal{H}^{r}(\mathcal{I}): \|g\|_{\mathcal{H}^r(\mathcal{I}),1}^2 \leq \lambda^{-1} D^2\}\} \lesssim c_0 \lambda^{-d/(2r)} \varepsilon^{-d/r},
\end{align*}
for all small $\varepsilon>0$. This bound leads to
\begin{align*}
\log \mathcal{N}_{\infty}(\varepsilon, \mathcal{G}^{\prime}_D) \leq \log \mathcal{N}_{\infty}(\varepsilon, \{g \in \mathcal{H}^{r}(\mathcal{I}): \|g\|_{\mathcal{H}^r(\mathcal{I}),1}^2 \leq \lambda^{-1} D^2\}\} \lesssim c_0 \lambda^{-d/(2r)} \varepsilon^{-d/r}.
\end{align*}
To pass now to a bound for $\mathcal{N}_{\infty}(\varepsilon,\{f_g,g \in \mathcal{G}^{\prime}_D\})$ observe that for any $g_1, g_2 \in \mathcal{G}^{\prime}_D$ the domain of integration shrinks to zero, so, by \ref{A4},
\begin{align*}
\left|f_{g_1}(\mathbf{x})-f_{g_2}(\mathbf{x})\right| & \leq \frac{1}{m} \sum_{j=1}^m \left| \int_{C_n^{1/2}g_1(\mathbf{T}_j)}^{C_n^{1/2}g_2(\mathbf{T}_j)}| \{\psi(x_j+u)-\psi(x_j)\}| \dd u \right|  \leq M_1 \frac{1}{m} \sum_{j=1}^m \left| \int_{C_n^{1/2}g_1(\mathbf{T}_j)}^{C_n^{1/2}g_2(\mathbf{T}_j)} \dd u \right|
\\ & \leq M_1 C_n^{1/2} \sup_{\mathbf{t} \in \mathcal{I}} \left|g_1(\mathbf{t})-g_2(\mathbf{t})\right|,
\end{align*}
from which it follows that\footnote{A cover for $\mathcal{G}^{\prime}$ of radius $\varepsilon/C_n^{1/2}$ up to constant yields a cover of radius $\varepsilon$ for $\{f_g, g \in \mathcal{G}^{\prime}_D\}$.}
\begin{align}
\label{eq:boundcovn}
\log \mathcal{N}_{\infty}(\varepsilon,\{f_g,g \in \mathcal{G}^{\prime}_D\}) \lesssim c_0 C_n^{d/2r} \lambda^{-d/(2r)} \varepsilon^{-d/r}, \quad \varepsilon>0,
\end{align}
which is the bound on the $(\log)$ covering number needed for (iii).
 
With bounds in (i)--(iii) we can now apply Lemma~19.36 of \citet{VDV:1998}. For each $\delta>0$, the bracketing integral is defined by
\begin{align*}
\mathcal{J}_{\infty}(\delta, \{f_g,g \in \mathcal{G}_D^{\prime}\})  = \int_{0}^{\delta} \sqrt{\log \mathcal{N}_{\infty}(u,\{f_g,g \in \mathcal{G}_D^{\prime}\})} \dd u.
\end{align*}
The integral is finite for all $\delta>0$ as long as $2r>d$. Lemma~19.36 of \citet{VDV:1998} ensures that \eqref{eq:boundi1expect} holds provided that
\begin{align*}
\frac{\mathcal{J}_{\infty}(C_n^{3/4} \lambda^{-d/(8r)},\{f_g,g \in \mathcal{G}^{\prime}_D\})}{n^{1/2}C_n} \to 0 \quad \text{and} \quad \lambda^{d/(4r)}\frac{\mathcal{J}_{\infty}(C_n^{3/4} \lambda^{-d/(8r)},\{f_g,g \in \mathcal{G}^{\prime}_D\})}{n^{1/2}C_n}	\to 0,
\end{align*}
as $n \to \infty$. Since $\lambda \sim C h_{\max}^{2r} \to 0$, which is ensured by our limit conditions, the second limit follows from the first. To verify the first limit, using  \eqref{eq:boundcovn} and evaluating the bracketing integral, we obtain
\begin{align*}
\mathcal{J}_{\infty}(C_n^{3/4} \lambda^{-d/(8r)}, \{f_g,g \in \mathcal{G}^{\prime}_D\}) &  = c_0 C_n^{d/(4r)} \lambda^{-d/(4r)}  C_n^{\frac{3}{4}(1-\frac{d}{2r})} \lambda^{-\frac{d}{8r}(1-\frac{d}{2r})}.
\end{align*}
Therefore, recalling that $C_n = n^{-1} + h_{\max}^{2r}$ so that $n^{-1/4} C_n^{-1/4} \leq 1$ and $\lambda \leq 1$ for all large $n$, we obtain
\begin{align*}
\frac{\mathcal{J}_{\infty}(C_n^{3/4} \lambda^{-d/(8r)},  \{f_g,g \in \mathcal{G}^{\prime}_D\})}{n^{1/2}C_n}  \leq \frac{1}{n^{1/4} \lambda^{3d/(8r)} C_n^{d/(8r)}} \sim \frac{c_0}{n^{1/4} h_{\max}^{d}},
\end{align*}
using again $\lambda \sim C h_{\max}^{2r}$. The latter tends to zero, as $n \to \infty$ on account of our assumed $n h_{\max}^{4d} \to \infty$. Thus, by Lemma~19.36 of \citet{VDV:1998}, \eqref{eq:boundi1expect} holds and the proof is complete.
\end{proof}

\begin{proof}[Proof of Proposition~\ref{prop:4}]

Equip $\mathcal{H}^r(\mathcal{I})$ with its standard inner product $\langle \cdot, \cdot \rangle_{\mathcal{H}^r(\mathcal{I})}$. By definition, $\mathcal{H}^r(\mathcal{I}) \subset L_2(\mathcal{I})$.  Let $\iota: \mathcal{H}^r(\mathcal{I}) \to L_2(\mathcal{I})$ denote the embedding operator from $\mathcal{H}^r(\mathcal{I})$ into $L_2(\mathcal{I})$. By the Rellich-Kondrachov theorem \citep[Theorem 6.3, PART II]{Adams:2003}, $\iota$ is compact. Let $\iota^{\star}: L_2(\mathcal{I}) \to \mathcal{H}^r(\mathcal{I})$ denote its adjoint (which is the well-known RKHS integral operator, see \citet[Proposition 10.28]{Wend:2005}) and set $S = \iota^{\star} \iota: \mathcal{H}^r(\mathcal{I}) \to \mathcal{H}^r(\mathcal{I})$. The linear operator $S$ is nonnegative, self-adjoint and compact, because it is the composition of the bounded operator $\iota^{\star}$ \citep[Proposition 10.28]{Wend:2005} and the compact operator $\iota$. Moreover, its kernel is empty because $Sf = 0$ implies $0= \langle Sf, f\rangle _{\mathcal{H}^r(\mathcal{I})} = \|f\|_{L_2(\mathcal{I})}^2$ and so $f=0$ as an element of $L_2(\mathcal{I})$ and hence also as an element of $\mathcal{H}^r(\mathcal{I}) \subset L_2(\mathcal{I})$. Thus, by the spectral theorem for compact operators \citep[Theorem 16.6]{Robinson2020}, there exist countable orthonormal $\mathcal{H}^r(\mathcal{I})$ basis functions $\{e_k\}_{k=1}^{\infty}$ and positive $\{\mu_k\}_{k=1}^{\infty}$ such that $\mu_1 \geq \mu_2\geq \ldots >0$, $ \mu_k \to 0$ and 
\begin{align*}
    Se_k = \mu_k e_k, \quad (k=1, \ldots ). 
\end{align*}
Put $\sigma_k = \mu_k^{-1} $ and $\phi_k = \sigma_k^{1/2} e_k$. It is clear that $\sigma_k \to \infty$. Moreover, by definition, $\langle\phi_j, \phi_k\rangle_{L_2(\mathcal{I})} = \sigma_k^{1/2} \sigma_j^{1/2} \langle e_k, e_j \rangle_{L_2(\mathcal{I})} =  \sigma_k^{1/2} \sigma_j^{1/2} \langle \iota e_k,  \iota e_j\rangle_{L_2(\mathcal{I})}  = \sigma_k^{1/2} \sigma_j^{1/2} ( Se_k, e_j)_{\mathcal{H}^r(\mathcal{I})} = \sigma_k^{1/2} \sigma_j^{1/2} \mu_k \delta_{kj} = \delta_{kj} $, where we have used the definition of the embedding operator $\iota$ and its adjoint $\iota^{\star}$. Thus, the functions $\{\phi_k\}_{k=1}^{\infty}$ are orthonormal in $L_2(\mathcal{I})$. 

We next show that they are complete. It suffices to show that if $g \in L_2(\mathcal{I})$ satisfies $\langle g, \phi_k\rangle_{L_2(\mathcal{I})} = 0$ for all $k$ then $g = 0$. But the collection $\{\phi_k\}_{k=1}^{\infty}$ is complete in $\mathcal{H}^r(\mathcal{I})$, because  $\{e_k\}_{k=1}^{\infty}$ is complete. Therefore, for any $f \in \mathcal{H}^r(\mathcal{I})$ there exists a sequence $f_n \in \Span(\{\phi_k\}_{k=1}^{\infty})$ such that $f_n \to f$ in $\mathcal{H}^r(\mathcal{I})$. But this also implies that $f_n \to f$ in $L_2(\mathcal{I})$, because the $\mathcal{H}^r(\mathcal{I})$ norm is stronger. Thus, we have $\langle g,f \rangle_{L_2(\mathcal{I})} = 0$ for all $f \in \mathcal{H}^r(\mathcal{I})$. By the density of $\mathcal{H}^r(\mathcal{I})$, it follows that $g = 0$ and hence $\{\phi_k\}_{k=1}^{\infty}$ is a complete orthonormal system of $L_2(\mathcal{I})$.

Now, $\langle \phi_j, \phi_k \rangle_{\mathcal{H}^r(\mathcal{I})} = \delta_{jk} + \sum_{|\alpha|=r} (r!/\alpha!) \langle \partial^{\alpha} e_j,  \partial^{\alpha} e_k \rangle_{L_2(\mathcal{I})} = \delta_{jk} + (e_j, e_k)_{\BL_r(\mathcal{I})} $. At the same time, $\langle \phi_j, \phi_k \rangle_{\mathcal{H}^r(\mathcal{I})}  = \sigma_k^{1/2} \sigma_j^{1/2} \langle e_j, e_k\rangle_{\mathcal{H}^r(\mathcal{I})} = \sigma_k \delta_{jk}$ by orthonomality. Therefore,
\begin{align*}
    (e_j,e_k)_{\BL_r(\mathcal{I})} = \delta_{jk}(\sigma_k-1).
\end{align*}
Set $\gamma_j = \sigma_j -1$ for all $j\geq 1$. Observe that $\gamma_j \geq 0$, because $\gamma_j = |e_j|_{\BL_r(\mathcal{I})}^2 \geq 0$. We need to show that $\gamma_j =0$ for $j=1, \ldots, M$ and that there exist $0<C_1 \leq C_2 < \infty$ such that $C_1 j^{2r/d} \leq \gamma_j \leq C_2 j^{2r/d}$ for $j \geq M+1$. The first assertion is clear, as $|f|_{\BL_r(\mathcal{I})} = 0$ if, and only if, $f \in \mathcal{P}_r(\mathcal{I})$ (recall that $\mathcal{I}$ is connected). This gives $\sigma_1 = \ldots = \sigma_M = 1$ whence  $\gamma_1   = \ldots = \gamma_M = 0$. For $\gamma_j$ with $j \geq M+1$, we recall that in Hilbert spaces the approximation numbers are bounded above and below by $\sigma_j^{1/2}$, see \citet[Definition 2, p. 108]{triebel1995}. These approximation numbers can be read off from \citet[Theorem 4.10.1]{triebel1995}, where in that author's notation, $n = d, p = 2, s=0$ and $\Omega = \mathcal{I}$. The conditions of that theorem are satisfied because $\mathcal{I}$ has a Lipschitz boundary and hence also satisfies the cone condition. Therefore, there exist positive numbers $0<c_1 \leq c_2<\infty$ such that
\begin{align}
    c_1 j^{2r/d} \leq \sigma_j \leq c_2 j^{2r/d}, \quad (j=M+1, \ldots ).
\end{align}
Up to constants, the same bounds hold for the $\gamma_j$ for $j \geq M+1$, because $\sigma_j \to \infty$, which implies that the small perturbation by $-1$ may be incorporated into new constants $0<C_1 \leq C_2 < \infty$. This completes the proof of part \ref{i}.

To prove part \ref{ii}, we observe that $\langle \phi_j, \phi_k \rangle_{\mathcal{H}^r(\mathcal{I}),\lambda} = \delta_{jk} + \lambda \gamma_j \delta_{jk}$. Consequently, the functions $\{\phi_j/(1+\lambda \gamma_j)^{1/2}\}_{j=1}^{\infty}$ are orthonormal. Their completeness follows from the completeness of the $\{e_j\}_{j=1}^{\infty} \subset \mathcal{H}^r(\mathcal{I})$. Lastly, the $\{\phi_j\}_{j=1}^{M}$ form a basis for $\mathcal{P}_r(\mathcal{I})$, because they are orthonormal, and thus linearly independent, and they have dimension $M = \dim(\mathcal{P}_r)$. The conclusion holds because, on finite-dimensional spaces, all bases have the same cardinality.

We lastly prove \ref{iii}. By the fundamental existence theorem for reproducing kernels \citep[Section 2]{Ar:1950}, $\mathbf{x} \mapsto \mathcal{R}_{\mathcal{H}^r(\mathcal{I}),\lambda}(\mathbf{x}, \mathbf{y}) \in \mathcal{H}^r(\mathcal{I})$ for every $\mathbf{y} \in \mathcal{I}$. Thus, we use Parseval's theorem to write
\begin{align*}
\mathcal{R}_{\mathcal{H}^r(\mathcal{I}),\lambda}(\mathbf{x}, \mathbf{y}) = \sum_{j=1}^{\infty} a_j(\mathbf{y}) \frac{\phi_j(\mathbf{x})}{(1+\lambda \gamma_j)^{1/2}},
\end{align*}
where we use the notation $\{a_j(\mathbf{y})\}_{j=1}^{\infty}$ to indicate that the coefficients depend on $\mathbf{y}$. By the same theorem, these coefficients satisfy $\alpha_j(\mathbf{y}) = \langle\mathcal{R}_{\mathcal{H}^r(\mathcal{I}),\lambda}(\cdot, \mathbf{y}), \phi_j \rangle_{r, \mathcal{I}, \lambda}/(1+\lambda \gamma_j)^{1/2}$. Recalling that $\phi_j \in \mathcal{H}^r(\mathcal{I})$, by the reproducing property, $\langle\mathcal{R}_{\mathcal{H}^r(\mathcal{I}),\lambda}(\cdot, \mathbf{y}), \phi_j \rangle_{\mathcal{H}^r(\mathcal{I}), \lambda} = \phi_j(\mathbf{y})$ whence $\alpha_j(\mathbf{y}) = \phi_j(\mathbf{y})/(1+\lambda \gamma_j)^{1/2}$ and 
\begin{align*}
\mathcal{R}_{\mathcal{H}^r(\mathcal{I}),\lambda}(\mathbf{x}, \mathbf{y}) = \sum_{j=1}^{\infty} \frac{\phi_j(\mathbf{y})\phi_j(\mathbf{x})}{1+\lambda \gamma_j}.
\end{align*}
To see the convergence of this series in $L_2(\mathcal{I} \times \mathcal{I})$ note that, by orthogonality and integral comparison, for every $n > m$,
\begin{align*}
\left\| \sum_{j = m+1}^n \frac{\phi_j \phi_j}{1+\lambda \gamma_j} \right\|_{L_2(\mathcal{I} \times \mathcal{I})}^2 = \sum_{j=m+1}^n \frac{1}{(1+\lambda \gamma_j)^2} & \leq  \sum_{j=m+1}^n \frac{1}{(1+\lambda C_1 j^{2r/d})^2} 
\\ & \leq \frac{1}{(1+\lambda C_1(m+1)^{2r/d})^2} + \int_{m+1}^{\infty} \frac{\dd x}{(1+\lambda C_1 x^{2r/d})^2}.
\end{align*}
Since $2r/d >1 >1/2$, the integral is finite for every value of $\lambda>0$ and both terms tend to zero as $m \to \infty$. Thus, for every $\varepsilon>0$ we can find $N \in \mathbbm{N}$ such that, for every $n>m \geq N$,
\begin{align*}
\left\| \sum_{j = m+1}^n \frac{\phi_j \phi_j}{(1+\lambda \gamma_j)} \right\|_{L_2(\mathcal{I} \times \mathcal{I})}^2 < \varepsilon.
\end{align*} 
It follows that $ S_n := \sum_{j = 1}^n \phi_j(\mathbf{x}) \phi_j(\mathbf{y})(1+\lambda \gamma_j)^{-1}$ is a Cauchy sequence in the complete space $L_2(\mathcal{I} \times \mathcal{I})$. Hence $S_n$ converges and $\mathcal{R}_{\mathcal{H}^r(\mathcal{I}),\lambda}$, its limit, is a well-defined element of $L_2(\mathcal{I} \times \mathcal{I})$.
\end{proof}

The proof of Theorem~\ref{thm:3} requires the careful bounding of the norm of the integral operator with $\mathcal{R}_{\mathcal{H}^r(\mathcal{I}),\lambda}$ as its kernel, that is,
\begin{align*}
(\mathcal{T}_{\mathcal{R}_{\mathcal{H}^r(\mathcal{I}),\lambda}}f)(\mathbf{y}) =  \int_{\mathcal{I}} \mathcal{R}_{\mathcal{H}^r(\mathcal{I}),\lambda} \left(\mathbf{x}, \mathbf{y} \right) f(\mathbf{x}) \dd \mathbf{x}, \quad f \in L_2(\mathcal{I}).
\end{align*}
Corollary~\ref{cor:1} shows that its norm is exactly equal to $1$.

\begin{cor}
\label{cor:1}
Under the conditions of Proposition~\ref{prop:4},
\begin{align*}
\left\|\mathcal{T}_{\mathcal{R}_{\mathcal{H}^r(\mathcal{I}),\lambda}}\right\| = \sup_{\|f\|_{L_2(\mathcal{I})}=1} \langle \mathcal{T}_{\mathcal{R}_{\mathcal{H}^r(\mathcal{I}),\lambda}} f, f \rangle_{L_2(\mathcal{I})} = 1,
\end{align*}
and $\langle \mathcal{T}_{\mathcal{R}_{\mathcal{H}^r(\mathcal{I}),\lambda}} f, f \rangle_{L_2(\mathcal{I})} = 1$ for $f \in L_2(\mathcal{I})$ if, and only if, $f \in \mathcal{P}_r(\mathcal{I})$.
\end{cor}

\begin{proof}[Proof of Corollary~\ref{cor:1}]

We begin by observing that, by the symmetry of $\mathcal{R}_{\mathcal{H}^r(\mathcal{I}),\lambda}$, i.e., $\mathcal{R}_{\mathcal{H}^r(\mathcal{I}),\lambda}(\mathbf{x},\mathbf{y}) =\mathcal{R}_{\mathcal{H}^r(\mathcal{I}),\lambda}(\mathbf{y},\mathbf{x})$ for all $\mathbf{x},\mathbf{y} \in \mathcal{I}$, the operator $\mathcal{T}_{\mathcal{R}_{\mathcal{H}^r(\mathcal{I}),\lambda}}$ is self-adjoint:
\begin{align*}
    \langle \mathcal{T}_{\mathcal{R}_{\mathcal{H}^r(\mathcal{I}),\lambda}} f, g \rangle_{L_2(\mathcal{I})} &= \int_{\mathcal{I}} \left[ \int_{\mathcal{I}} \mathcal{R}_{\mathcal{H}^r(\mathcal{I}), \lambda}(\mathbf{x},\mathbf{y}) f(\mathbf{x}) \dd \mathbf{x} \right] g(\mathbf{y}) \dd \mathbf{y} 
     \\ & =  \int_{\mathcal{I}} f(\mathbf{x}) \left[ \int_{\mathcal{I}} \mathcal{R}_{\mathcal{H}^r(\mathcal{I}),\lambda}(\mathbf{y},\mathbf{x})  g(\mathbf{y})    \dd \mathbf{y} \right] \dd \mathbf{x} 
     \\ & = \langle f, \mathcal{T}_{\mathcal{R}_{\mathcal{H}^r(\mathcal{I}),\lambda}} g \rangle_{L_2(\mathcal{I})}.
\end{align*}
for all $f,g \in L_2(\mathcal{I})$. Therefore \citep[see, e.g.][Proposition~2.13]{Conway:1990},  
\begin{align*}
\|\mathcal{T}_{\mathcal{R}_{\mathcal{H}^r(\mathcal{I}),\lambda}}\| = \sup_{\|f\|_{L_2(\mathcal{I})} = 1} |  \langle \mathcal{T}_{\mathcal{R}_{\mathcal{H}^r(\mathcal{I}),\lambda}}f,f \rangle_{L_2(\mathcal{I})}|.
\end{align*}
Now, for $f \in L_2(\mathcal{I})$, use Parseval's theorem to express it as
\begin{align*}
f = \sum_{j=1}^{\infty} f_j \phi_j
\end{align*}
where $\{\phi_j\}_{j=1}^{\infty}$ is the orthonormal basis for $L_2(\mathcal{I})$ of Proposition~\ref{prop:4}. Moreover,  $\|f \|_{L_2(\mathcal{I})}^2 = \sum_{j=1}^{\infty} f_j^2$. By the $L_2(\mathcal{I} \times \mathcal{I})$ expansion of $\mathcal{R}_{r, \mathcal{I}, \lambda}$ established also in Proposition~\ref{prop:4} and the orthonormality of $\{\phi_j\}_{j=1}^\infty$ we find
\begin{align*}
\left|\langle \mathcal{T}_{\mathcal{R}_{\mathcal{H}^r(\mathcal{I}),\lambda}}f, f \rangle_{L_2(\mathcal{I})} \right| = \sum_{j=1}^{\infty} \frac{f_j^2}{1+\lambda \gamma_j} \leq \sum_{j=1}^{\infty} f_j^2 = \left\| f \right\|^2_{L_2(\mathcal{I})},
\end{align*}
as $\lambda \gamma_j \geq 0$ for all $j \geq 1$. Dividing both sides by $\|f\|_{L_2(\mathcal{I})}^2$ and using the linearity of the operator $\mathcal{T}_{\mathcal{R}_{\mathcal{H}^r(\mathcal{I}), \lambda}}$ and that of the inner product, we obtain
\begin{align}
\label{eq:intbound}
\sup_{\left\|f\right\|_{L_2(\mathcal{I})} = 1}|\langle \mathcal{T}_{\mathcal{R}_{\mathcal{H}^r(\mathcal{I}),\lambda}}f, f \rangle_{L_2(\mathcal{I})} |= \sup_{f \in L_2(\mathcal{I})}\frac{|\langle \mathcal{T}_{\mathcal{R}_{\mathcal{H}^r(\mathcal{I}),\lambda}}f, f \rangle_{L_2(\mathcal{I})}|}{\left\|f \right\|_{L_2(\mathcal{I})}^2}  \leq 1,
\end{align}
which proves the first assertion.

We next show that equality prevails in \eqref{eq:intbound} if, and only if, $f \in \mathcal{P}_r(\mathcal{I}) \subset L_2(\mathcal{I})$. Take $f \in \mathcal{P}_r(\mathcal{I})$. Clearly,  $f \in \mathcal{H}^r(\mathcal{I})$ and so, by the reproducing property, we obtain
\begin{align*}
f(\mathbf{x}) = \langle f , \mathcal{R}_{\mathcal{H}^r(\mathcal{I}),\lambda}(\mathbf{x}, \cdot) \rangle_{r, \lambda} & = \int_{\mathcal{I}} f(\mathbf{y}) \mathcal{R}_{\mathcal{H}^r(\mathcal{I}),\lambda}(\mathbf{x}, \mathbf{y}) \dd \mathbf{y} + \lambda \underbrace{(f, \mathcal{R}_{\mathcal{H}^r(\mathcal{I}),\lambda}(\mathbf{x},\cdot))_{\BL_r(\mathcal{I})}}_{=0}
\\ & =  \int_{\mathcal{I}} f(\mathbf{y})  \mathcal{R}_{\mathcal{H}^r(\mathcal{I}),\lambda}(\mathbf{x}, \mathbf{y}) \dd \mathbf{y}
\\ & = (\mathcal{T}_{\mathcal{R}_{\mathcal{H}^r(\mathcal{I}),\lambda}} f)(\mathbf{x}),
\end{align*}
for every $\mathbf{x} \in \mathcal{I}$, as all order $r$ derivatives of $f \in \mathcal{P}_r(\mathcal{I})$ vanish. Hence,
\begin{align}
\label{eq:poly}
\left|\langle \mathcal{T}_{\mathcal{R}_{\mathcal{H}^r(\mathcal{I}),\lambda}}f, f \rangle_{L_2(\mathcal{I})} \right| =  |\langle f, f \rangle_{L_2(\mathcal{I})}|  =  \left\| f\right\|_{L_2(\mathcal{I})}^2,
\end{align}
so that for every $f \in \mathcal{P}_r(\mathcal{I})$, \eqref{eq:intbound} holds with equality.  We now establish the other direction, namely that if \eqref{eq:intbound} holds with equality for some $f \in L_2(\mathcal{I})$, then $f \in \mathcal{P}_r(\mathcal{I})$. To see this, observe that if~\eqref{eq:intbound} were to hold with equality, the $L_2(\mathcal{I} \times \mathcal{I})$ expansion of $\mathcal{R}_{r,\mathcal{I},\lambda}$ would yield
\begin{align*}
\left|\langle \mathcal{T}_{\mathcal{R}_{\mathcal{H}^r(\mathcal{I}),\lambda}}f, f \rangle_{L_2(\mathcal{I})} \right| = \sum_{j=1}^{\infty} \frac{f_j^2}{1+\lambda \gamma_j} = \sum_{j=1}^M f_j^2 + \sum_{ j = M+1}^{\infty} \frac{f_j^2}{1+\lambda \gamma_j}=  \|f\|_{L_2(\mathcal{I})}^2,
\end{align*}
since $\gamma_1 = \ldots = \gamma_M = 0$. But, by Parseval's theorem, $\|f\|_{L_2(\mathcal{I})}^2 =  \sum_{j} f_j^2$, hence the above gives 
\begin{align*}
\sum_{j = M+1}^{\infty} f_j^2 \left(1-\frac{1}{1+\lambda \gamma_j} \right) =  0.
\end{align*}
Now, $\lambda \gamma_j>0$ for all $j \geq M+1$, hence $1-(1+\lambda \gamma_j)^{-1}>0$ so the equality holds if, and only if, $f_j = 0$ for all $j\geq M+1$, which leads to $f = \sum_{j = 1}^M f_j \phi_j$, so that, by Proposition~\ref{prop:4}, $f \in \mathcal{P}_r(\mathcal{I})$ (it is in the linear span of $\phi_j \in \mathcal{P}_{r}(\mathcal{I})$, $j=1, \ldots, M$). The proof is complete. \qedhere

\end{proof}

Before delving into the proof of Theorem~\ref{thm:3} we take note of an important inequality involving $\mathcal{R}_{r,\mathcal{I},\lambda}$ that follows from Proposition~\ref{prop:3}. The function $\mathbf{x} \mapsto \mathcal{R}_{\mathcal{H}^r(\mathcal{I}),\lambda}(\mathbf{x},\mathbf{y}) \in \mathcal{H}^r(\mathcal{I})$ for every $\mathbf{y} \in \mathcal{I}$ and, by the reproducing property, $\mathcal{R}_{\mathcal{H}^r(\mathcal{I}),\lambda}(\mathbf{x},\mathbf{y}) = \langle \mathcal{R}_{\mathcal{H}^r(\mathcal{I}),\lambda}(\mathbf{x},\cdot), \mathcal{R}_{\mathcal{H}^r(\mathcal{I}),\lambda}(\cdot,\mathbf{y}) \rangle_{\mathcal{H}^r(\mathcal{I}),\lambda}$. So, by symmetry and Proposition~\ref{prop:3},
\begin{align*}
\|\mathcal{R}_{\mathcal{H}^r(\mathcal{I}),\lambda}(\mathbf{x},\cdot)\|_{\mathcal{H}^r(\mathcal{I}),\lambda}^2 & = \langle \mathcal{R}_{\mathcal{H}^r(\mathcal{I}),\lambda}(\mathbf{x},\cdot), \mathcal{R}_{\mathcal{H}^r(\mathcal{I}),\lambda}(\cdot,\mathbf{x}) \rangle_{r,\mathcal{I},\lambda} = \mathcal{R}_{\mathcal{H}^r(\mathcal{I}),\lambda}(\mathbf{x},\mathbf{x}) 
\\ & \leq c_0 \lambda^{-d/(4r)} \|\mathcal{R}_{\mathcal{H}^r(\mathcal{I}),\lambda}(\mathbf{x},\cdot)\|_{\mathcal{H}^r(\mathcal{I}),\lambda}.
\end{align*}
If $\|\mathcal{R}_{r,\mathcal{I},\lambda}(\mathbf{x},\cdot)\|_{\mathcal{R}_{\mathcal{H}^r(\mathcal{I}),\lambda}}>0$, divide by $\|\mathcal{R}_{\mathcal{H}^r(\mathcal{I}),\lambda}(\mathbf{x},\cdot)\|_{\mathcal{H}^r(\mathcal{I}),\lambda}$ to find
\begin{align}
\label{eq:rkineq}
\sup_{\mathbf{x} \in \mathcal{I}: \ \|\mathcal{R}_{\mathcal{H}^r(\mathcal{I}),\lambda}(\mathbf{x},\cdot)\|_{\mathcal{H}^r(\mathcal{I}),\lambda}>0}\|\mathcal{R}_{\mathcal{H}^r(\mathcal{I}),\lambda}(\mathbf{x},\cdot)\|_{\mathcal{H}^r(\mathcal{I}),\lambda} \leq c_0 \lambda^{-d/(4r)}.
\end{align}
Since $\lambda>0$, \eqref{eq:rkineq} is trivially true also when $\|\mathcal{R}_{\mathcal{H}^r(\mathcal{I}),\lambda}(\mathbf{x},\cdot)\|_{\mathcal{H}^r(\mathcal{I}),\lambda} = 0$. In the proof of Theorem~\ref{thm:3} below we assume for convenience that the $\mathbf{T}_{ij}$ are uniformly distributed on $\mathcal{I}$ with Lebesgue density
\begin{align*}
f(\mathbf{x}) = \begin{cases} \dfrac{1}{\mu(\mathcal{I})}, & \mathbf{x} \in \mathcal{I},\\ 0, & \text{otherwise}, \end{cases}
\end{align*}
where $\mu$ denotes the Lebesgue measure. The general case of a density bounded away from zero and infinity may be treated with similar arguments. We will also assume without loss of generality $\mu(\mathcal{I}) = 1$. Under these assumptions, for every $g \in L_2(\mathcal{I})$:
\begin{align*}
\mathbb{E}\left\{\frac{1}{n} \sum_{i=1}^n \frac{1}{m_i} \sum_{j=1}^{m_i} |g(\mathbf{T}_{ij})|^2\right\} = \frac{1}{n} \sum_{i=1}^n \frac{1}{m_i} \sum_{j=1}^{m_i} \mathbb{E}\{|g(\mathbf{T}_{ij})|^2\} = \int_{\mathcal{I}} |g(\mathbf{t})|^2 \dd \mathbf{t} = \|g\|_{L_2(\mathcal{I})}^2.
\end{align*}
We make use of this observation in the proof of Theorem~\ref{thm:3} below.

\begin{proof}[Proof of Theorem~\ref{thm:3}]

Define 
\begin{align*}
L_n(g) = \frac{1}{n} \sum_{i=1}^n \frac{1}{m_i} \sum_{j=1}^{m_i} \rho(Y_{ij}-g(\mathbf{T}_{ij})) + \lambda  |g|_{\BL_r(\mathbbm{R}^d)}^2,
\end{align*}
and let $C_n = n^{-1}+ (nm\lambda^{d/(2r)})^{-1} + \lambda$. As in the proof of Theorem~\ref{thm:2}, the result will be established if we show that there exists a $D \geq 1$ such that
\begin{align}
\label{eq:toprove}
\inf_{g \in \mathcal{S}_D} L_n(\mu_0+C_n^{1/2}g)>L_n(\mu_0),
\end{align}
with high probability, for all large $n$. To that end, use the fundamental theorem of calculus to decompose 
\begin{align*}
L_n(\mu_0+C_n^{1/2}g)-L_n(\mu_0) = I_1(g) + I_2(g) + I_3(g),
\end{align*}
where
\begin{align*}
I_1(g) & = \frac{1}{n} \sum_{i=1}^n \frac{1}{m_i} \sum_{j=1}^{m_i} \int_{0}^{C_n^{1/2}g(\mathbf{T}_{ij})} \{\psi(\epsilon_{ij}+u)-\psi(\epsilon_{ij})\} \dd u + \lambda C_n |g|_{\BL_r(\mathbbm{R}^d)}^2,
\\ I_2(g) & = \frac{C_n^{1/2}}{n} \sum_{i=1}^n \frac{1}{m_i} \sum_{j=1}^m \psi(\epsilon_{ij})g(\mathbf{T}_{ij}),
\\ I_3(g) & = 2 \lambda C_n^{1/2} (\mu_0,g)_{\BL_r(\mathbbm{R}^d)}.
\end{align*}
It will be shown that
\begin{align}
\inf_{g \in  \mathcal{S}_D } \mathbb{E}\{I_1(g) \} \geq c_0 D^2C_n, \label{bound:I1r}
\\ \sup_{g \in \mathcal{G}_D} |I_1(g) - \mathbb{E}\{I_1(g) \}| = o_{\mathbb{P}}(1) C_n, \label{bound:I1mdr}
\\ \sup_{g \in \mathcal{G}_D} |I_2(g)| = O_{\mathbb{P}}(1) D C_n, \label{bound:I2r}
\\ \sup_{g \in \mathcal{G}_D} |I_3(g)| = O(1) D C_n,  \label{bound:I3r}
\end{align}
which in combination lead to
\begin{align*}
\inf_{g \in \mathcal{S}_D} L_n(\mu_0+C_n^{1/2}g)-L_n(\mu_0) \geq c_0 D^2 C_n + O_{\mathbb{P}}(1) DC_n + o_{\mathbb{P}}(1) C_n.
\end{align*}
Here, $\mathcal{S}_D$ and $\mathcal{G}_D$ are defined as in the proof of Theorem~\ref{thm:2}. For sufficiently large $D$ the right-hand side of this inequality will be positive with high probability as $n \to \infty$ and \eqref{eq:toprove} will hold.

The bound \eqref{bound:I3r} can be proven in a similar manner as in the proof of Theorem~\ref{thm:2}, hence we begin by proving \eqref{bound:I2r}. Recall that $\mathcal{H}^r(\mathcal{I})$ is a reproducing kernel Hilbert space under $\langle \cdot, \cdot \rangle_{\mathcal{H}^{r}(\mathcal{I}),\lambda}$ with reproducing kernel $\mathcal{R}_{\mathcal{H}^r(\mathcal{I}),\lambda}$. Therefore, by the reproducing property and linearity of the inner product,
\begin{align*}
I_2(g) & = \frac{C_n^{1/2}}{n} \sum_{i=1}^n \frac{1}{m_i} \sum_{j=1}^{m_i} \psi(\epsilon_{ij}) \langle \mathcal{R}_{\mathcal{H}^r(\mathcal{I}),\lambda}(\mathbf{T}_{ij},\cdot), g \rangle_{\mathcal{H}^r(\mathcal{I}),\lambda} \\ &  = C_n^{1/2} \bigg \langle \frac{1}{n} \sum_{i=1}^n \frac{1}{m_i} \sum_{j=1}^{m_i} \psi(\epsilon_{ij}) \mathcal{R}_{\mathcal{H}^r(\mathcal{I}),\lambda} (\mathbf{T}_{ij},\cdot), g  \bigg \rangle_{\mathcal{H}^r(\mathcal{I}),\lambda}.
\end{align*}
Consequently, by the Schwarz inequality,
\begin{align*}
\sup_{g \in \mathcal{G}_D} |I_2(g)| \leq C_n^{1/2} \left\| \frac{1}{n} \sum_{i=1}^n \frac{1}{m_i} \sum_{j=1}^{m_i} \psi(\epsilon_{ij}) \mathcal{R}_{\mathcal{H}^r(\mathcal{I}),\lambda}(\mathbf{T}_{ij},\cdot) \right\|_{\mathcal{H}^r(\mathcal{I}),\lambda} D,
\end{align*}
where we have used $\|g\|_{\mathcal{H}^r(\mathcal{I}),\lambda}^2 \leq \|g\|_{\BL_r(\mathbbm{R}^d), \lambda}^2$, by extension of the domain integration on the derivative term from $\mathcal{I}$ to $\mathbbm{R}^d$. Now, by assumptions \ref{A5} (i.i.d. errors), \ref{B1} (independence) and \ref{A7} ($\mathbb{E}\{\psi(\epsilon_1(\mathbf{t}))\} = 0$ for all $\mathbf{t} \in \mathcal{I})$ we have
\begin{align*}
\mathbb{E}\{\psi(\epsilon_{ij}) |\mathcal{T}\} = \int_\mathbbm{R}\psi(y) \dd F_{\mathbf{T}_{ij}}(y) = 0,
\end{align*}
as well as
\begin{align*}
\mathbb{E}\{\psi(\epsilon_{ij})\psi(\epsilon_{i^{\prime} k}) |\mathcal{T}\} &= \int_{\mathbbm{R}} \int_\mathbbm{R}\psi(y) \psi(x) \dd F_{\mathbf{T}_{ij}}(y) \dd F_{\mathbf{T}_{i^{\prime} k}}(x) 
\\ & = \left[ \int_{\mathbbm{R}} \psi(y) \dd F_{\mathbf{T}_{ij}}(y) \right] \left[ \int_{\mathbbm{R}} \psi(x) \dd F_{\mathbf{T}_{i^{\prime} k}}(x) \right] = 0
\end{align*}
for all $j,k$ whenever $i \neq i^{\prime}$. Therefore,
\begin{align*}
& \mathbb{E}\left\{ \left\| \frac{1}{n} \sum_{i=1}^n \frac{1}{m_i} \sum_{j=1}^{m_i} \psi(\epsilon_{ij}) \mathcal{R}_{\mathcal{H}^r(\mathcal{I}),\lambda}(\mathbf{T}_{ij},\cdot) \right\|_{\mathcal{H}^r(\mathcal{I}),\lambda}^2 \bigg \vert \mathcal{T}\right\} 
\\ & \quad =  \sum_{i=1}^n \frac{1}{n^2m_i^2} \sum_{j=1}^{m_i} \sum_{\ell=1}^{m_i} \mathcal{R}_{\mathcal{H}^r(\mathcal{I}),\lambda}(\mathbf{T}_{ij}, \mathbf{T}_{i \ell}) \mathbb{E}\{ \psi(\epsilon_{ij}) \psi(\epsilon_{i \ell}) | \mathcal{T} \} = Z_1 + Z_2,
\end{align*}
say, with 
\begin{align*}
Z_1 & = \frac{1}{n^2} \sum_{i=1}^n \frac{1}{m_i^2} \sum_{j=1}^{m_i}  \mathcal{R}_{\mathcal{H}^r(\mathcal{I}),\lambda}(\mathbf{T}_{ij}, \mathbf{T}_{ij}) \mathbb{E}\{| \psi(\epsilon_{ij})|^2 | \mathcal{T} \},\\ 
Z_2 & = \frac{1}{n^2} \sum_{i=1}^n \frac{1}{m_i^2} \sum_{j \neq \ell } \mathcal{R}_{\mathcal{H}^r(\mathcal{I}),\lambda}(\mathbf{T}_{ij}, \mathbf{T}_{i \ell}) \mathbb{E}\{ \psi(\epsilon_{ij}) \psi(\epsilon_{i \ell}) | \mathcal{T} \}.
\end{align*}
We bound these two terms as follows. For $Z_1$, using the finiteness of $\mathbb{E}\{| \psi(\epsilon_{ij})|^2 | \mathcal{T} \}$ entailed in assumption \ref{A7},  $\mathcal{R}_{\mathcal{H}^r(\mathcal{I}),\lambda}(\mathbf{T}_{ij}, \mathbf{T}_{ij}) = \|\mathcal{R}_{\mathcal{H}^r(\mathcal{I}),\lambda}(\mathbf{T}_{ij}, \cdot)\|_{\mathcal{H}^r(\mathcal{I}),\lambda}^2$ and \eqref{eq:rkineq}, we find
\begin{align*}
Z_1 \leq \frac{c_0}{n^2 \lambda^{\frac{d}{2r}}} \sum_{i=1}^n \frac{1}{m_i} = O\left( \frac{1}{n m \lambda^{d/(2r)}}\right),
\end{align*}
since, by definition of the harmonic mean, $m^{-1} = n^{-1} \sum_{i=1}^n m_i^{-1}$. To derive the order of $Z_2$, recall that the $\mathbf{T}_{ij}$ are i.i.d. and independent of $\epsilon_i$, which are identically distributed. Hence, using our simplifying assumption that density of the $\mathbf{T}_{ij}$ is uniform on $\mathcal{I}$ with $\mu(\mathcal{I}) = 1$,
\begin{align*}
E\{Z_2\} = \frac{1}{n^2} \sum_{i=1}^n \frac{m_i(m_i-1)}{m_i^2} \int_{\mathcal{I}} \int_{\mathcal{I}} \mathcal{R}_{\mathcal{H}^r(\mathcal{I}),\lambda}(\mathbf{x},\mathbf{y}) \mathbb{E}\{\psi(\epsilon_1(\mathbf{x})) \psi(\epsilon_1(\mathbf{y})) \} \dd \mathbf{x} \dd \mathbf{y},
\end{align*}
or, after using Fubini's theorem and the notation of Corollary~\ref{cor:1},
\begin{align*}
E\{Z_2\} & = \frac{1}{n^2} \sum_{i=1}^n \frac{m_i(m_i-1)}{m_i^2}\mathbb{E}\{ \langle \mathcal{T}_{\mathcal{R}_{\mathcal{H}^r(\mathcal{I}),\lambda}} \psi(\epsilon_1), \psi(\epsilon_1)\rangle_{L_2(\mathcal{I})}\}
\\ & =  \frac{1}{n^2} \sum_{i=1}^n \frac{m_i(m_i-1)}{m_i^2} \mathbb{E}\left\{ \frac{ \langle \mathcal{T}_{\mathcal{R}_{\mathcal{H}^r(\mathcal{I}),\lambda}} \psi(\epsilon_1), \psi(\epsilon_1)\rangle_{L_2(\mathcal{I})}}{\|\psi(\epsilon_1)\|^2_{L_2(\mathcal{I})}} \|\psi(\epsilon_1)\|^2_{L_2(\mathcal{I})}\right\}
\\ & \leq \frac{1}{n^2} \sum_{i=1}^n \frac{m_i(m_i-1)}{m_i^2} \mathbb{E}\left\{ \sup_{\|f\|_{L_2(\mathcal{I})}=1} |\langle \mathcal{T}_{\mathcal{R}_{\mathcal{H}^r(\mathcal{I}),\lambda}} f, f\rangle_{L_2(\mathcal{I})}|\|\psi(\epsilon_1)\|^2_{L_2(\mathcal{I})}\right\}
\\ & \leq \frac{1}{n^2} \sum_{i=1}^n \frac{m_i(m_i-1)}{m_i^2} \mathbb{E}\left\{ \|\psi(\epsilon_1)\|^2_{L_2(\mathcal{I})}\right\}
\\ & = O\left(\frac{1}{n}\right),
\end{align*}
where we have used Corollary~\ref{cor:1} and assumption \ref{A7} in the last step. Observing now that the bound on $Z_1$ is deterministic, combining the bounds on $Z_1$ and $Z_2$ and iterating expectations leads to
\begin{align*}
\mathbb{E}\left\{ \left\| \frac{1}{n} \sum_{i=1}^n \frac{1}{m_i} \sum_{j=1}^{m_i} \psi(\epsilon_{ij}) \mathcal{R}_{\mathcal{H}^r(\mathcal{I}),\lambda}(\mathbf{T}_{ij},\cdot) \right\|_{\mathcal{H}^r(\mathcal{I}),\lambda}^2\right\} = O\left( \frac{1}{n} +\frac{1}{n m \lambda^{d/(2r)}}\right).
\end{align*}
The definition of $C_n^{1/2}$ and Markov's inequality now directly lead to 
\eqref{bound:I2r}.

We next prove \eqref{bound:I1r}. Recall from the proof of Theorem~\ref{thm:2} that under our limit assumptions we have
\begin{align*}
\lim_{n \to \infty} C_n^{1/2} \sup_{\mathbf{x} \in \mathcal{I}} |g(\mathbf{x})| = 0,
\end{align*}
uniformly for $g \in \mathcal{G}_D$. Conditioning on $\mathcal{T}$ and employing  assumption \ref{A7}, we obtain
\begin{align*}
\mathbbm{E}\{I_1(g) | \mathcal{T}\} & = \frac{1}{n} \sum_{i=1}^n \frac{1}{m_i} \sum_{j=1}^{m_i} \int_{0}^{C_n^{1/2}g(\mathbf{T}_{ij})} \mathbb{E}\{\psi(\epsilon_{ij}+u)| \mathcal{T}\} \dd u + \lambda C_n  |g|_{\BL_r(\mathbbm{R}^d)}^2
\\ & = \frac{1}{n} \sum_{i=1}^n \frac{1}{m_i} \sum_{j=1}^{m_i} \int_{0}^{C_n^{1/2}g(\mathbf{T}_{ij})} \{\delta(\mathbf{T}_{ij})u + o(u)\} \dd u + \lambda C_n |g|_{\BL_r(\mathbbm{R}^d)}^2
\\ & \geq \inf_{\mathbf{t} \in \mathcal{I}} \delta(\mathbf{t}) \frac{C_n}{2n} \sum_{i=1}^n \frac{1}{m_i} \sum_{j=1}^{m_i} |g(\mathbf{T}_{ij})|^2 (1+o(1))+ \lambda C_n  |g|_{\BL_r(\mathbbm{R}^d)}^2
\\ & \geq c_0 \frac{C_n}{n} \sum_{i=1}^n \frac{1}{m_i} \sum_{j=1}^{m_i} |g(\mathbf{T}_{ij})|^2+ \lambda C_n  |g|_{\BL_r(\mathbbm{R}^d)}^2
\end{align*}
for some strictly positive $c_0$. Taking now expectations with respect to $\mathcal{T}$ we find
\begin{align*}
\mathbbm{E}\{I_1(g)\} & \geq \min(c_0,1) C_n \|g\|_{\BL_r(\mathbbm{R}^d),\lambda}^2,
\end{align*}
which leads to $\inf_{g \in \mathcal{S}_D} \mathbb{E}\{I_1(g) \} \geq c_0 D^2C_n$, for some strictly positive $c_0>0$, and this is precisely \eqref{bound:I1r}.

To complete the proof, we establish \eqref{bound:I1mdr}. To that end, write
\begin{align*}
I_1(g)  - \mathbb{E}\{I_1(g)\} & = \frac{1}{n} \sum_{i=1}^n \frac{1}{m_i} \sum_{j=1}^{m_i} \int_{0}^{C_n^{1/2}g(\mathbf{T}_{ij})} \{\psi(\epsilon_{ij}+u)-\psi(\epsilon_{ij})\} \dd u
\\ & \quad -\mathbb{E}\left\{ \frac{1}{n} \sum_{i=1}^n \frac{1}{m_i} \sum_{j=1}^{m_i} \int_{0}^{C_n^{1/2}g(\mathbf{T}_{ij})} \{\psi(\epsilon_{ij}+u)-\psi(\epsilon_{ij})\}  \dd u \right\}
\\ & =: \frac{1}{n} \sum_{i=1}^n Z_{i}(g).
\end{align*}
Clearly, the $\{Z_i\}_{i=1}^n$ are independent mean-zero processes. They are bounded, as, by our previous argument, the domain of integration shrinks to zero so that, by assumption \ref{A4} and Proposition~\ref{prop:3},
\begin{align}
\max_{1 \leq i \leq n}|Z_{i}(g)| \leq 2 M_1 C_n^{1/2} \sup_{g \in \mathcal{G}_D} \|g\|_{\infty} \leq 2 M_1  C_n^{1/2} D\lambda^{-d/(4r)} \label{eq:boundU}.
\end{align}
Furthermore, by the Schwarz inequality in both its sum and integral forms,
\begin{align*}
&\frac{1}{n} \sum_{i=1}^n \frac{1}{m_i^2}\mathbb{E} \left\{ \left|\sum_{j=1}^{m_i} \int_{0}^{C_n^{1/2}g(\mathbf{T}_{ij})} \{\psi(\epsilon_{ij}+u)-\psi(\epsilon_{ij})\} \dd u \right|^2 \bigg \vert \mathcal{T}\right\}
\\ & \leq \frac{1}{n} \sum_{i=1}^n \frac{1}{m_i}\mathbb{E} \left\{ \sum_{j=1}^{m_i} \left|\int_{0}^{C_n^{1/2}g(\mathbf{T}_{ij})} \{\psi(\epsilon_{ij}+u)-
\psi(\epsilon_{ij})\} \dd u \right|^2 \bigg \vert \mathcal{T}\right\}
\\ & \leq \frac{1}{n} \sum_{i=1}^n \frac{1}{m_i}\sum_{j=1}^{m_i} C_n^{1/2}|g(\mathbf{T}_{ij})|  \left|\int_{0}^{C_n^{1/2}g(\mathbf{T}_{ij})} \mathbb{E} \left\{  \left|\psi(\epsilon_{ij}+u)-
\psi(\epsilon_{ij})\right|^2 \bigg \vert \mathcal{T} \right\}\dd u  \right|
\\ & \leq c_0 \frac{C_n^{1/2} \lambda^{-d/(4r)}}{n} \sum_{i=1}^n \frac{1}{m_i} \sum_{j=1}^{m_i} \left|\int_{0}^{C_n^{1/2}g(\mathbf{T}_{ij})} \mathbb{E} \left\{  \left|\psi(\epsilon_{ij}+u)-
\psi(\epsilon_{ij})\right|^2  \vert \mathcal{T} \right\}\dd u  \right|.
\end{align*}
Now, observe that, by independence,
\begin{align*}
\mathbb{E} \left\{  \left|\psi(\epsilon_{ij}+u)-
\psi(\epsilon_{ij})\right|^2  \vert \mathcal{T} \right\} & = \int_{\mathbbm{R}} |\psi(y+u) - \psi(y)|^2 \dd F_{\mathbf{T}_{ij}}(y)
\\ & \leq \sup_{\mathbf{t} \in \mathcal{I}} \int_{\mathbbm{R}} |\psi(y+u) - \psi(y)|^2 \dd F_{\mathbf{t}}(y)
\\ & = \sup_{\mathbf{t} \in \mathcal{I}}\mathbb{E} \left\{\left|\psi(\epsilon_1(\mathbf{t})+u)-\psi(\epsilon_1(\mathbf{t}))\right|^2 \right\},
\end{align*}
where $F_{\mathbf{t}}$ denotes the CDF of $\epsilon_1(\mathbf{t})$. Thus, by assumption \ref{A6}, we find
\begin{align*}
& \frac{1}{n} \sum_{i=1}^n \frac{1}{m_i^2}\mathbb{E} \left\{ \left|\sum_{j=1}^{m_i} \int_{0}^{C_n^{1/2}g(\mathbf{T}_{ij})} \{\psi(\epsilon_{ij}+u)-\psi(\epsilon_{ij})\} \dd u \right|^2 \bigg \vert \mathcal{T}\right\}
\\ & \leq c_0 \frac{C_n^{1/2} \lambda^{-d/(4r)}}{n} \sum_{i=1}^n \frac{1}{m_i} \sum_{j=1}^{m_i} \left|\int_{0}^{C_n^{1/2}g(\mathbf{T}_{ij})} |u| \dd u  \right| \leq c_0 \frac{C_n^{3/2} \lambda^{-d/(4r)}}{n} \sum_{i=1}^n \frac{1}{m_i} \sum_{j=1}^{m_i} |g(\mathbf{T}_{ij})|^2.
\end{align*}
Recalling that the $\mathbf{T}_{ij}$ are uniformly distributed in $\mathcal{I}$ with $\mu(\mathcal{I})=1$, taking expectations with respect to $\mathcal{T}$, we can bound the (unconditional) expectation with
\begin{align}
\label{eq:varU}
\frac{1}{n} \sum_{i=1}^n \frac{1}{m_i^2}\mathbb{E} \left\{ \left|\sum_{j=1}^{m_i} \int_{0}^{C_n^{1/2}g(\mathbf{T}_{ij})} \{\psi(\epsilon_{ij}+u)-\psi(\epsilon_{ij})\} \dd u \right|^2 \right\} \lesssim  C_n^{3/2} \lambda^{-d/(4r)}.
\end{align}
The bound is uniform in $g \in \mathcal{G}_D$. Finally, using again that the domain of integration shrinks to zero, assumption \ref{A4} yields
\begin{align*}
|Z_{i}(g) - Z_{i}(g^{\prime})| \lesssim C_n^{1/2} \sup_{\mathbf{t} \in \mathcal{I}}|g(\mathbf{t})-g^{\prime}(\mathbf{t})|.
\end{align*}
It follows as in \eqref{eq:boundcovn} in the proof of Theorem~\ref{thm:2} that the log of the $\varepsilon$-covering number in the supremum norm is of order $C_n^{d/(2r)} \lambda^{-d/(2r)} \varepsilon^{-d/r}$. This implies that the generalized entropy with bracketing integral in Theorem~8.13 of \citet{vandeGeer:2000} may be bounded to the order of
\begin{align*}
C_n^{d/(4r)} \lambda^{-d/(4r)} \int_{0}^{C_n^{3/4} \lambda^{-d/(8r)}} u^{-d/(2r)} \dd u \lesssim  C_n^{d/(4r)} \lambda^{-d/(4r)}  C_n^{3/4(1-d/(2r))} \lambda^{-d/8r+d^2/16r^2}.
\end{align*}
The definition of $C_n$, in particular, $C_n \geq n^{-1}$ and $C_n \geq \lambda$, and our limit assumptions now jointly give
\begin{align*}
\frac{C_n^{d/(4r)} \lambda^{-d/(4r)}  C_n^{3/4(1-d/(2r))} \lambda^{-d/8r+d^2/16r^2}}{n^{1/2} C_n} \leq \frac{1}{n^{1/4} \lambda^{d/(4r)}} \to 0,
\end{align*}
as $n \to \infty$. Hence Theorem~8.13 of \citet{vandeGeer:2000} applies and gives the existence of a global $c_0$ with the property that for every $\varepsilon>0$ there exists an $N = N(\varepsilon)$ such that for all $n \geq N$,
\begin{align*}
\mathbb{P} \left( \sup_{g \in \mathcal{G}_D} \left| \frac{1}{n} \sum_{i=1}^n Z_{i}(g) \right| \geq \varepsilon C_n \right) \leq  c_0 \exp \left[ - \varepsilon^2 c_0 n C_n^{1/2} \lambda^{d/(4r)} \right].
\end{align*}
For $n \to \infty$, the right side tends to zero as, $nC_n^{1/2} \lambda^{d/(4r)} \geq n^{1/2} \lambda^{d/(4r)} \to \infty$ implied by our limit assumptions. This shows that $|I_1(g) - \mathbb{E}\{I_1(g)\}| = o_{\mathbb{P}}(1) C_n$ uniformly, thereby establishing \eqref{bound:I1mdr}. The result of the theorem follows.

\end{proof}

\section{Supplement to Section~\ref{sec:5}: Numerical Results} \label{sec: appendix B}

\begin{figure}[htpb]
    \centering
    \includegraphics[width=0.4\textwidth]{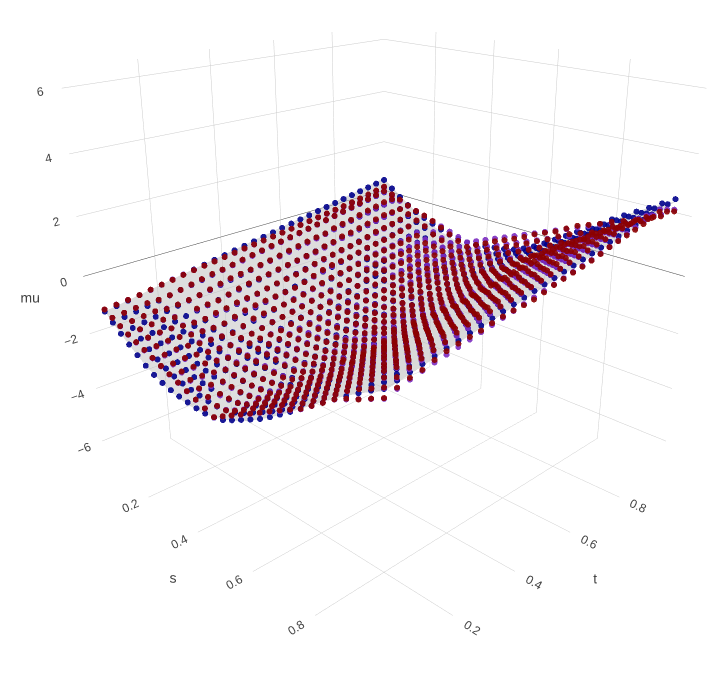} 
    \hspace{2cm}
    \includegraphics[width=.4\textwidth]{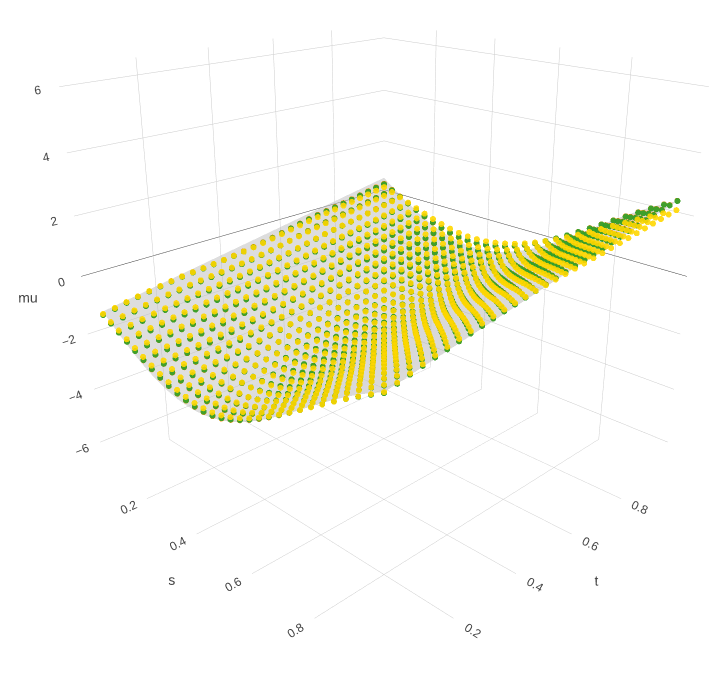}
    \includegraphics[width=0.4\textwidth]{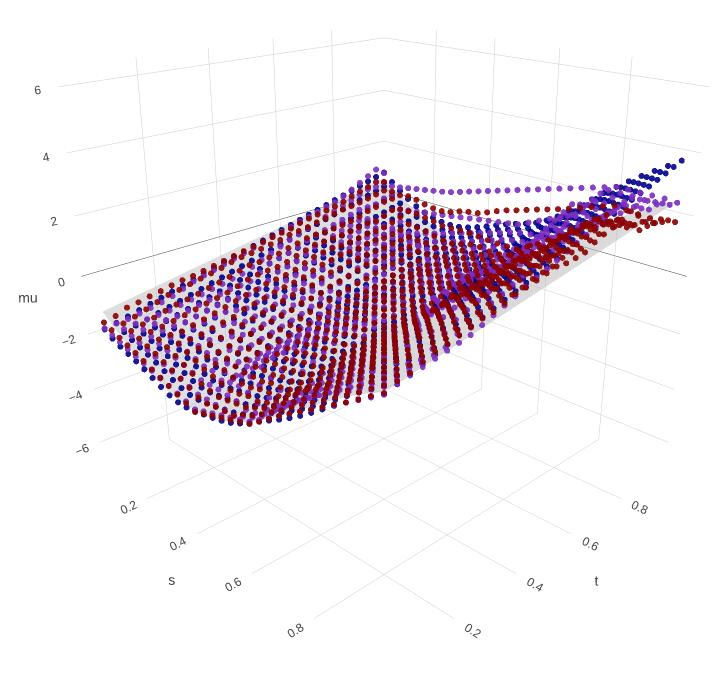} 
    \hspace{2cm}
    \includegraphics[width=.4\textwidth]{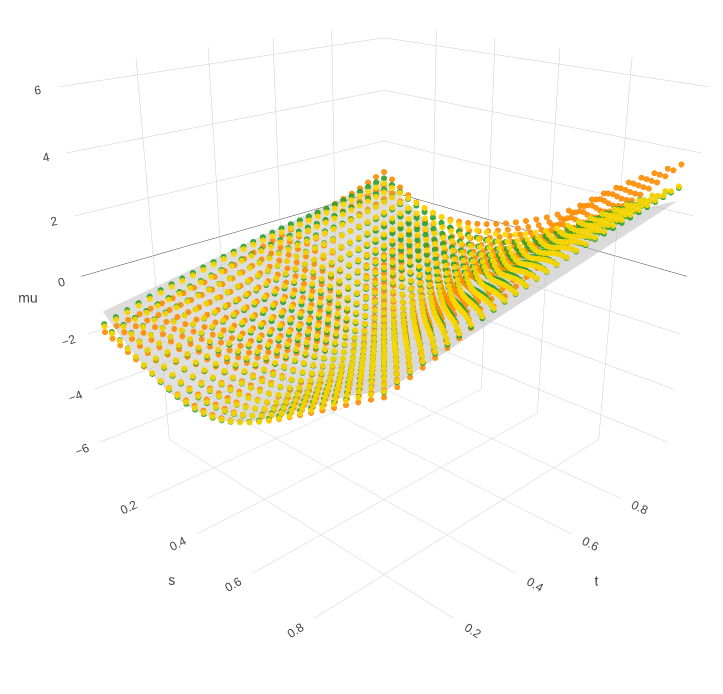} 
    \includegraphics[width=0.4\textwidth]{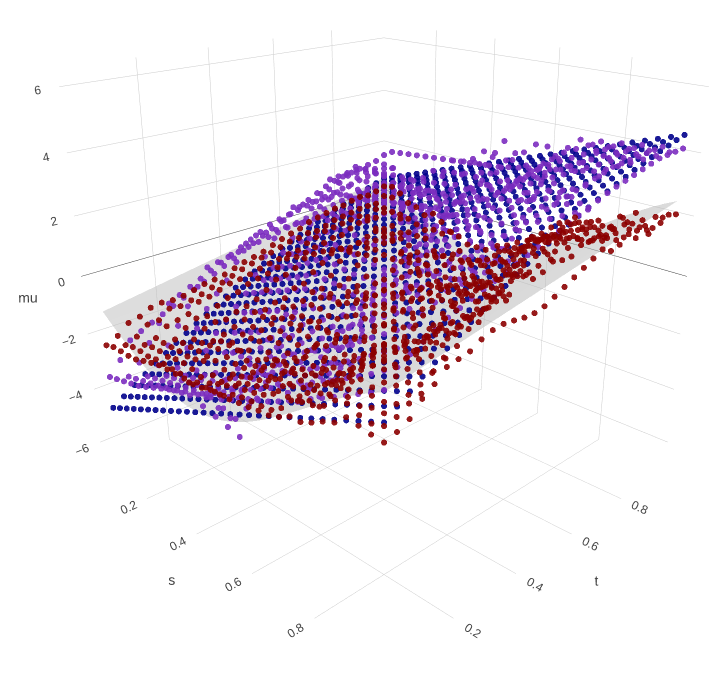} 
    \hspace{2cm}
    \includegraphics[width=.4\textwidth]{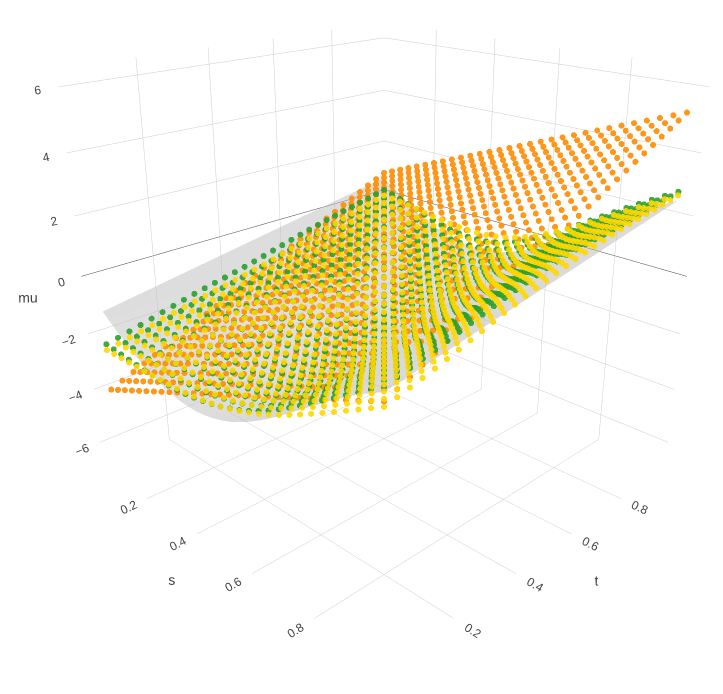} 
    %\caption{Estimation for two-dimensional functional data from Section~\ref{sec:5}. A single simulation run with $n = 50$, $p = 50$ in \textbf{Model$^2$~1} (top) and \textbf{Model$^2$~2} (bottom). The true location function $\mu$ (yellow surface) with its estimates using the naive methods \textsf{N$_1$} (blue) and \textsf{N$_2$} (orange) in left-hand panels, and the thin-plate spline estimates \textsf{S} (red) and \textsf{H} (green) in right-hand panels.}
    \caption{Estimation for two-dimensional functional data locations from Section~\ref{sec:5}. A single simulation run with $n = 50$, $p = 50$ and data generated from \textbf{Model$^2$~1} (top), \textbf{Model$^2$~2} (central) and \textbf{Model$^2$~3} (bottom). The true location function $\mu$ (gray surface) with its estimates using the competitor methods \textsf{NP} (blue), \textsf{N$_1$} (violet) and \textsf{N$_2$} (dark red) in the left-hand panels, and the thin-plate spline estimates \textsf{S} (orange), \textsf{H} (green) and \textsf{A} (yellow) in the right-hand panels.}
    \label{fig:d2}
\end{figure}

Examples of the true regression function $\mu_0$ using the non-robust least squares fit (\textsf{S}) and the robust Huber fit (\textsf{H}) or Absolute fit (\textsf{A}) for functional data generated under \textbf{Models$^2$~1--3} are displayed in Figure~\ref{fig:d2}. We see that especially for the heavier-tailed \textbf{Models}$^2$~\textbf{2}--\textbf{3}, the non-spline estimators in the left-hand panel, as well as the non-robust spline estimator \textsf{S}, fail to accurately capture the shape of the true surface. In contrast, the robust spline estimators \textsf{H} and \textsf{A} recover the true shape of $\mu_0$ without evident systematic departures.

The computation times of the considered methods can be found in Tables~\ref{Tab:times_model_d1} and~\ref{Tab:times_d2_regular_cont}. The proposed spline estimators are more expensive than the naive and competing non-robust methods, reflecting the additional cost of robust estimation. Among the spline methods, \textsf{S} and \textsf{L} are the cheapest, while the cost of \textsf{H} grows more slowly with $n$ and $p$ and becomes comparable to that of \textsf{S} for the largest balanced settings. The $L_1$ estimator \textsf{A} is the most computationally demanding. Both \textsf{H} and \textsf{A} are solved via quadratic programming (OSQP); in particular, \textsf{A} solves the $L_1$ problem directly, without the smoothing approximation required by an IRLS implementation.

% In dimension $d=2$, the computation times are considerably higher than in the one-dimensional case, reflecting the larger optimization problems, but exhibit the same qualitative ordering: the proposed spline estimators are more expensive than the naive and competitor methods, with \textsf{S} and \textsf{L} being the cheapest and \textsf{A} the most computationally demanding.

\begin{table}[htpb]
\centering
\caption{One-dimensional case, average computation times (in seconds), averaged over the one-dimensional \textbf{Models$^1$~1--3}; 1000 runs of the simulation study.}
\label{Tab:times_model_d1}
\resizebox{\textwidth}{!}{\begin{tabular}{cc|ccc|ccc|ccc}
\hline
 & & \multicolumn{3}{c}{$n=25$} & \multicolumn{3}{c}{$n=50$} & \multicolumn{3}{c}{$n=100$} \\
 & & $p=10$  & $p=25$  & $p=50$  & $p=10$  & $p=25$  & $p=50$  & $p=10$  & $p=25$  & $p=50$ \\
\hline
\parbox[t]{2mm}{\multirow{8}{*}{\rotatebox[origin=c]{90}{\textbf{Balanced}}}} & \textsf{S} & 0.058 & 0.489 & 3.823 & 0.146 & 1.856 & 14.966 & 0.525 & 7.797 & 79.598 \\
 & \textsf{A} & 4.552 & 24.701 & 99.373 & 10.391 & 63.898 & 408.323 & 31.690 & 224.264 & 1652.084 \\
 & \textsf{H} & 1.248 & 2.840 & 9.038 & 1.756 & 5.296 & 30.042 & 2.993 & 12.100 & 101.973 \\
 & \textsf{L} & 0.067 & 0.542 & 4.099 & 0.165 & 1.950 & 15.445 & 0.574 & 8.026 & 80.997 \\
 & \textsf{P} & 0.002 & 0.002 & 0.003 & 0.003 & 0.004 & 0.006 & 0.004 & 0.007 & 0.013 \\
 & \textsf{NP} & 0.337 & 0.733 & 2.341 & 0.560 & 2.334 & 8.487 & 1.573 & 8.355 & 34.927 \\
 & \textsf{N$_1$} & 0.000 & 0.000 & 0.000 & 0.000 & 0.000 & 0.000 & 0.000 & 0.000 & 0.001 \\
 & \textsf{N$_2$} & 0.003 & 0.004 & 0.004 & 0.005 & 0.005 & 0.007 & 0.005 & 0.006 & 0.014 \\ \hline
\parbox[t]{2mm}{\multirow{8}{*}{\rotatebox[origin=c]{90}{\textbf{Unbalanced}}}} & \textsf{S} & 0.578 & 2.535 & 7.715 & 1.762 & 7.710 & 28.140 & 5.373 & 28.042 & 125.147 \\
 & \textsf{A} & 19.363 & 82.755 & 302.109 & 39.998 & 304.748 & 1096.093 & 187.014 & 1103.313 & 2908.577 \\
 & \textsf{H} & 5.264 & 10.636 & 24.172 & 9.176 & 25.414 & 74.490 & 20.431 & 84.048 & 216.913 \\
 & \textsf{L} & 0.839 & 3.201 & 8.760 & 2.320 & 8.750 & 29.847 & 6.300 & 29.747 & 124.964 \\
 & \textsf{P} & 0.005 & 0.009 & 0.014 & 0.008 & 0.015 & 0.027 & 0.014 & 0.030 & 0.056 \\
 & \textsf{NP} & 0.282 & 0.723 & 2.222 & 0.546 & 2.241 & 8.162 & 1.610 & 8.153 & 33.319 \\
 & \textsf{N$_1$} & 0.000 & 0.000 & 0.000 & 0.000 & 0.000 & 0.000 & 0.000 & 0.000 & 0.000 \\
 & \textsf{N$_2$} & 0.004 & 0.004 & 0.004 & 0.005 & 0.005 & 0.006 & 0.005 & 0.007 & 0.008 \\ \hline
\end{tabular}}
\end{table}

% latex table generated in R 4.5.2 by xtable 1.8-4 package
% Sun Oct 11 17:18:15 2026
\begin{table}[htpb]
\centering
\caption{Average computation times (in seconds), averaged over the bivariate \textbf{Models$^2$~1--3}. 500 runs of the simulation study.} 
\label{Tab:times_d2_regular_cont}
\resizebox{\textwidth}{!}{\begin{tabular}{c|cccc|cccc|cccc}
   \hline
 & \multicolumn{4}{c}{$n=10$} & \multicolumn{4}{c}{$n=25$} & \multicolumn{4}{c}{$n=50$} \\
 & $p=10$  & $p=50$  & $p=100$  & $p=200$  & $p=10$  & $p=50$  & $p=100$  & $p=200$  & $p=10$  & $p=50$  & $p=100$  & $p=200$ \\
   \hline
  \textsf{S}     &     0.040 &     0.727 &     5.594 &    45.978 &     0.066 &     3.753 &    29.480 &   389.026 &     0.153 &    14.559 &   128.409 &  2406.099 \\ 
  \textsf{A}     &     1.672 &    19.508 &   163.364 &  1535.588 &     4.262 &    95.181 &  1001.460 &  5071.293 &     9.963 &   349.739 &  2847.272 & 10740.257 \\ 
  \textsf{H}     &     1.069 &     4.088 &    17.388 &   151.375 &     1.351 &     9.253 &    70.793 &   511.455 &     1.912 &    26.695 &   253.074 &  1411.141 \\ 
  \textsf{L}     &     0.043 &     0.854 &     6.394 &    52.418 &     0.075 &     3.997 &    30.547 &   363.862 &     0.171 &    15.082 &   122.494 &  2259.472 \\ 
  \textsf{NP}    &     0.390 &     3.786 &    14.276 &    54.197 &     1.158 &    22.345 &    84.841 &   338.874 &     3.992 &    93.085 &   368.974 &  1495.821 \\ 
  \textsf{N$_1$} &     0.000 &     0.000 &     0.000 &     0.000 &     0.000 &     0.000 &     0.000 &     0.000 &     0.000 &     0.000 &     0.001 &     0.001 \\ 
  \textsf{N$_2$} &     0.033 &     0.035 &     0.037 &     0.044 &     0.028 &     0.031 &     0.033 &     0.045 &     0.037 &     0.043 &     0.053 &     0.115 \\ 
   \hline
\end{tabular}}
\end{table}

\section{Supplement to Section \ref{sec:application}: Robustness to Outlying Data}
\label{sec:appendix_C}

\begin{figure}
    \centering
    \includegraphics[width=0.95\textwidth]{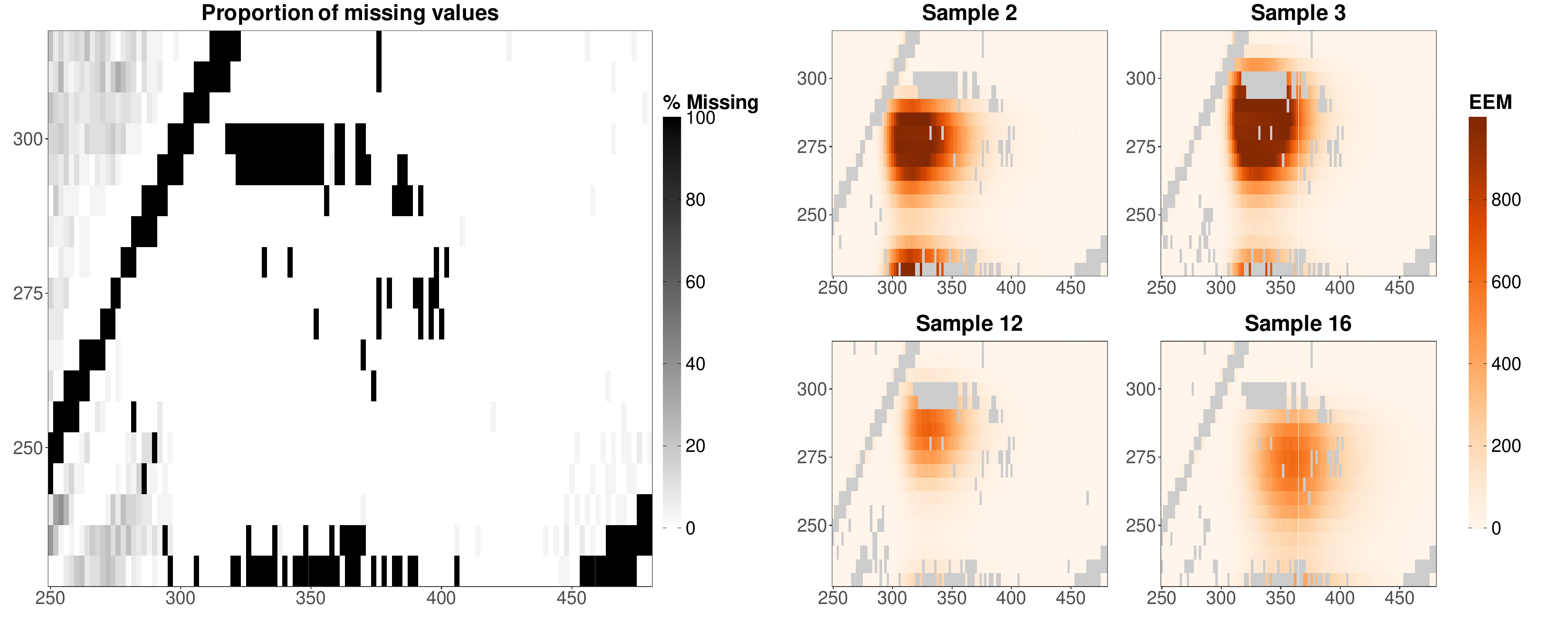}
    \caption{Left: The proportion of missing values for each excitation-emission wavelength pair (darker for larger missing proportion). Right: Four representative EEM samples, with missing values indicated by gray tiles.}
    \label{fig:Dorrit_NA_and_Examples}
\end{figure}

Consider the Dorrit fluorescence dataset \citep{Riu:2003}, consisting of 27 Excitation-Emission Matrices (EEMs) of samples containing mixtures of four fluorophores: hydroquinone, tryptophan, phenylalanine, and dopa. Emission wavelengths range from 250~nm to 482~nm at 2~nm intervals, and excitation wavelengths range from 230~nm to 315~nm at 5~nm intervals, yielding a total of $(116 \times 18)$ emission-excitation pairs. Each EEM can be viewed as a two-dimensional fluorescence intensity surface, representing an observation of functional data over a two-dimensional domain. Some measurements are missing across samples and locations due to the removal of regions affected by scattering during data preprocessing. Figure~\ref{fig:Dorrit_NA_and_Examples} (left) displays the proportion of missing values in each location. As a result of the preprocessing procedure, some locations have missing values across all samples.

This dataset has been extensively studied in the robust statistics literature \citep{Huber:2024, Hirari:2026}, as it contains different types of outliers. Some samples exhibit unusually high concentrations of specific fluorophores, resulting in censored measurements at the highest fluorescence intensities (samples 2, 3, and 5, classified by \citealp{Hirari:2026} as ``rowwise outliers''). Further anomalous measurements also occur at specific excitation-emission combinations, typically close to the locations associated with missing data (``cellwise outliers'' in the terminology of \citealp{Hirari:2026}). Those effects are restricted to single or a few nearby locations and may differ across samples. Figure~\ref{fig:Dorrit_NA_and_Examples} (right panel) shows four representative EEM samples. Samples 2 and 3 exhibit markedly higher fluorescence intensities than the remaining samples.

\begin{sidewaysfigure}
    \centering
    \includegraphics[width=0.85\textheight]{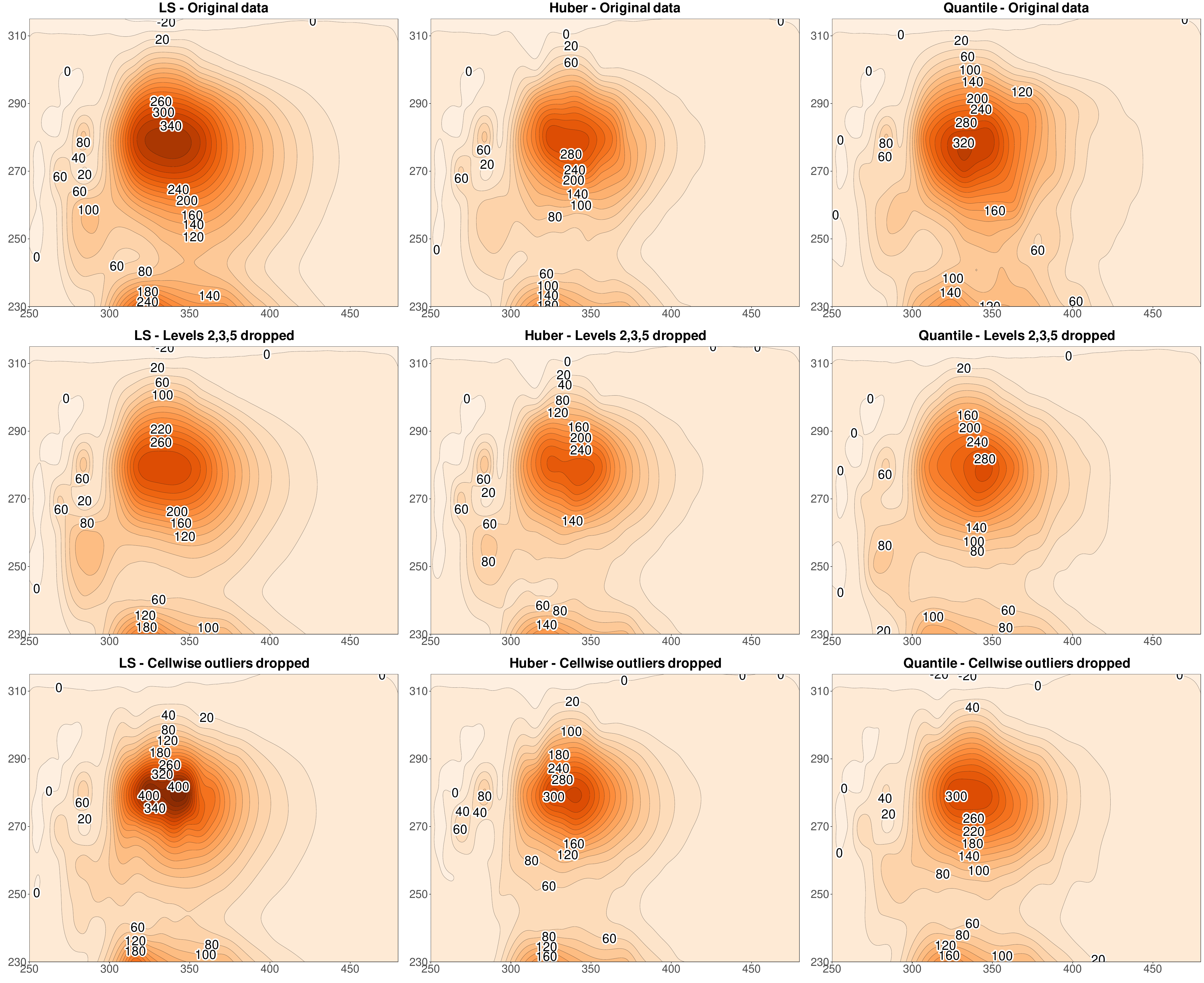}%
    \hspace{0.5cm}
    \raisebox{0.5\height}{%
        \includegraphics[width=0.055\textheight]{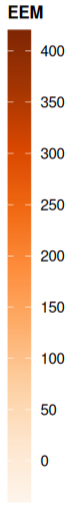}%
    }
    \caption{Contour plots of the location estimators for the EEM from the Dorrit data, obtained using the Least Squares (LS), Huber and Quantile losses. The top row shows the estimators for the original data, the central row are the estimators for the data without samples 2, 3 and 5 and the bottom row the estimators for the data without any outlier.}
    \label{fig:Dorrit_estimators}
\end{sidewaysfigure}

\begin{sidewaysfigure}
    \centering
    \includegraphics[width=0.85\textheight]{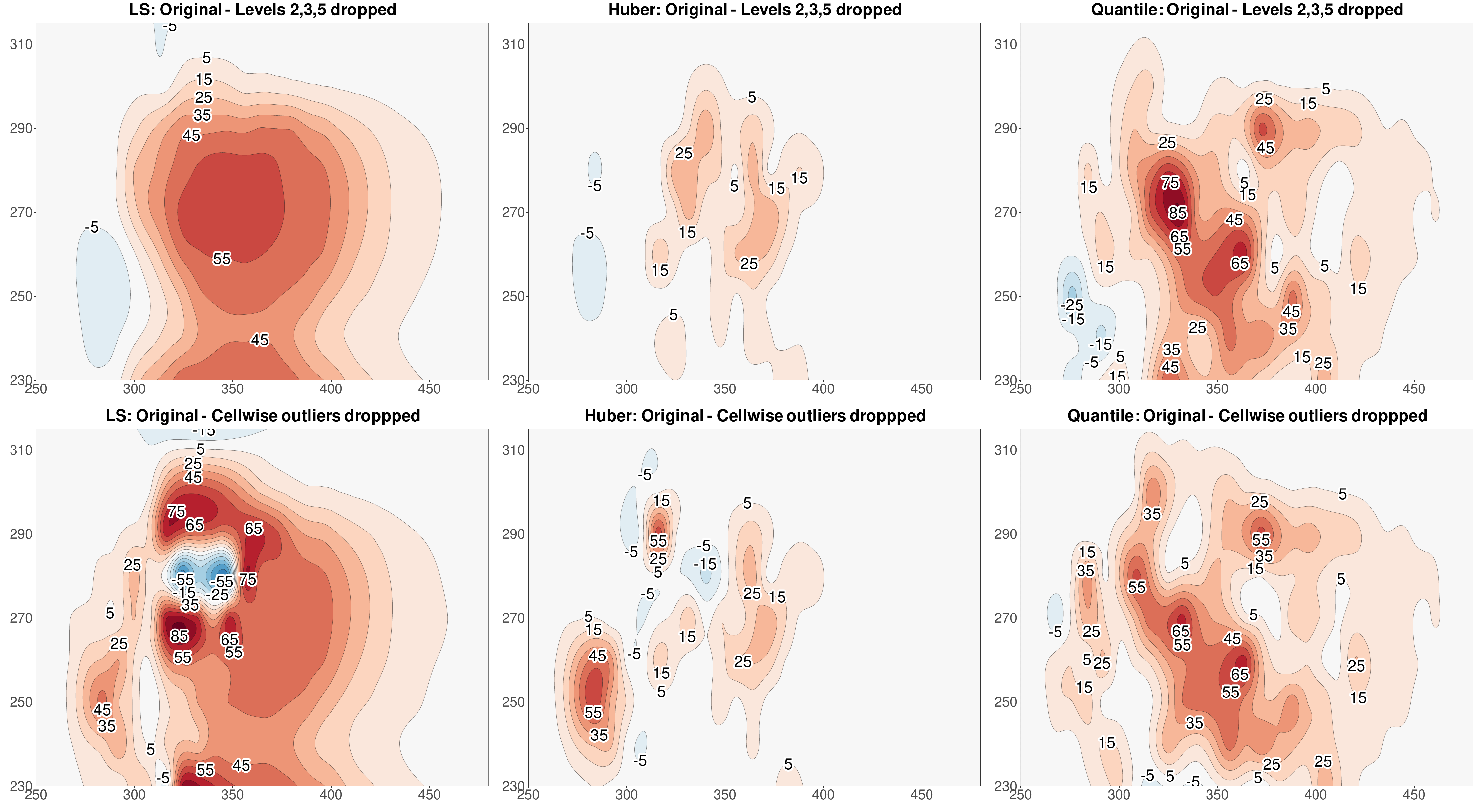}%
    \hspace{0.5cm}
    \raisebox{0.20\height}{%
        \includegraphics[width=0.065\textheight]{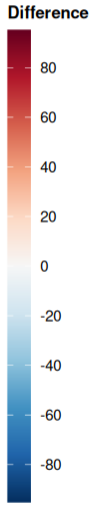}%
    }
    \caption{Contour plots of the differences among location estimators estimated using the original EEM data and the ones obtained on the outlier-free datasets, obtained using the Least Squares (LS), Huber, and Quantile losses. The top row shows the estimators for the original data, the central row are the estimators for the data without samples 2, 3 and 5 and the bottom row the estimators for the data without any outlier.}
    \label{fig:differeces_Dorrit}
\end{sidewaysfigure}

Our analysis is organized into three steps. First, we compute the spline location estimates using the Least Squares, Huber, and Absolute (Quantile) loss (we restrict the analysis to a $(29 \times 9)$ subset of the available emission and excitation wavelengths to reduce the computational burden). The resulting estimators are displayed in the first row of Figure~\ref{fig:Dorrit_estimators}.

Second, we construct two outlier-free datasets using different strategies. First, we remove three samples, 2, 3, and 5, identified as extreme in the existing chemometrics literature \citep{Huber:2024,Hirari:2026}. The right panel of Figure~\ref{fig:Dorrit_NA_and_Examples} illustrates the extreme samples 2 and 3 alongside two non-extreme samples 12 and 16. Alternatively, we construct a second dataset by identifying outliers based on the residuals of the EEMs with respect to the quantile-loss spline estimator. Specifically, at each of the $(29 \times 9)$ emission-excitation wavelength pairs, we flag observations whose residuals fall beyond the whisker fences of a boxplot computed across the 27 functional observations. 

Finally, we compute the location estimators on both outlier-free datasets. 
Figure~\ref{fig:Dorrit_estimators} displays the contour plots of the three spline estimators for the original dataset, as well as after removing the two types of outliers.  Figure~\ref{fig:differeces_Dorrit} reports the differences between the estimators from the original dataset and those computed on the two outlier-free datasets.

The Least Squares estimator is the most affected by the removal of outliers. Although the overall shape of the estimated surfaces remains stable, noticeable differences in their magnitude can be observed. This is confirmed by the difference plots in Figure~\ref{fig:differeces_Dorrit}, where the Least Squares estimator exhibits the largest deviations. Comparing Huber and Absolute loss estimators, the former appears more robust to outlier removal in this setting, in agreement with our findings from Section~\ref{sec:5}.

\paragraph*{Acknowledgment} The work of S.~Nagy and M.~Cavazzutti was supported by the ERC CZ grant LL2407 of the Ministry of Education, Youth and Sport of the Czech Republic.

\bibliographystyle{apalike}
\bibliography{papers}
\end{document}